\documentclass[3p]{elsarticle}
\usepackage{graphicx, algorithmic, algorithm2e,  psfrag, amssymb, amsfonts, amsmath, pstricks-add, bm, float, comment, color, hyperref, placeins}
\usepackage{subcaption}
\usepackage{caption}
\usepackage{mathtools}
\usepackage{amsmath}

\usepackage{tikz}
\usetikzlibrary{matrix} 
\usetikzlibrary{arrows} 
\usetikzlibrary{calc} 
\usetikzlibrary{shapes}
\usetikzlibrary{fit}
\usetikzlibrary{3d}
\usetikzlibrary{positioning}
\usetikzlibrary{shadows}
\usepackage{amsthm}

\tikzstyle{block} = [draw,rectangle,thick,minimum height=2em,minimum width=2em]
\tikzstyle{sum} = [draw,circle,inner sep=0mm,minimum size=2mm]
\tikzstyle{connector} = [->,thick]
\tikzstyle{line} = [thick]
\tikzstyle{branch} = [circle,inner sep=0pt,minimum size=1mm,fill=black,draw=black]
\tikzstyle{guide} = []
\tikzset{>=latex}
\usepackage{pgfplots}
\pgfplotsset{compat=newest}

\journal{~~}

\newcommand{\Cm}{\mathbf{C}}
\newcommand{\dm}{\mathbf{d}}
\newcommand{\ff}{\mathbf{f}}

\newcommand{\gm}{\mathbf{g}}

\newcommand{\hb}{\mathbf{h}}

\newcommand{\Km}{\mathbf{K}}

\newcommand{\ngg}{n_{\mathbf{g}}}
\newcommand{\np}{n_{\mathbf{\ppm}}}
\newcommand{\ny}{n_{\mathbf{\Ym}}}

\newcommand{\ppm}{{\pmb\theta}}
\newcommand{\uu}{\mathbf{u}}

\newcommand{\Ym}{\pmb{\xi}}

\newcommand{\ie}{\textit{i.e.},}
\newcommand{\eg}{\textit{e.g.},}

\newcommand{\Exp}{\mathbb{E}}
\newcommand{\Var}{\mathrm{Var}}

\newcommand{\rhoo}{\boldsymbol{\rho}}
\newcommand{\kappaa}{\pmb\kappa}
\newcommand{\thetaa}{\pmb\theta}
\DeclareMathOperator*{\argmin}{argmin} 

\begin{document}

\begin{frontmatter}

\title{Bi-fidelity Stochastic Gradient Descent for Structural Optimization under Uncertainty}

\author{Subhayan De}
\ead{Subhayan.De@colorado.edu}
\author{Kurt Maute}
\ead{maute@colorado.edu}
\author{Alireza Doostan\corref{mycorrespondingauthor}}
\address{Smead Aerospace Engineering Sciences Department, University of Colorado, Boulder, CO 80309, USA}


\cortext[mycorrespondingauthor]{Corresponding author}
\ead{Alireza.Doostan@colorado.edu}


\begin{abstract}
The presence of uncertainty in material properties and geometry of a structure is ubiquitous. The design of robust engineering structures, therefore, needs to incorporate uncertainty in the optimization process. 
Stochastic gradient descent (SGD) method can alleviate the cost of optimization under uncertainty, which includes statistical moments of quantities of interest in the objective and constraints. However, the design may change considerably during the initial iterations of the optimization process which impedes the convergence of the traditional SGD method and its variants. 
In this paper, we present two SGD based algorithms, where the computational cost is reduced by employing a \textit{low-fidelity} model in the optimization process. In the first algorithm, most of the stochastic gradient calculations are performed on the low-fidelity model and only a handful of gradients from the \textit{high-fidelity} model are used per iteration, resulting in an improved convergence. 
In the second algorithm, we use gradients from \textit{low-fidelity} models to be used as \textit{control variate}, a variance reduction technique, to reduce the variance in the search direction. These two bi-fidelity algorithms are illustrated first with a conceptual example. Then, the convergence of the proposed bi-fidelity algorithms is studied with two numerical examples of shape and topology optimization and compared to popular variants of the SGD method that do not use \textit{low-fidelity} models. The results show that the proposed use of a bi-fidelity approach for the SGD method can improve the convergence. Two analytical proofs are also provided that show the linear convergence of these two algorithms under appropriate assumptions. 

\end{abstract}


\end{frontmatter}



%
%
\section{Introduction}\label{sec:intro}

In simulation-based engineering, models, often in the form of discretized (partial) differential equations, are used for purposes such as analysis, design space exploration, uncertainty quantification, and design optimization. In the context of structural optimization, such as shape and topology optimization, the models need to be simulated many times throughout the optimization process \cite{spillers2009structural}. Structures are often subjected to uncertainties in the material properties, geometry, and external loads \cite{hasselman2001quantification,uncertainty2006,bulleit2008uncertainty}. Hence, for robust design of these structures, such uncertainties must be accounted for in the optimization process.
The most commonly used method to compute the stochastic moments of the design criteria for Optimization under Uncertainty (OuU) is random sampling based Monte Carlo approach. 
In this approach, at every iteration, statistics calculated from a number of random samples are used as the objective and constraints for the optimization. However, the number of random samples often needs to be large to get a small approximation error.
As a result, this approach increases the computational burden even further as one needs to solve the governing equations many times at every iteration of the optimization \cite{diwekar1997efficient,sahinidis2004optimization,diwekar2008optimization,de2017efficient}. Stochastic collocation \cite{babuvska2007stochastic,nobile2008sparse} or polynomial chaos expansion \cite{roger2003stochastic,xiu2002wiener} methods can also be utilized to estimate these statistics, but the required number of random samples increases rapidly with the number of optimization parameters. Note that, sparse polynomial chaos expansions \citep{Doostan09b,doostan2011non,blatman2010adaptive,hampton2016compressive,hampton2018basis} can be used to reduce the computational cost. However, for design optimization problems, where uncertainty is represented by a large number of random variables the computational cost may remain unbearable. 

In deterministic optimization problems, $\ppm\in \mathbb{R}^{\np}$ denote the vector of design parameters and the objective function $f(\thetaa): \mathbb{R}^{\np} \rightarrow \mathbb{R}$, \textit{e.g.}, strain energy of a structure, depends on $\thetaa$. For constrained optimization problems, let $\gm(\ppm): \mathbb{R}^{\np} \rightarrow \mathbb{R}^{\ngg}$ be $\ngg$ real-valued constraint functions, \textit{e.g.}, allowable mass of a structure. 
The constraints are satisfied if $g_i(\ppm) \le 0$ for $i=1,\dots,\ngg$. 
Hence, the optimization problem can be written as 
\begin{equation}\label{eq:det_opt}
    \mathop{\min~}\limits_{\ppm}f(\ppm) \mbox{   subject to } g_i(\ppm) \le {0} \text{    for }i=1,\dots,\ngg.
\end{equation}

In the presence of uncertainty, a reformulation of the optimization problem in \eqref{eq:det_opt} is generally used. 
Let $\Ym\in\mathbb{R}^{\ny}$ be the vector of random variables, with known probability distribution function, characterizing the system uncertainty. The objective function $f(\ppm;\Ym): \mathbb{R}^{\np} \times \mathbb{R}^{\ny} \rightarrow \mathbb{R}$ now also depends on the realized values of $\Ym$. Similarly, for constrained optimization problems, $\ngg$ real-valued constraint functions $\gm(\ppm;\Ym): \mathbb{R}^{\np} \times \mathbb{R}^{\ny} \rightarrow \mathbb{R}^{\ngg}$ in general depend on $\Ym$. The optimization problem is defined using the \textit{expected risk} $R(\ppm) = \Exp[f(\ppm;\Ym)]$ and expected constraint value $C_i(\ppm) = \Exp[g_i(\ppm;\Ym)]$ as follows \cite{sandgren2002robust,zang2005review,calafiore2008optimization,de2017efficient}
\begin{equation}\label{eq:stoch_opt}
    \mathop{\min~}\limits_{\ppm}R(\ppm) \mbox{   subject to } C_i(\ppm) \le {0} \text{    for }i=1,\dots,\ngg,
\end{equation}
where $\Exp[\cdot]$ denotes the mathematical expectation of its argument.

The stochastic gradient descent methods illustrated in this paper utilize evaluations of the gradients of $R(\ppm)$ and $\Cm(\ppm)$. 
%
In optimization, a combination of these gradients defines a direction along which a search is performed. Each of the investigated methods modifies these directions with the goal to improve stability in the optimization so that a converged solution is reached more reliably in a smaller number of optimization steps. 

\subsection{Multi-fidelity Models with Applications to OuU}
In most engineering problems (or applications), multiple models are often available to describe the system. Some of these models are able to describe the behavior with a higher level of accuracy but are generally associated with high computational cost; in the sequel referred to as \textit{high-fidelity models}. Models with lower computational cost, on the other hand, are often (not always) less accurate and termed \textit{low-fidelity models}. For a structural system analyzed by the finite element method, \textit{low-fidelity models} can be obtained, for instance, using coarser grid discretizations of the governing equations. Multi-fidelity methods exploit the availability of these different models to accelerate design optimization, parametric studies, and uncertainty quantification; see, \textit{e.g.}, Fern{\'a}ndez-Godino \textit{et al}. \cite{fernandez2016review}, Peherstorfer \textit{et al}. \cite{peherstorfer2018survey}, and the references therein. In these methods, most of the computation is performed using the \textit{low-fidelity models} and the \textit{high-fidelity models} are utilized to {\it correct} the low-fidelity predictions. 


The use of multi-fidelity models is especially helpful in optimization, where one needs to solve the governing equations multiple times \cite{park2017remarks}. Booker \textit{et al.}~\cite{booker1999rigorous} used pattern search to construct low-fidelity surrogate models. Forrester \textit{et al.}~\cite{forrester2007multi} used \textit{co-kriging} \cite{myers1982matrix,kennedy2001bayesian,qian2008bayesian} for constructing surrogates that are then used for design optimization of an aircraft wing. 
Similar applications of multi-fidelity models for optimization of aerodynamic design are performed in Keane \cite{keane2003wing}, Huang \textit{et al}. \cite{huang2006sequential}, Choi \textit{et al.}~\cite{choi2008multifidelity}, and Robinson \textit{et al.}~\cite{robinson2008surrogate}. Reduced order models constructed using Pad\'{e} approximation, Shanks transformation, Krylov subspace methods, derivatives of eigenmodes with respect to design variables, and proper orthogonal decomposition are used in Chen \textit{et al}.~\cite{chen2000static}, Kirsch \cite{kirsch2000combined}, Hurtado \cite{hurtado2002reanalysis}, Sandbridge and Haftka \cite{sandridge1989accuracy}, and Weickum \textit{et al}.~\cite{weickum2006multi} to perform optimization of structural systems. 
Eldred and Dunlavy \cite{eldred2006formulations} compared different data-fit, multi-fidelity, and reduced order models to construct surrogates for optimization. 
Yamazaki \textit{et al}. \cite{yamazaki2010design} used a gradient enhanced \textit{kriging} to perform design optimization as well as uncertainty quantification at the optimal design point. 
Another optimization method known as \textit{space mapping} \cite{bandler1994space,bakr2000review,bakr2001introduction,koziel2016rapid} uses the low-fidelity models to solve an approximate optimization problem, where the input parameter space is mapped onto a different space to construct multi-fidelity models. Fischer \textit{et al}. \cite{fischer2017bayesian} used a Bayesian approach for estimating weights of the low-fidelity models to be used for optimization in a  multi-fidelity optimization setting. 

Recently, the multi-fidelity approach has also been applied for efficient uncertainty quantification. There, the low-fidelity models are used to reduce the computational cost associated with the solution of governing equations of large-scale physical systems in the presence of high-dimensional uncertainty \cite{koutsourelakis2009accurate,ng2012multifidelity,narayan2014stochastic,perdikaris2015multi,peherstorfer2016multifidelity,doostan2016bi,parussini2017multi,hampton2018practical,fairbanks2018bi,peherstorfer2018survey,skinner2019reduced}. 

For OuU, Jin \textit{et al.}~\cite{jin2003use} used different metamodeling approaches, \textit{e.g.}, polynomial regression, krigging \cite{martin2005use} to ease the computational burden of the optimization iterations. 
Kroo \textit{et al}. \cite{kroo2010multifidelity} employed multi-fidelity models for multiobjective optimization, where these objectives provide a balance between performance and risk associated with the uncertainty in the problem.
Keane \cite{keane2012cokriging} used \textit{co-kriging} for robust design optimization with an application to shape optimization of a gas-turbine blade.
Allaire \textit{et al}. \cite{allaire2010bayesian} used a Bayesian approach for risk based multidisciplinary optimization, where models of multiple fidelity are used to merge information on uncertainty. In multidisciplinary optimization, Christensen \cite{christensen2012multifidelity} increased the fidelity of the model for a particular discipline from which the contribution to the uncertainty of the quantity of interest is large.
Eldred and Elman \cite{eldred2011design} and Padr\'{o}n \textit{et al}. \cite{padron2016multi} used stochastic expansion methods, \eg~stochastic collocation and polynomial chaos expansion methods to construct high- and low-fidelity surrogate models for design OuU.
March and Wilcox \cite{march2011gradient,march2012constrained,march2012provably} proposed a multi-fidelity based trust-region algorithm for optimization using gradients from low-fidelity models only. 
Ng and Wilcox \cite{ng2014multifidelity} used control variate approach for the design of an aircraft wing under uncertainty.

\subsection{Stochastic Gradient Descent Methods}
 Gradient descent methods \cite{boyd2004convex,nocedal2006numerical} for solving \eqref{eq:det_opt} are the preferred choice if $f(\ppm)$ and $\gm(\ppm)$ are differentiable with respect to $\ppm$. 
To solve the optimization problem under uncertainty in \eqref{eq:stoch_opt},
a stochastic version of gradient descent \cite{robbins1951stochastic} that has a smaller per iteration computational cost compared to a Monte Carlo approach with large number of random samples can be used. Recently, the stochastic gradient descent (SGD) method has seen increasing use in training of neural networks \cite{bottou2010SGDML}, where the number of optimization variables is large. In this method, the gradients used at every iteration are obtained using either only one or a small number of random samples of $\Ym$. However, the standard SGD method converges slowly \cite{nemirovski2009robust}. Hence, to improve the convergence various modifications to the standard SGD method have been proposed recently. Among these, Adaptive Gradient (AdaGrad)~\citep{duchi2011adagrad} and  Adaptive Moment (Adam) \citep{kingma2014adam} retard the movement in {directions} with historically large gradient magnitudes and are useful for problems that lack convexity. Adadelta~\citep{zeiler2012adadelta} removes the need for an explicitly specified learning rate, \ie~step size; however, small initial gradients affect this algorithm adversely \cite{de2019topology}. 
Another variant of SGD method, the Stochastic Average Gradient (SAG) algorithm \cite{rouxetalSAG2012} updates a single gradient using one random sample of the uncertain parameters per iteration and keeps the rest the same as in the last iteration. Then, the optimization parameters are updated using an average of all computed gradients. However, for design optimization, this approach of using past gradients may lead to poor convergence \cite{de2019topology}.
The Stochastic Variance Reduced Gradient (SVRG) algorithm \cite{johnson2013accelerating} uses a \textit{control variate} \cite{ross2013simulation} to reduce the variance of the stochastic gradients. This approach, as introduced in Section \ref{sec:sgd}, can also suffer from poor convergence in OuU if the same control variate is used for a large number of iterations \cite{defazio2018ineffectiveness,de2019topology}. 

In this paper, to achieve a better convergence and to ease the computational burden of OuU, 
we formulate two bi-fidelity based variants of the SAG and SVRG algorithms to overcome their respective shortcomings. 
We use the word \textit{bi-fidelity} instead of \textit{multi-fidelity} as we are only using one high- and one low-fidelity models. 
The first algorithm, Bi-fidelity Stochastic Average Gradient (BF-SAG), evaluates most of the gradients using a low-fidelity model and then applies a gradient descent step using an average gradient. In the second algorithm, Bi-fidelity Stochastic Variance Reduced Gradient (BF-SVRG), we propose a control variate approach, where the mean of the control variate is estimated using only low-fidelity model evaluations resulting in the reduction of the computational cost. 
In addition, we use the correlation between the high- and low-fidelity gradients to reduce the variance in the estimated gradients. Further, we prove linear convergence of these two proposed algorithms for strongly convex objectives and gradients that are Lipschitz continuous.

We illustrate the proposed bi-fidelity algorithms using three numerical examples. For the first example, we choose a simple fourth order polynomial as the high-fidelity model and a second-order approximation of it as the low-fidelity model. This example  shows that the use of a low-fidelity model allows for evaluating more gradients per iteration using the same computational budget and thereby leads to a faster convergence. In the second example, we use an example of a square plate with a hole subjected to uniaxial tension. 
We optimize the shape of the hole to minimize the maximum principal stress in the plate. 
For the third example, we apply the proposed algorithms to a topology optimization problem using the \textit{Solid Isotropic Material with Penalisation} (SIMP) method \cite{bendsoe1989optimal,bendose2003topology,sigmund2013topology}.  In all examples, we observe convergence improvements of the proposed bi-fidelity strategies over their standard counterparts, \textit{e.g.}, SAG and SVRG. 

The rest of the paper is organized as follows. In the next section, we briefly discuss some popular variants of the SGD method. In Section \ref{sec:method}, we introduce the proposed bi-fidelity based optimization algorithms with convergence properties and computational cost analyses. In Section \ref{sec:examples}, we illustrate various aspects of the algorithms with three numerical examples. Finally, we conclude our paper with a brief discussion of future research directions of bi-fidelity based SGD methods. 

\RestyleAlgo{boxruled}

\section{Background}

In this section, 
we discuss the stochastic gradient method and its variants, namely, stochastic average gradient (SAG) and stochastic gradient reduced gradient decent (SVRG). A related concept of variance reduction using a control variate is also briefly discussed. 

\subsection{Stochastic Gradient Descent (SGD) Method and its Variants} 
\label{sec:sgd}
With SGD, for unconstrained problems the expected risk is minimized as mentioned in Section \ref{sec:intro}. For constrained optimization problems, an unconstrained formulation of \eqref{eq:stoch_opt} can be used employing a penalty formulation and constraint violation defined as $g^+_j(\ppm;\Ym)=0$ for $g_j(\ppm;\Ym)\leq0$ and $g^+_j(\ppm;\Ym)=g_j(\ppm;\Ym)$ for $g_j(\ppm;\Ym)>0$, $j=1,\dots,\ngg$. The optimization problem is then formulated as follows
\begin{align}
\label{eq:opt_def2}
\mathop{\min~}\limits_{\ppm}J(\ppm)=R(\ppm) + \sum_{j=1}^{\ngg}{\kappa_j} C_j(\ppm),
\end{align}
where the objective $J(\ppm)$ is a combination of the expected risk $R(\ppm)$ and the expected squared constraint violations $C_j(\ppm) = \Exp\left[(g_j^+(\ppm;\Ym))^2\right]$ for $j=1,\dots,\ngg$. 
{In (\ref{eq:opt_def2}),} $\kappaa$ is a user-specified penalty parameter vector. 
For large values of $\kappaa$, (\ref{eq:opt_def2}) has similar solutions to (\ref{eq:det_opt}). 

The objective in \eqref{eq:opt_def2} is often estimated using a set of realizations $\{\Ym_i\}_{i=1}^N$ of the random vector $\Ym$ in a standard Monte Carlo approach. Using these realizations the objective is written as
\begin{equation}
J_N(\ppm) = \frac{1}{N} \sum_{i=1}^N \left(f(\ppm;\Ym_i) + \sum_{j=1}^{\ngg}\kappa_j (g_j^+(\ppm;\Ym_i))^2 \right).
\end{equation}
The basic SGD method uses a single realization of $\Ym$ from its set of $N$ realizations to perform the update on $\ppm$ at the $k$th iteration utilizing the gradient $\nabla f(\ppm_k;\Ym)$ \citep{bottou2018optimization}. For a constrained optimization problem, the search direction $\hb_k$ is estimated as the combination of the gradients of the cost function $f(\ppm_k;\Ym)$ and the square of constraint violations $(\gm^+(\ppm_k;\Ym))^2$ as 
\begin{equation}\label{eq:grad}
\begin{split}
&\hb_k:=\hb(\ppm_k;\Ym_i)=\nabla{f}(\ppm_k;\Ym_i) + \sum_{j=1}^{\ngg}\kappa_j \nabla\left( g_j^+(\ppm_k;\Ym_i) \right)^2, 
\end{split}
\end{equation}
where $\Ym_i$ is selected uniformly at random from its set of $N$ realizations.
The parameter update is then applied as follows
\begin{equation}\label{eq:sgd}
\begin{split}
&\ppm_{k+1}=\ppm_k-\eta \hb_{k},\\
\end{split}
\end{equation}
where $\eta $ is the step size, also known as the \textit{learning rate}. 
These steps are illustrated in Algorithm \ref{alg:sgd}. Following \eqref{eq:grad} and \eqref{eq:sgd}, the SGD method performs only one gradient calculation per iteration; hence, its computational cost per iteration is relatively small. 
However, the descent direction $\nabla J$ is not followed at every iteration. Instead, the descent is achieved in expectation as the expectation of the stochastic gradient is the same as the gradient of the objective $J(\ppm)$. As a result, the convergence of the SGD method can be very slow \cite{bottou2010SGDML}. A straightforward extension of the SGD method is to use a small batch of random samples to compute the search direction $\hb_k$. This version is known as \textit{mini-batch gradient descent} \citep{ruder2016overview,bottou2018optimization}.
\begin{algorithm}[htb]
	\begin{algorithmic}
		\STATE Given $\eta$.
		\STATE Initialize $\ppm_1$.
		\FOR {$k=1,2,\dots$}
		\STATE Compute $\hb_k:=\hb(\ppm_k;\Ym)$. {[see Eqn. (\ref{eq:grad})]}
		\STATE Set $\ppm_{k+1} \leftarrow \ppm_{k} - \eta\hb_k$. 
		\ENDFOR
	\end{algorithmic}
	\caption{\textit{Stochastic gradient descent} \citep{bottou2018optimization}}
	\label{alg:sgd}
\end{algorithm}
In the past few years, several modifications of the SGD method have been proposed to improve its convergence \cite{bottou2018optimization}. Two of them, namely, the Stochastic Average Gradient (SAG) and the Stochastic Variance Reduced Gradient (SVRG) algorithms that are relevant to this paper are discussed next.

\subsection{\texorpdfstring{Stochastic Average Gradient (SAG) Algorithm}{Stochastic Gradient Descent: SAG}}
\label{subsec:SAG}
One popular variant of the SGD method is the SAG algorithm~\citep{rouxetalSAG2012}, which updates the gradient information for one random sample at every iteration and keeps the old gradients for other samples. 
The parameters are then updated using \eqref{eq:sgd} with the search direction $\hb_k$ defined as 
\begin{equation}\label{eq:sag}
\begin{split}
\hb_k & = \frac{1}{N}\sum_{i=1}^{N}\dm_{k,i};\\
\dm_{k,i} &= \begin{cases}
\hb(\ppm_k;\Ym_{i})\qquad \text{if $i=t\in \{1,2,\dots,N\}$};\\
\dm_{k-1,i}\qquad\qquad \text{otherwise},\\
\end{cases}
\end{split}
\end{equation}
where $t$ is selected uniformly at random from $\{1,2,\dots,N\}$ and at the start of the algorithm $\dm_{0,i}=\mathbf{0}$ for $i=1,\dots,N$. 
This use of previous gradient information accelerates the convergence of the optimization compared to the standard SGD method defined in Algorithm \ref{alg:sgd}. However, relying on past gradients can be impede convergence rate and/or stability for OuU as the design often changes drastically during the early iterations of the optimization process \cite{de2019topology}.
In this paper, we use a batch version of this algorithm, where $N_h>1$ gradients are updated at every iteration (see Algorithm \ref{alg:sag}). 
\begin{algorithm}[htb]
	\begin{algorithmic}
		\STATE Given $\eta$ and $N_h$.
		\STATE Initialize $\ppm_1$.
		\STATE Initialize $\dm_{0,i} = \bm{0}$ for $i=1,\dots,N$.
		\FOR {$k=1,2,\dots,$}
		\STATE Draw $\{t_j\}_{j=1}^{N_h}$ uniformly at random from $\{1,\cdots,N \}$.
		\FOR{$i=1,2,\dots,N$}
		\IF{$i\in \{t_j\}_{j=1}^{N_h}$}
		\STATE Compute $\dm_{k,i} := {\hb}(\ppm_{k};\Ym_{i})$
		\ELSE
		\STATE {$\dm_{k,i}:=\dm_{k-1,i}$.}
		\ENDIF
		\ENDFOR
		\STATE $\ppm_{k+1} \leftarrow \ppm_k - \frac{\eta}{N}\sum_{i=1}^{N}\dm_{k,i}$. [see Eqn. (\ref{eq:sag})]
		\ENDFOR
	\end{algorithmic}
	\caption{Batch Implementation of \textit{SAG}}
	\label{alg:sag}
\end{algorithm}

\subsection{Variance Reduction using Control Variates}
In this subsection, we briefly discuss control variates, a variance reduction technique that we use in one of our proposed algorithms. It has also been used in the SVRG algorithm, a variant of the SGD method described in Section \ref{sec:sgd}.
A \textit{control variate} can be used to estimate the expected value of a random variable $X$ via Monte Carlo averaging, while reducing the variance of the estimate \cite{ross2013simulation,hammersley2013monte}. Here, another random variable $Y$ is introduced such that it is correlated with $X$, is cheaper to simulate than $X$, and either the expected value of $Y$, $\mathbb{E}[Y]$, is known or can be estimated accurately and relatively cheaply.
Using $Y$, $\Exp[Y]$, and the standard Monte Carlo simulation, the expected value of 
\begin{equation}
\label{eqn:cv}
  Z =  X-\alpha\left(Y-\Exp[Y] \right),
\end{equation}
is estimated as an unbiased estimator of $\Exp[X]$. In (\ref{eqn:cv}), $\alpha$ is a control variate parameter and when set to 
\begin{equation}
\label{eqn:opt_cv}
    \alpha^* = \frac{\sigma_{XY}}{\sigma^2_Y},
\end{equation}
the sample average estimate of $\Exp[Z]$ achieves the minimum mean squared error in the estimate of $\Exp[X]$. In (\ref{eqn:opt_cv}), $\sigma_{XY}$ is the covariance between $X$ and $Y$ and $\sigma^2_Y$ is the variance of $Y$.
When $\alpha =\alpha^*$ is used in (\ref{eqn:cv}) to estimate $\Exp[Z]$ with $N$ Monte Carlo samples of $X$ and $Y$, the variance of sample average of $Z$ is given by $\frac{1}{N}(1-\rho_{XY}^2)\sigma^2_X$, where $\rho_{XY}=\frac{\sigma_{XY}}{\sigma_{X}\sigma_{Y}}$ is the correlation between $X$ and $Y$, and $\sigma^2_X$ is the variance of $X$. 
Since, $-1\leq \rho_{XY}\leq 1$, the reduction in variance of the estimate based on $Z$ can be seen when compared to the variance of the standard Monte Carlo estimate of $\Exp[X]$ given by $\frac{1}{N}\sigma^2_X$. In an ideal situation, $\rho_{XY}\rightarrow 1$ and as a result $\frac{1}{N}(1-\rho_{XY}^2)\sigma^2_X\rightarrow 0$. A bi-fidelity algorithm proposed in Section \ref{sec:bfsvrg} employs this variance reduction technique, where $X$ and $Y$ consist of high- and low-fidelity model gradients, respectively. 

For random vectors $\mathbf{X}$ and $\mathbf{Y}$, \eqref{eqn:cv} is replaced by 
\begin{equation}\label{eq:cv_vector}
    \mathbf{Z} = \mathbf{X} - \boldsymbol{\alpha} (\mathbf{Y}-\Exp[\mathbf{Y}]),
\end{equation}
where $\boldsymbol{\alpha}$ is a coefficient matrix. The optimal value of the coefficient matrix  $\boldsymbol{\alpha}$ that minimizes the trace of the covariance matrix of $\mathbf{Z}$ is given by
\begin{equation}\label{eq:opt_cv_matrix}
    \boldsymbol{\alpha}^* =  \mathbb{V}_{\!_\mathbf{Y}}^{-1}\mathbb{C}_{\!_\mathbf{XY}},
\end{equation}
where $\mathbb{V}_{\!_\mathbf{Y}}$ is the covariance matrix of $\mathbf{Y}$ and $\mathbb{C}_{\!_\mathbf{XY}}$ is the cross-covariance between $\mathbf{X}$ and $\mathbf{Y}$. 

\subsection{Stochastic Variance Reduced Gradient (SVRG) Algorithm}
The second variant of the SGD method that we present here is the SVRG algorithm \citep{johnson2013accelerating}. In this algorithm, { a variance reduction method is introduced by maintaining} a  parameter estimate ${\ppm}_{\mathrm{prev}}$ at every inner iteration that is updated only during the outer iteration. 
Using this parameter estimate ${\ppm}_{\mathrm{prev}}$ and $N_h$ samples of $\Ym$, the mean of $\hb(\ppm_\mathrm{prev},\Ym)$ is estimated as
\begin{equation}\label{eq:svrg_cv}
\widehat{\hb}(\ppm_\mathrm{prev}) = \frac{1}{N_h}\sum_{i=1}^{N_h} \hb(\ppm_\mathrm{prev};\Ym_i).
\end{equation}
%
Note that, in \eqref{eq:svrg_cv} $N_h$ can be smaller than $N$ in \eqref{eq:sag}.
Next, the update rule in \eqref{eq:sgd} is applied with the search direction $\hb_k$ defined as
\begin{equation}\label{eq:svrg}
\hb_k = \hb(\ppm_k;\Ym_{t})-\hb(\ppm_\mathrm{prev};\Ym_{t})+\widehat{\hb}(\ppm_\mathrm{prev})
\end{equation}
for a chosen $t$ uniformly at random, 
{\textit{i.e.}, $\hb(\ppm_\mathrm{prev};\Ym)$ is used here as a control variate with $\boldsymbol\alpha = \mathbf{I}$ in \eqref{eq:cv_vector}, where $\mathbf{I}$ is the identity matrix.
These steps are illustrated in Algorithm \ref{alg:svrg}.
\renewcommand\algorithmiccomment[1]{%
 {#1}%
}
\begin{algorithm}[htb]
	\begin{algorithmic}
		\STATE Given $\eta$, $m$, and $N_h$.
		\STATE Initialize $\widetilde\ppm_1$.
		\FOR[({outer iteration})] {$j=1,2,\dots,N_{oit}$}
		\STATE Set ${\ppm}_{\mathrm{prev}} = \widetilde{\ppm}_{j}$.
		\STATE Set $\bar{\hb}(\ppm_\mathrm{prev}) = \frac{1}{N_h}\sum_{i=1}^{N_h} \hb(\ppm_\mathrm{prev};\Ym_i)$.
		\STATE Set $\ppm_1 = \ppm_{\mathrm{prev}}$.
		\FOR[({inner iteration})] {$k=1,2,\dots,m$}
		\STATE Uniformly at random choose $t \in \{1,2,\dots,N_h\}$.
		\STATE Set $\ppm_{k+1}\rightarrow\ppm_k-\eta \big[\hb(\ppm_k;\Ym_{t})-\hb(\ppm_\mathrm{prev};\Ym_{t})+\bar{\hb}(\ppm_\mathrm{prev})\big]$. 
		\STATE [see Eqn. (\ref{eq:svrg})]
		\ENDFOR
		\STATE $\widetilde\ppm_{j} \leftarrow \ppm_{m+1}$. 
		\ENDFOR
	\end{algorithmic}
	\caption{\textit{SVRG} \citep{johnson2013accelerating}}
	\label{alg:svrg}
\end{algorithm}

\section{Methodology: Proposed Algorithms}\label{sec:method}

In this section, inspired by the SAG and SVRG algorithms, we propose two variants of SGD using a combination of high- and low-fidelity gradient evaluations. 
\subsection{Bi-fidelity Stochastic Average Gradient (BF-SAG) Algorithm}\label{sec:bfsag}
Similar to batch implementation of the SAG algorithm, at every iteration of the proposed BF-SAG algorithm, we update $N_h$ gradients using the high-fidelity model. In addition, we update $N_l\gg N_h$ gradients using the low-fidelity model. 
Unlike in the SAG algorithm, by using many low-fidelity model evaluations to update most of the gradients, in addition to the high-fidelity model evaluations, we reduce the dependency on previous designs as in the optimization process designs may go through drastic changes over a few iterations.
Similar to the SAG algorithm in \eqref{eq:sag},  
%
parameter update in \eqref{eq:sgd} is performed at the $k$th iteration using $N$ stochastic gradients $\dm_{k,i}$.
Specifically, the search direction is defined as
\begin{equation}\label{eq:bfsag}
\begin{split}
\hb_k & = \frac{1}{N}\sum_{i=1}^{N}\dm_{k,i};\\
\dm_{k,i} &= \begin{cases}
\hb_\mathrm{low}(\ppm_k;\Ym_{i})\qquad \text{if $i\in \{t_l\}_{l=1}^{N_l}$};\\
\hb_\mathrm{high}(\ppm_k;\Ym_{i})\qquad \text{if $i\in \{t_h\}_{h=1}^{N_h}$};\\
\dm_{k-1,i}\qquad\qquad \text{otherwise},\\
\end{cases}
\end{split}
\end{equation}
where $t_l$ and $t_h$ are selected uniformly at random from $\{1,2,\dots,N\}$ with $\{t_l\}\cap\{t_h\}=\emptyset$.
The implementation of these steps is summarized in Algorithm \ref{alg:bfsag}.

\subsubsection{Computational cost}
Let the ratio of the computational effort for a low-fidelity model compared to a high-fidelity model be $\gamma<1$. Note that, $\gamma$ may include the cost of generating a high-fidelity gradient estimate from a low-fidelity gradient, \textit{e.g.}, interpolating a gradient computed from a coarse grid model on a fine grid. 
The computational cost in terms of the cost of high-fidelity gradient evaluations for $N_{it}$ number of iterations of the BF-SAG algorithm can be given by 
\begin{equation}
\begin{split}
    C_\mathrm{BF-SAG} &:= N_{it}\left[N_h +\gamma N_l\right].\\
\end{split}
\end{equation}
Next, we compare the per-iteration cost of the BF-SAG algorithm relative to that of a batch implementation of the SAG algorithm, where we update $N^\prime_h$ gradients per iteration. The ratio of their respective per-iteration cost is $\frac{N_h+\gamma N_l}{N^\prime_h}$. 
In general, $N^\prime_h > N_h$ and $\gamma$ is a very small number leading to per-iteration cost efficiency. Note that, the per-iteration costs used to evaluate the ratio do not necessarily lead to similar accuracy in a fixed number of iteration. 
\begin{algorithm}[htb]
	\begin{algorithmic}
		\STATE Given $\ppm_0$, $\eta$, $N_l$, and $N_h$.
		\STATE Initialize $\ppm_1=\ppm_0$ and $\dm = \bm{0}$.
		\FOR {$k=1,2,\dots,N_{it}$}
		\STATE Draw $N_l+N_h$ samples $\{t_b\}_{b=1}^{N_l+N_h}$ uniformly from $\{1,\cdots,N\}$. 
		\STATE Define $\{t_l\}_{l=1}^{N_l} \equiv \{t_b\}_{b=1}^{N_l}$ and $\{t_h\}_{h=1}^{N_h} \equiv\{t_b\}_{b=N_l+1}^{N_l+N_h}$. 
		\FOR{$i=1,2,\dots,N$} 
		\IF{$i\in \{t_l\}_{l=1}^{N_l}$}
		\STATE Compute $\dm_{k,i} := {\hb}_\mathrm{low}(\ppm_{k};\Ym_{t_b})$.
		\ELSIF{$i\in \{t_h\}_{h=1}^{N_h}$}
		\STATE Compute $\dm_{k,i} := {\hb}_\mathrm{high}(\ppm_{k};\Ym_{t_b})$.
		\ELSE
		\STATE Set $\dm_{k,i}=\dm_{k-1,i}$.
		\ENDIF
		\ENDFOR
		\STATE $\ppm_{k+1} \leftarrow \ppm_k - \frac{\eta}{N}\sum_{i=1}^N\dm_{k,i}$.
		\ENDFOR
	\end{algorithmic}
	\caption{\textit{Bi-fidelity Stochastic Average Gradient (BF-SAG)}}
	\label{alg:bfsag}
\end{algorithm}

\subsubsection{Convergence}
In this subsection, we present a result for linear convergence of the proposed BF-SAG algorithm. The corresponding assumptions and result that are presented in the following theorem are inspired from the results of \cite{schmidt2017minimizing}.

\newtheorem{theorem}{Theorem}
\begin{theorem}\label{theorem1}
Assume the objective function $J(\ppm)$ obtained from low- and high-fidelity models are strongly convex with constants $\mu_\mathrm{low}$ and $\mu_\mathrm{high}$, respectively. Also, assume that the corresponding gradients are Lipschitz continuous with constants $L_\mathrm{low}$ and $L_\mathrm{high}$, respectively. Let $\ppm^*=\argmin\limits_{\ppm}J(\ppm)$ and initialize the gradient history vector $\dm$ to zero. For some constants $\mu,L^2> 0$ that depend on the constants $\mu_\mathrm{low}$, $\mu_\mathrm{high}$ and $L_\mathrm{low}^2$, $L_\mathrm{high}^2$, respectively, Algorithm \ref{alg:bfsag} achieves a linear convergence as
\begin{equation}\label{eq:bfsag_proof}
    \Exp[\lVert \ppm_{k+1} - \ppm^*\rVert^2\lvert \ppm_{0}]  \leq (1 -  \mu^2/L^2)^k \lVert \ppm_{0} - \ppm^* \rVert ^2.
\end{equation}
\end{theorem}
The proof of this theorem along with the details of \eqref{eq:bfsag_proof} and the definitions of $\mu_\mathrm{low}$, $L_\mathrm{high}$, $L_\mathrm{low}$, and $L_\mathrm{high}$ are given are given in \ref{proof_thm1}. 
Note that, the structural optimization settings that we consider here do not satisfy the conditions of this theorem, \textit{e.g.}, strong convexity. However, the empirical results presented in Section \ref{sec:examples} illustrate convergence for the non-convex problems considered in this paper. 



%
%
\subsection{Bi-fidelity Stochastic Variance Reduced Gradient (BF-SVRG)}\label{sec:bfsvrg}
In our second algorithm, we exploit the availability of
low- and high-fidelity model gradients by using a control variate similar to its use in the SVRG algorithm \citep{johnson2013accelerating}, resulting in a bi-fidelity extension of SVRG named here BF-SVRG.
In particular, we employ the gradient ${\hb}_\mathrm{low}$ evaluated using the low-fidelity model and a previous estimate of the parameters $\thetaa_\mathrm{prev}$ as a control variate ($Y$ in \eqref{eqn:cv}). To estimate the mean of the control variate $\widehat{\hb}_\mathrm{low}$, we use $N_l$ random samples as follows
\begin{equation}\label{eq:cv_mean}
    \widehat\hb_\mathrm{low} = \frac{1}{N_l}\sum_{i=1}^{N_l} \hb_\mathrm{low}(\ppm_\mathrm{prev};\Ym_i).
\end{equation}
Note that, the number of random samples $N_l$ required to estimate $\widehat\hb_\mathrm{low}$ can be large \cite{rubinstein2016simulation} but the low-fidelity model is used keeping the cost small. 
Also, as in the standard SVRG (Algorithm \ref{alg:svrg}), $\ppm_\mathrm{prev}$ is updated only at every $m$ inner iteration and the same $\widehat{\hb}_\mathrm{low}$ is kept for all these inner iterations. At every inner iteration, $N_h\ll N_l$ random variables are used to estimate gradients using the high-fidelity model. Note that, these $N_h$ random variables are not necessarily a subset of the $N_l$ random samples used in \eqref{eq:cv_mean}). The estimated gradients are used to evaluate the mean as follows
\begin{equation}\label{eq:h_mean}
    \widehat\hb_\mathrm{high} = \frac{1}{N_h} \sum_{b=1}^{N_h} \hb_\mathrm{high}(\thetaa_k;\Ym_{b}).
\end{equation}
The gradient descent step in \eqref{eq:sgd} is performed with the search direction defined as 
\begin{equation}\label{eq:bfsvrg}
\hb_k = \widehat\hb_\mathrm{high}-\frac{\boldsymbol\alpha}{N_h}\sum_{b=1}^{N_h}\left(\hb_\mathrm{low}(\thetaa_\mathrm{prev};\Ym_{b})-\widehat\hb_\mathrm{low}\right), 
\end{equation}
where $\boldsymbol{\alpha}$ is a coefficient matrix and the same set of random samples $\{\Ym_b\}_{b=1}^{N_h}$ from \eqref{eq:h_mean} is used. 
The estimation of optimal coefficient matrix using \eqref{eq:opt_cv_matrix} requires computation of the inverse of the covariance matrix of $\hb_\mathrm{low}(\ppm_\mathrm{prev};\Ym)$. To avoid this computationally expensive step, herein we assume $\boldsymbol{\alpha}$ is a diagonal matrix with optimal diagonal entries given by \cite{wang2013variance}
\begin{equation}
    \alpha^*_{ii} = \frac{\mathbb{C}_{ii}}{\mathbb{V}_{ii}},\quad i=1,\dots,n_{\ppm},
\end{equation}
where $\mathbb{V}$ is the covariance matrix of $\hb_\mathrm{low}(\thetaa_\mathrm{prev};\Ym)$ and $\mathbb{C}$ is the cross-covariance between $\hb_\mathrm{high}(\thetaa_k;\Ym)$ and $\hb_\mathrm{low}(\thetaa_\mathrm{prev};\Ym)$. 
We use the $N_h$ random samples used in \eqref{eq:h_mean} and \eqref{eq:bfsvrg} to compute the diagonal entries of  $\boldsymbol{\alpha}$ as follows
\begin{equation}
    {\alpha}_{ii} = \frac{\left[\sum_{b=1}^{N_h}\left(\hb_\mathrm{low}(\thetaa_\mathrm{prev};\Ym_{b})-\widehat\hb_\mathrm{low}\right)\odot \left(\hb_\mathrm{low}(\thetaa_\mathrm{prev};\Ym_{b})-\widehat\hb_\mathrm{low}\right)\right]_{ii}}{\left[\sum_{b=1}^{N_h}\left(\hb_\mathrm{high}(\thetaa_k;\Ym_{b})-\widehat\hb_\mathrm{high}\right)~\!\odot \left(\hb_\mathrm{low}(\thetaa_\mathrm{prev};\Ym_{b})-\widehat\hb_\mathrm{low}\right)\right]_{ii}},\quad i=1,\dots,n_{\ppm},
\end{equation}
where $\odot$ denotes a Hadamard product (\textit{i.e.}, elementwise multiplication).
Note that, the SVRG algorithm uses the high-fidelity gradient from the past iteration as the control variate and an identity matrix as the $\boldsymbol{\alpha}$. On the other hand, the proposed BF-SVRG algorithm uses the low-fidelity gradient as the control variate and uses a diagonal coefficient matrix $\boldsymbol{\alpha}$, which is not necessarily an identity matrix. 
A description of these steps is shown in Algorithm \ref{alg:bfsgdcv}.

\subsubsection{Computational Cost}
From Algorithm \ref{alg:bfsgdcv}, the computational cost of the BF-SVRG algorithm can be given by
\begin{equation}\label{eq:cost_bfsvrg}
\begin{split}
    C_{BF-SVRG} &= N_{oit}\left[\gamma N_l + \gamma mN_h+mN_h \right]\\
    & = N_{oit}\left[\gamma N_l + (\gamma+1)mN_h\right],
\end{split}
\end{equation}
where $N_{oit}$ is the number of outer iterations and $\gamma$ is the ratio of the computational cost of the low- over high-fidelity gradient evaluations, as in Section \ref{sec:bfsag}.
Next, we compare the per outer iteration cost of the BF-SVRG algorithm with the SVRG algorithm, where both use $m$ inner iterations, but the SVRG algorithm uses $N_h^\prime$ gradient evaluations to estimate the mean of the control variate $\widehat\hb(\ppm_\mathrm{prev})$ in \eqref{eq:svrg_cv} and a batch of $N_h^{{\prime\prime}}$ random samples to estimate $\hb(\ppm_k;\Ym)$ and $\hb(\ppm_\mathrm{prev};\Ym)$ in \eqref{eq:svrg}. Using \eqref{eq:cost_bfsvrg}, the ratio of per outer iteration cost of the BF-SVRG and SVRG algorithms is given by $\frac{\gamma N_l+(\gamma+1)mN_h}{N_h^\prime+2N_h^{{\prime\prime}} m}$. 
In general, $N_h^\prime\approx N_l\gg N_h$ and $N_h^{{\prime\prime}}\approx N_h$, thus BF-SVRG leads to cost efficiency. Note that, the per-iteration costs used to evaluate the ratio do not necessarily lead to similar accuracy in a fixed number of iteration. 

\begin{algorithm}[htb!]
	\begin{algorithmic}
		\STATE Given $\ppm_0$, $\eta$, $N_h$, $N_l$, and $m$.
		\STATE Initialize $\widetilde\ppm_1=\ppm_0$.
		\FOR[({outer iteration})] {$j=1,2,\dots,N_{oit}$}
		\STATE Set ${\ppm}_{\mathrm{prev}} = \widetilde{\ppm}_{j}$.
		\STATE Set $\widehat\hb_\mathrm{low} = \frac{1}{N_l}\sum_{i=1}^{N_l} \hb_\mathrm{low}(\ppm_\mathrm{prev};\Ym_i),\quad N_l$ sufficiently large.
		\STATE Set $\ppm_1 = \ppm_{\mathrm{prev}}$.
		\FOR[({inner iteration})] {$k=1,2,\dots,m$}
		\STATE Draw $N_h\ll N_l$ samples $\{\Ym_b\}_{b=1}^{N_h}$ 
		\STATE Compute $\widehat\hb_\mathrm{high} = \frac{1}{N_h} \sum_{b=1}^{N_h} \hb_\mathrm{high}(\thetaa_k;\Ym_{b})$.
\STATE Set $\boldsymbol{\alpha}$ as a diagonal matrix with diagonal entries 
${\alpha}_{ii} = \frac{\left[\sum_{b=1}^{N_h}\left(\hb_\mathrm{low}(\thetaa_\mathrm{prev};\Ym_{b})-\widehat\hb_\mathrm{low}\right)\odot \left(\hb_\mathrm{low}(\thetaa_\mathrm{prev};\Ym_{b})-\widehat\hb_\mathrm{low}\right)\right]_{ii}}{\left[\sum_{b=1}^{N_h}\left(\hb_\mathrm{high}(\thetaa_k;\Ym_{b})-\widehat\hb_\mathrm{high}\right)~\!\odot \left(\hb_\mathrm{low}(\thetaa_\mathrm{prev};\Ym_{b})-\widehat\hb_\mathrm{low}\right)\right]_{ii}},\quad i=1,\dots,n_{\ppm}$.
		\STATE Set $\ppm_{k+1}\rightarrow\ppm_k-\eta \left[\widehat\hb_\mathrm{high}-\frac{\boldsymbol\alpha}{N_h}\sum_{b=1}^{N_h}\left(\hb_\mathrm{low}(\thetaa_\mathrm{prev};\Ym_{b})-\widehat\hb_\mathrm{low}\right)\right]$.
		\ENDFOR
		\STATE $\widetilde\ppm_{j} \leftarrow \ppm_{m+1}$. 
		\ENDFOR
	\end{algorithmic}
	\caption{\textit{Bi-fidelity Stochastic Variance Reduced Gradient (BF-SVRG)}}
	\label{alg:bfsgdcv}
\end{algorithm}
\subsubsection{Convergence}
In this subsection, we present a theorem to show that the convergence rate of the proposed BF-SVRG algorithm is linear and depends on the correlation between the high- and low-fidelity gradients. The corresponding assumptions and result that are presented in the following theorem are inspired from the results of \cite{schmidt2017minimizing,johnson2013accelerating}. 
\begin{theorem}\label{theorem2}
Assume the objective function $J(\ppm)$ obtained from high-fidelity models is strongly convex with a constant $\mu_\mathrm{high}$ and the corresponding gradients are Lipschitz continuous with constant $L_\mathrm{high}$. Let $\ppm^*=\argmin\limits_{\ppm}J(\ppm)$. For some constant $\delta\geq 1$, Algorithm \ref{alg:bfsgdcv} with $m$ inner iterations and $j$ outer iteration achieves a linear convergence as 
\begin{equation}
    \Exp[\lVert \widetilde\ppm_{j} - \ppm^*\rVert^2\lvert \ppm_0] \leq\left(1-\frac{\mu^2_\mathrm{high}}{2L^2_\mathrm{high}\delta}\right)^{jm} \lVert \ppm_{0} - \ppm^* \rVert ^2,
\end{equation}
where $\widetilde\ppm_{j}$ is the value of the parameter after $j$ outer iterations.
\end{theorem}
The proof of this theorem is given in \ref{proof_thm2}. We note that the constant $\delta$ depends on the correlation between high- and low-fidelity gradients at a given $\bm\theta$ and has a minimum value of 1 as defined in \ref{proof_thm2}. A higher correlation between the high- and low-fidelity gradients leads to a $\delta$ closer to 1, hence a tighter bound on $\Exp[\lVert \widetilde\ppm_{j} - \ppm^*\rVert^2\lvert \ppm_0]$. 
Similar to Theorem \ref{theorem1}, the structural optimization settings that we consider here do not satisfy the conditions of this theorem, \textit{e.g.}, strong convexity. However, the empirical results presented in Section \ref{sec:examples} illustrate the convergence for non-convex problems considered in this paper. 
\section{Numerical Examples}
\label{sec:examples}

\subsection{Example I: Polynomial Regression}

To show the working principle of the proposed algorithms, we first consider the following fourth order polynomial as the high-fidelity model
\begin{equation}
    y_\mathrm{high}(x) = 2+5x + 1.75x^2+5x^3+6.5x^4,
\end{equation}
where $x\in[-1,1]$. The low-fidelity model near $x_0\in[-1,1]$ is given by a second order Taylor series expansion, \textit{i.e.},
\begin{equation}\label{eq:toy_low}
    y_\mathrm{low}(x;x_0) = y_{\mathrm{high}}(x_0) + \frac{\mathrm{d}y_{\mathrm{high}}}{\mathrm{d}x}\Bigg\rvert_{x_0} (x-x_0) +  \frac{1}{2}\frac{\mathrm{d}^2y_{\mathrm{high}}}{\mathrm{d}x^2}\Bigg\rvert_{x_0} (x-x_0)^2.
\end{equation}
We seek to approximate $y_{\mathrm{high}}(x)$ in a fourth order monomial basis as 
\begin{equation}
y_\mathrm{pred}(x,\bm\theta)=\sum_{i=0}^4\theta_{i}x^i,  \qquad \theta_i\in\mathbb{R},
\label{eqn:y_pred}
\end{equation}
via the standard least sqaures regression over the unknown coefficients $\thetaa=[\theta_0,\theta_1,\theta_2,\theta_3,\theta_4]^T$. To this end, we generate $N=1000$ noisy measurements of $y_\mathrm{high}$ using
\begin{equation}
    y_{\mathrm{obs}}(x_i) = y_{\mathrm{high}}(x_i)+e_i,\quad i=1,2,\dots,N,
\end{equation}
where $x_i$ is randomly selected from $[-1,1]$ and $e_i$ are independent, zero-mean normal random variables with variance $0.25$. We then seek to find an optimal $\bm\theta$ in (\ref{eqn:y_pred}) by solving the optimization problem 
\begin{equation}
    \min\limits_{\thetaa}~~J(\thetaa) =\frac{1}{N}\sum_{i=1}^N \Big(y_\mathrm{obs}(x_i)-y_{\mathrm{pred}}(x_i,\thetaa)\Big)^2.
\end{equation}
\subsubsection{Results}

We implement the proposed algorithms for this problem, where  the low-fidelity model is defined by searching for the $x_0\in\{-1.0,-0.75,\dots,0.75,1.0\}$ closest to $x_i$, for each $x_i$, and using a Taylor series expansion about $x_0$ (see \eqref{eq:toy_low}). The optimization is performed with an initial guess $\thetaa = [1.5,4,1,4,5]^T$. 
Figure \ref{fig:sagbfsagExI} shows a comparison of mean squared error (MSE) values for the SAG and the BF-SAG algorithms with $\eta=0.25$. Note that, when the total number of gradient evaluations per iteration is the same, \ie~$N_l+N_h=50$, the performance of the BF-SAG algorithm is similar or a little worse. However, most of the gradient evaluations are performed with the low-fidelity model, thus making the BF-SAG algorithm cheaper. On the other hand, if we increase the number of low-fidelity  gradient evaluations to 230, the convergence of the BF-SAG algorithm is faster. Note that, the number of high-fidelity gradient evaluation in the BF-SAG algorithm at every iteration is still smaller than that of the SAG algorithm.
This advantage will be exploited in the rest of the numerical examples in this paper.
\begin{figure}[htb!]
\centering
    \includegraphics[scale=0.3]{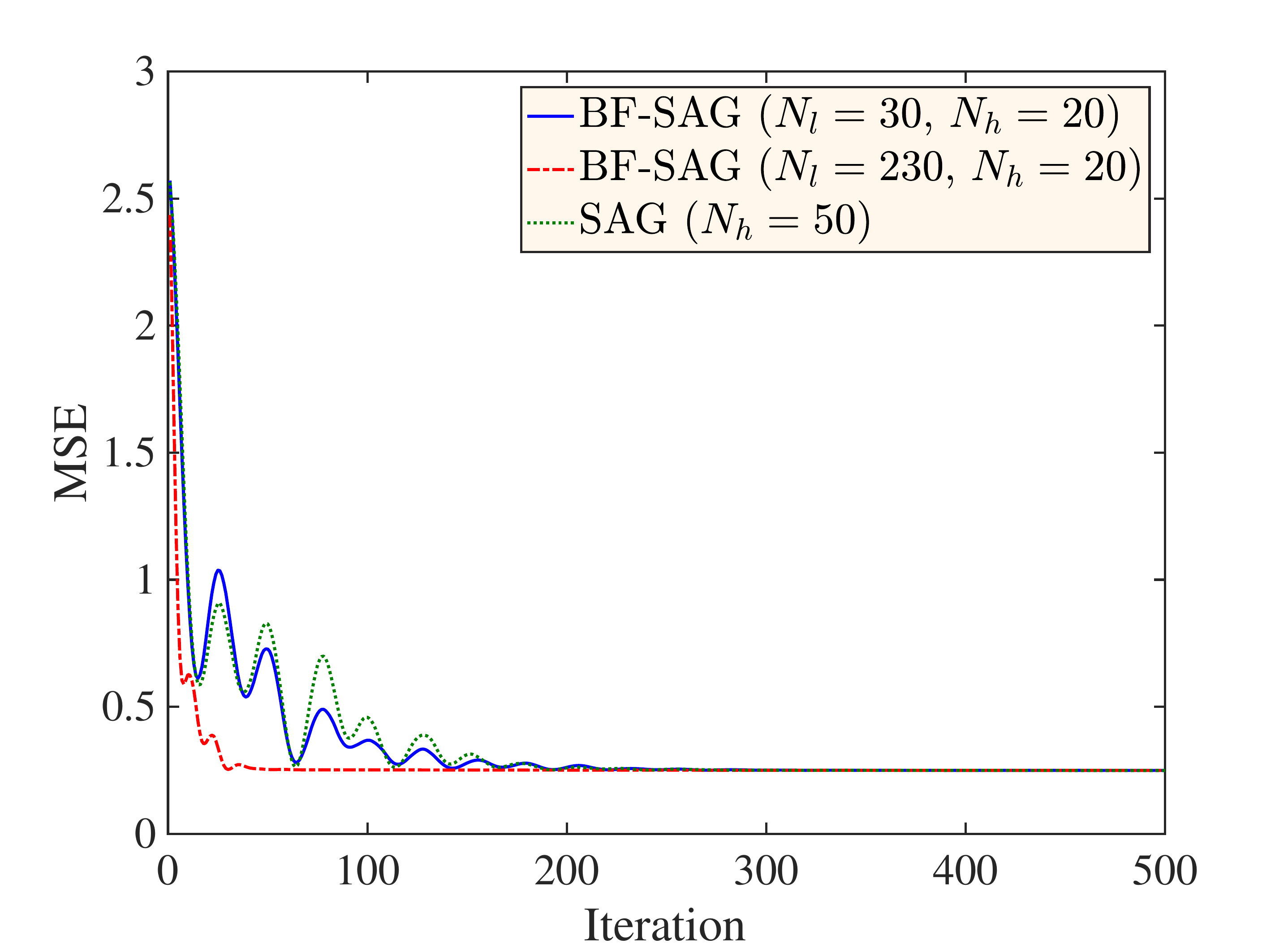}
\caption{Comparison of MSE reduction achieved by the SAG and BF-SAG algorithms for Example I.}\label{fig:sagbfsagExI}
\end{figure}

Similarly, we compare the SVRG algorithm with the BF-SVRG algorithm in Figure \ref{fig:svrgbfsvrgExI} with $\eta=0.25$ and $m=20$.  
When using the same number of gradient evaluations in each outer iteration ($320+20\times2=360$ for SVRG and $200+20\times8=360$ for BF-SVRG in Figure \ref{fig:svrgbfsvrgExI}), the performance of SVRG and BF-SVRG are comparable. However, the BF-SVRG algorithm uses significantly smaller number of high-fidelity gradients. The convergence is improved by using more low- and high-fidelity gradient evaluations to estimate the diagonal entries of the coefficient matrix $\boldsymbol{\alpha}$ as shown in Figure \ref{fig:svrgbfsvrgExI}. When compared to the BF-SAG algorithm, the BF-SVRG algorithm with $N_l=200$ and $N_h=16$ has similar convergence to the BF-SAG algorithm with $N_l=230$ and $N_h=20$ and both show small oscillations in the objective after some initial iterations. The relative expected error in the estimates of the optimization parameters estimated using 100 independent runs of the optimization algorithms is shown in Figure \ref{fig:ExI_conv_rate}, which shows that the proposed algorithms have a linear convergence after a few initial iterations. Note that, a linear convergence of the algorithms is not shown for the next two examples as the true value of the parameters $\ppm^*$ are unknown in those two examples. 
\begin{figure}[htb!]
\centering
    \centering
    \includegraphics[scale=0.3]{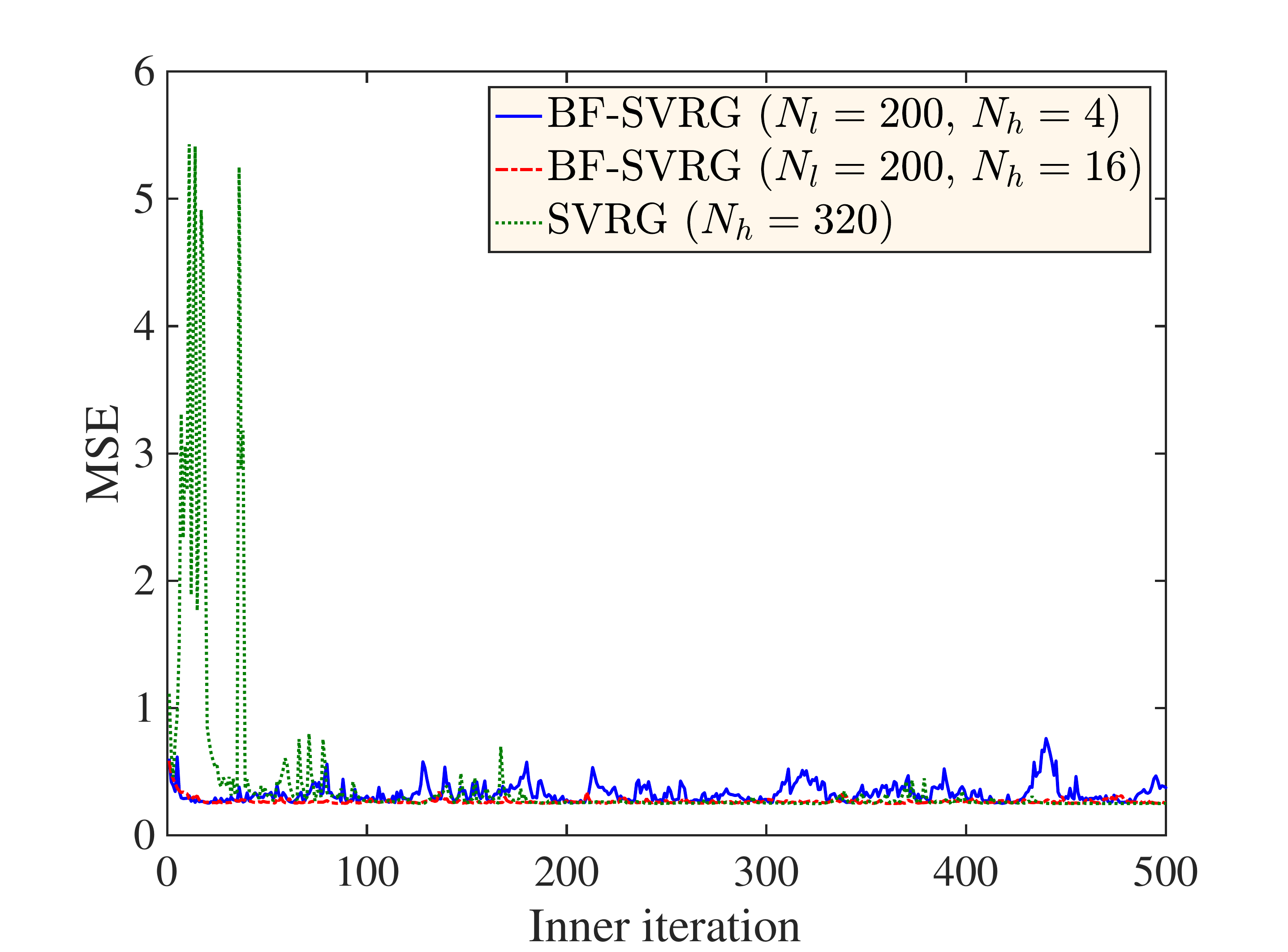}
\caption{Comparison of MSE reduction achieved by the SVRG and BF-SVRG algorithms for Example I.}\label{fig:svrgbfsvrgExI}
\end{figure}
\begin{figure}[htb!]
\centering
    \centering
    \includegraphics[scale=0.3]{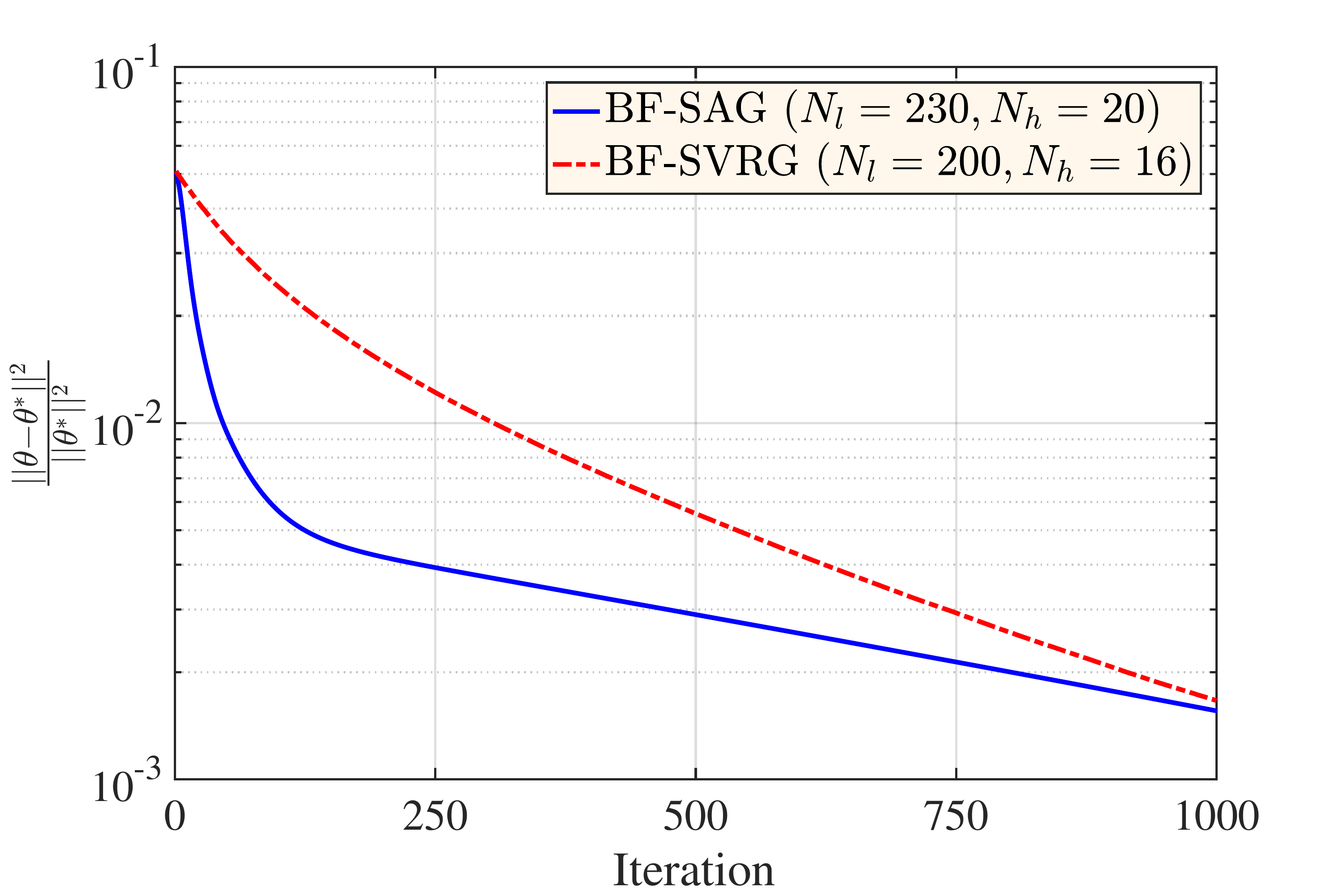}
\caption{Relative expected error in the estimates of the optimization parameters using the BF-SAG and BF-SVRG algorithms in Example I.}\label{fig:ExI_conv_rate}
\end{figure}
%

%
%
\subsection{Example II: Shape Optimization of a Plate with a Hole}

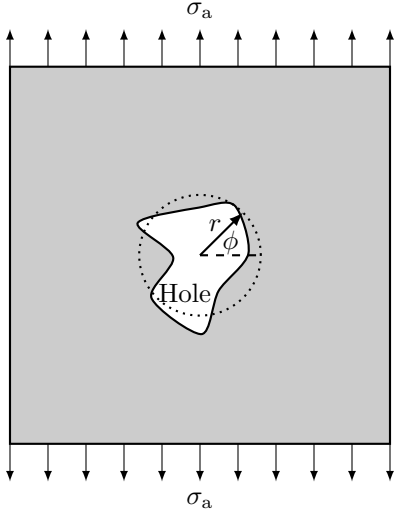
\begin{figure}
    \centering
    \begin{tikzpicture}
    \pgfmathsetseed{8}
    \draw[thick,fill=gray!40] (-2.5,-2.5) rectangle (2.5,2.5);
    \draw[ thick, fill =white] plot [smooth cycle, samples=8,domain={1:8},fill=white] (\x*360/8+2.5*rnd:0.25+1*rnd) node at (0,0) {};
    \foreach \x in {0,...,10}{
        \draw[-latex] (-2.5+0.5*\x,2.5) -- (-2.5+0.5*\x,3);}
    \foreach \x in {0,...,10}{
        \draw[-latex] (-2.5+0.5*\x,-2.5) -- (-2.5+0.5*\x,-3);}
    \node[draw=none] at (0,3.25) {$\sigma_\mathrm{a}$};
    \node[draw=none] at (0,-3.25) {$\sigma_\mathrm{a}$};
    
    \node[draw=none] at (-0.2,-0.5) {Hole};
    \draw[-latex,thick] (0,0) -- (0.5657,0.5657);
    \draw[thick,dashed] (0,0) -- (0.8,0);
    \node[draw=none] at (0.4,0.165) {$\phi$}; 
    \node[draw=none] at (0.2,0.4) {$r$}; 
    \draw[thick,dotted] (0,0) circle (0.8);
    \end{tikzpicture}
    \caption{The plate is subjected to tensile stress on two sides. The shape of the hole is parameterized.}
    \label{fig:holeplate}
\end{figure}

The second example is concerned with optimizing the shape of an elastic $2D$ square plate of dimension $20\times20$ with a hole located at the center; see Figure \ref{fig:holeplate}. The plate is subject to a uni-axial uniform stress, $\sigma_\mathrm{a}$, and the goal is to minimize the maximum principal stress, \ie\, the stress intensity factor.  The shape of the hole is described in the polar coordinate $(r,\phi)$ with the center of the coordinate system placed at the center of the plate. The radius of the hole is described in a harmonic basis,
\begin{equation}\label{eq:hole}
    \begin{split}
        & r = r_0 + \tau \sum_{i=1}^d \frac{1}{\sqrt{i}}\left[ \theta^s_i \sin(i\phi)+\theta^c_i \cos(i\phi) \right];\\
        & -1\leq \theta^s_i,\theta^c_i \leq 1,\qquad i=1,\dots,d,\\
    \end{split}
\end{equation}
where $\thetaa = [\theta^s_1,\dots,\theta^s_d,\theta^c_1,\dots,\theta^c_d]^T$ is the vector of optimization variables. The parameters in \eqref{eq:hole} are set to $r_0=1$, $\tau=0.15$, and $d=3$. Note that, these values will not result in a negative radius. 
We further add a contribution from the deviation of the area of the hole from a circle of radius one to the objective to avoid converging to the solution of a hole with much smaller radius. 
The optimization problem is given by
\begin{equation}
\begin{split}
  &  \min\limits_{\thetaa} ~~J(\thetaa) = \Exp\left[ \frac{\sigma_{\max}(\thetaa)}{\sigma_\mathrm{a}}\right]+\lambda (\pi-A_\mathrm{hole}(\thetaa))\\
 &    \text{{subject to}}~~ {-1}\leq \theta_j \leq {1}, \quad j=1,\dots,\np,
\end{split}
\end{equation}
where $\sigma_{\max}(\thetaa)$ is the maximum value of the principal stress in the plate; $A_\mathrm{hole}(\thetaa)$ is the area of the hole; and we choose $\lambda=50$. The box constraint on $\ppm$ is applied here by restricting any parameter update to be within $[-1,1]$, \ie if the update puts the parameter outside of $[-1,1]$ we replace the parameter with $-1$ or $1$, respectively. 
We assume uncertainty in the tensile stress $\sigma_\mathrm{a}$ applied at the two ends of the plate, which we model as
\begin{equation}
    \sigma_\mathrm{a}(\xi_\sigma) = \sigma_0(1+0.5\xi_\sigma),
\end{equation}
where $\xi_\sigma$ is a standard normal random variable and $\sigma_0=1$ is a constant.
Further, the elastic modulus of the plate $E$ and Poisson's ratio $\nu$ are assumed uncertain and given by
\begin{equation}
    \begin{split}
        E(\xi_E)&=E_0(1+0.05\xi_E);\\
        \nu(\xi_\nu)&=\nu_0(1+0.01\xi_\nu),\\
    \end{split}
\end{equation}
where $\xi_E$ is a standard normal random variable truncated on one side to keep $E(\xi_E)\in(0,\infty)$; $\xi_\nu$ is standard normal random variable truncated at both sides to get $\nu(\xi_\nu) \in (0.3,0.5)$; $E_0=1000$; and $\nu_0=0.3$. To compute the gradient of $J(\bm\theta)$ we use a differentiable approximation of $\sigma_{\max}/\sigma_0\approx \left( \sum_i (\sigma_i/\sigma_0)^p \right)^{1/p}$ for a large (even) $p$, where $\sigma_i$ are nodal principal stresses \cite{holmberg2013stress}. In our numerical experiments, we set $p=200$.

To solve for the stress in the plate we use the finite element package FEniCs \cite{AlnaesBlechta2015a,LoggMardalEtAl2012a}. We employ finite differencing to compute the gradients of $J(\bm\theta)$.
Note that, only a quarter of the plate is analyzed leveraging the symmetry of the problem.
A low-fidelity model with $\sim260$ degrees of freedom is constructed using a coarse mesh while a high-fidelity model with $\sim10250$ degrees of freedom uses refined mesh around the hole. For a circular hole with a radius of $r=1$, the coarse mesh gives a relative difference of 9.1369\% in $\sigma_{\max}/\sigma_\mathrm{a}$ as compared to the fine mesh. 
Figure \ref{fig:mesh} shows the typical meshes used for both high- and low-fidelity models. We remesh each of these models once for every $\ppm$ but keep the total number of degrees of freedom approximately same.

\begin{figure}[htb!]
\centering
\begin{subfigure}[t]{0.45\textwidth}
    \centering
    \includegraphics[scale=0.175]{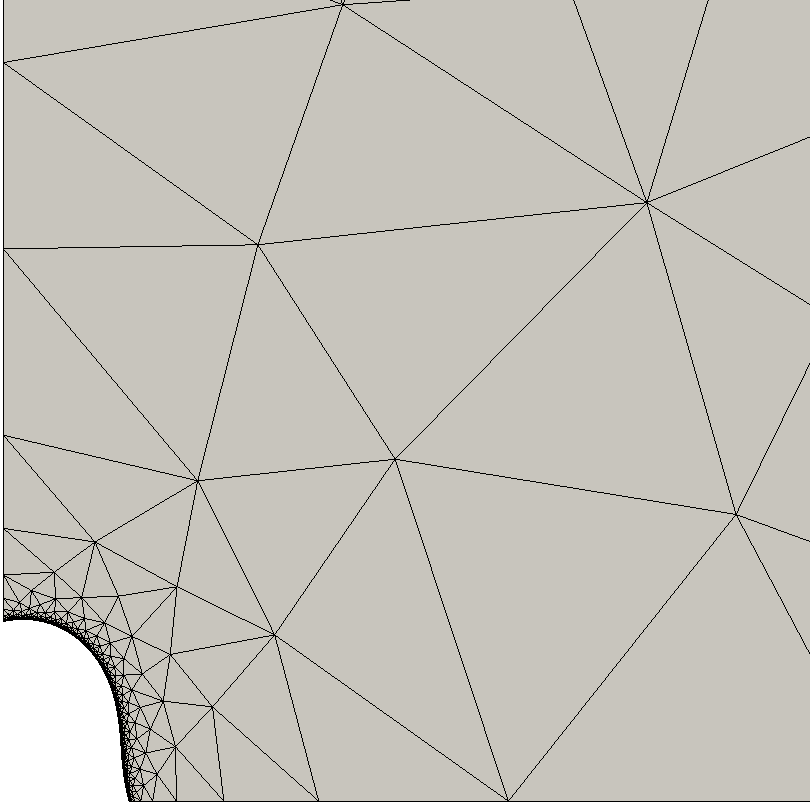}
    \caption{High-fidelity mesh.}\label{fig:HF}
\end{subfigure}
\begin{subfigure}[t]{0.45\textwidth}
    \centering
    \includegraphics[scale=0.175]{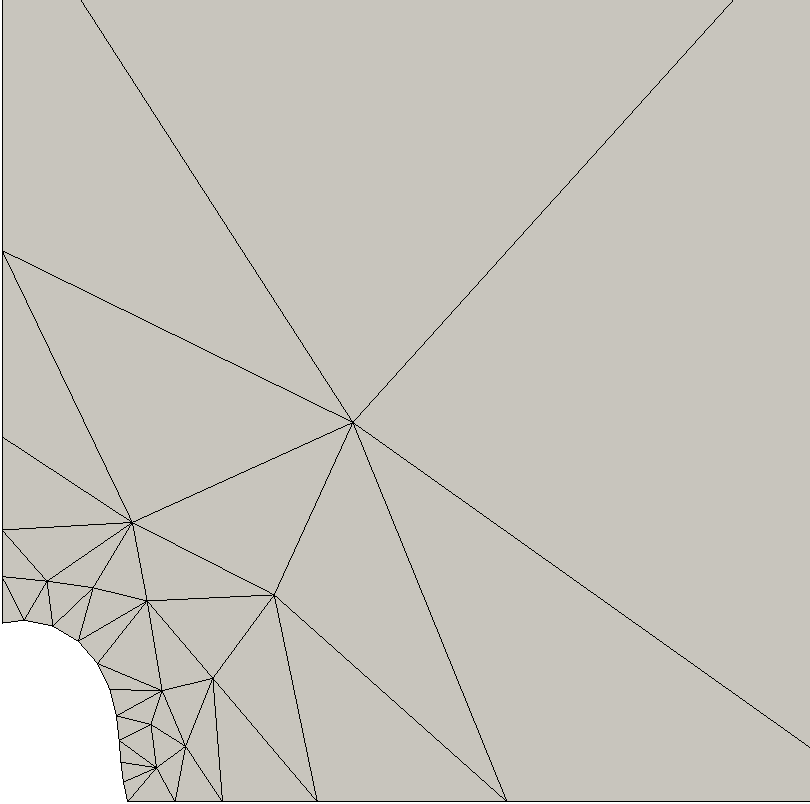}
    \caption{Low-fidelity mesh.}\label{fig:LF}
\end{subfigure}
\caption{Meshes used for high- and low-fidelity models (for an arbitrary hole that we also use as initial shape for the optimization).}\label{fig:mesh}
\end{figure}



Here, we assume the computational cost $C_T$ for a model with $N_\mathrm{d}$ degrees of freedom is proportional to $N_\mathrm{d}^\beta$, where $\beta$ depends on the solver used. Hence, we can write
\begin{equation}
    \gamma = \frac{C_{T,\mathrm{low}}}{C_{T,\mathrm{high}}} = \frac{N_\mathrm{d,low}^{\beta}}{N_\mathrm{d,high}^{\beta}},
\end{equation}
where $C_{T,\mathrm{low}}$ and $C_{T,\mathrm{high}}$ are computational costs for a low- and a high-fidelity gradient evaluations, respectively, including any interpolation costs; and $N_\mathrm{d,low}$ and $N_\mathrm{d,high}$ are degrees of freedom for a low- and a high-fidelity model, respectively.
Solving the plate problem with an algebraic multigrid solver and using the wall-clock data for meshes with different number of degrees of freedom, we determine via data fitting a $\gamma$ value approximately of $0.015$ in this example. Note that, $\gamma$ changes slightly as we remesh but the total number of degrees of freedom remains almost same and hence the change in $\gamma$ is insignificant. 
The cost calaculations are performed on a desktop with a quad-core Intel Xeon(R) processor W3550 3.07GHz and 12 GB of memory running Ubuntu. 

\subsubsection{Results}

We use the hole shown in Figure \ref{fig:optimized_hole} (dashed line) generated with $\thetaa^s=[-0.1660,0.4406,-0.9998]^T$ and $\thetaa^c=[-0.3953,-0.7065,-0.8153]^T$ as the initial guess for all the optimization algorithms. Figure \ref{fig:ex2_bfsag} compares the performance of the SAG and the BF-SAG algorithms for this example with a learning rate $\eta=0.02$. In our implementation of the SAG algorithm, we use $N_h=10$, \ie ~10 high-fidelity gradient evaluations per iteration to keep the computational cost reasonable on a desktop computer. In the BF-SAG algorithm, at every iteration, we use $N_h=5$ and $N_l=5$, \ie\ we update five gradients using the high-fidelity model as before but for the other five gradients we use the low-fidelity model.
The evolution of the objective for the SAG and the BF-SAG algorithms is shown in Figure \ref{fig:ex2iter_eval_bfsag}, which shows similar performance for both of these algorithms. 
Figure \ref{fig:ex2norm_cost_eval_bfsag} further shows that if we normalize the computational cost in terms of the high-fidelity model evaluation, we can reach the optimum objective using a fraction of the cost in the BF-SAG algorithm compared to the SAG algorithm, which only uses the high-fidelity gradients. 

\begin{figure}[htb!]
\centering
\begin{subfigure}[t]{\textwidth}
    \centering
    \includegraphics[scale=0.3]{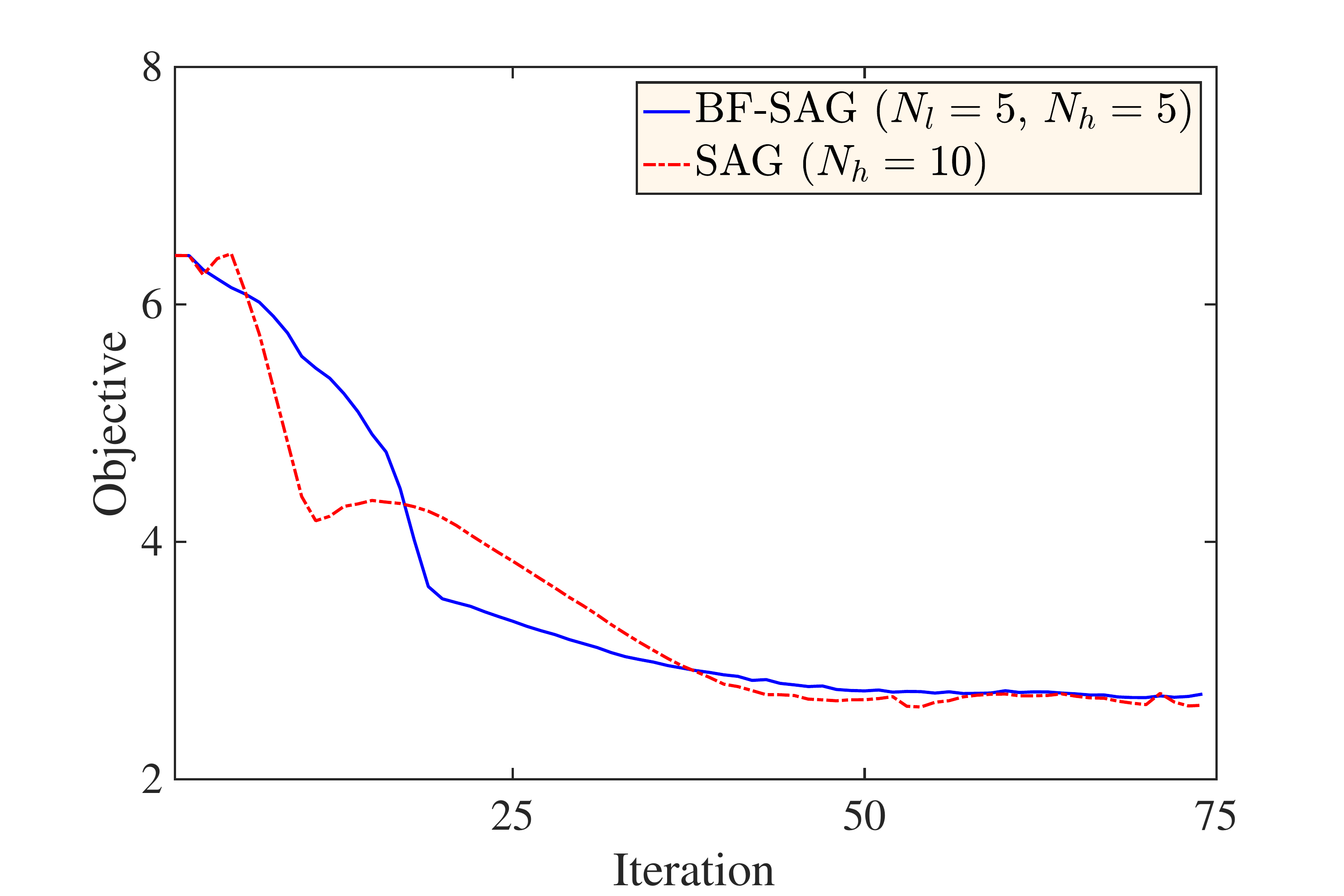}
    \caption{Objective vs. iteration.}\label{fig:ex2iter_eval_bfsag}
\end{subfigure}\\
\begin{subfigure}[t]{\textwidth}
    \centering
    \includegraphics[scale=0.3]{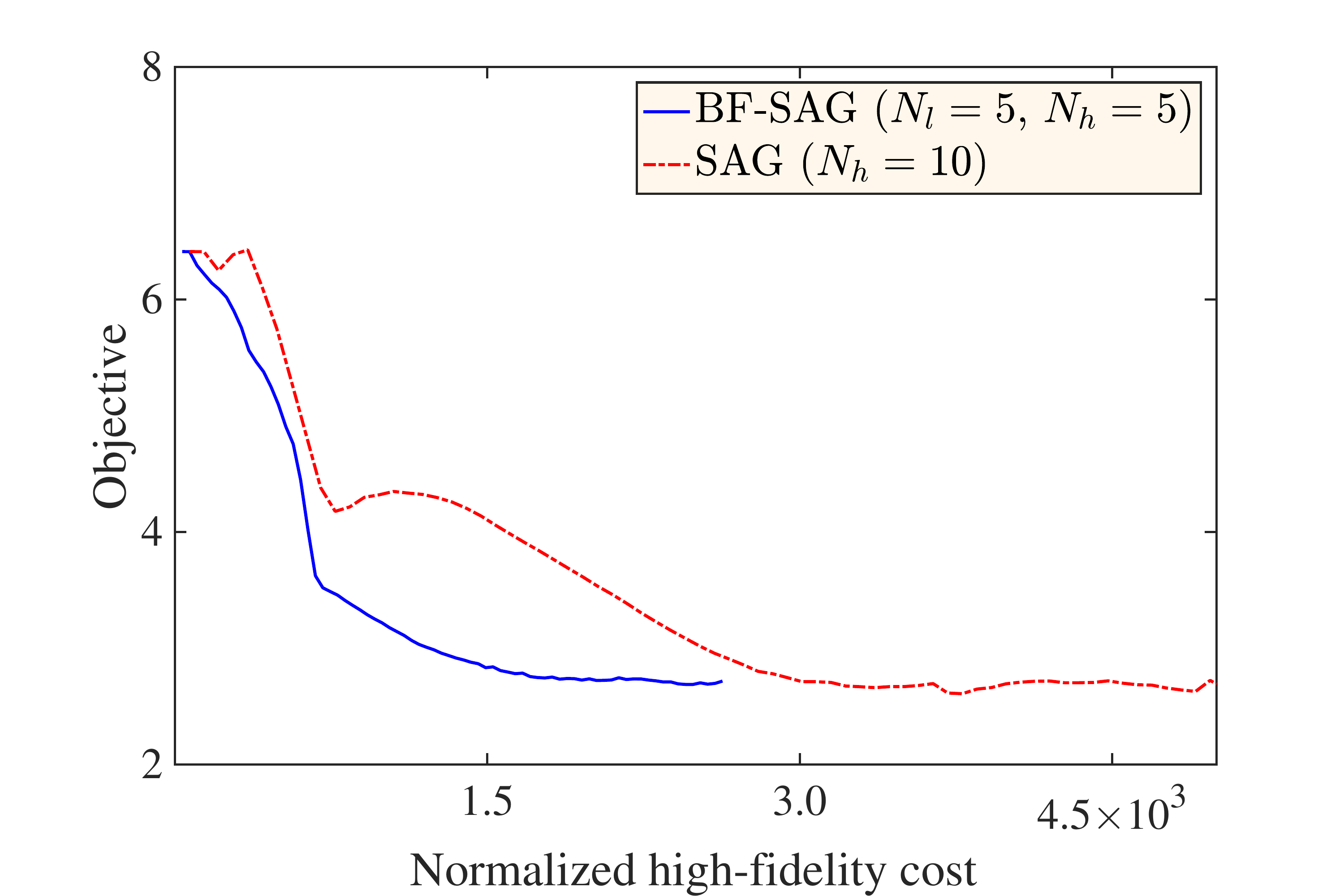}
    \caption{Objective vs. normalized high-fidelity cost.}\label{fig:ex2norm_cost_eval_bfsag}
\end{subfigure}
\caption{Comparison of the SAG and the proposed BF-SAG algorithms for Example II. Figure (a) shows that the BF-SAG and SAG algorithms have similar convergence. In Figure (b) normalized costs are calculated with respect to the cost of one high-fidelity model evaluation and shows that the computational cost of the BF-SAG algorithm is considerably smaller than the SAG algorithm.}\label{fig:ex2_bfsag}
\end{figure}

Next, we compare the SVRG and the proposed BF-SVRG algorithms for a learning rate $\eta=0.025$, where we keep the number of gradient evaluations for every outer iteration the same. Here, we consider the SVRG algorithm with $N_h=15$ and the BF-SVRG algorithm with $N_l=5$ and $N_h=2$. We use inner iteration $m=5$, which requires 25 gradient evaluations for every outer iteration for each of these algorithms. 
Figure \ref{fig:ex2iter_eval_bfsvrg} shows that we obtain a similar convergence, as the iteration progresses. 
However, in terms of the computational cost of evaluating the high-fidelity models, the BF-SVRG algorithm features a faster convergence as shown in Figure \ref{fig:ex2norm_cost_eval_bfsvrg}. 
Hence, the use of the BF-SAG and BF-SVRG algorithms are effective in reduction of the computational cost of the OuU in this example. 
Figure \ref{fig:optimized_hole} shows the final optimized shapes of the hole using the BF-SAG and BF-SVRG algorithms along with the initial shape. 


\begin{figure}[htb!]
\centering
\begin{subfigure}[t]{\textwidth}
    \centering
    \includegraphics[scale=0.3]{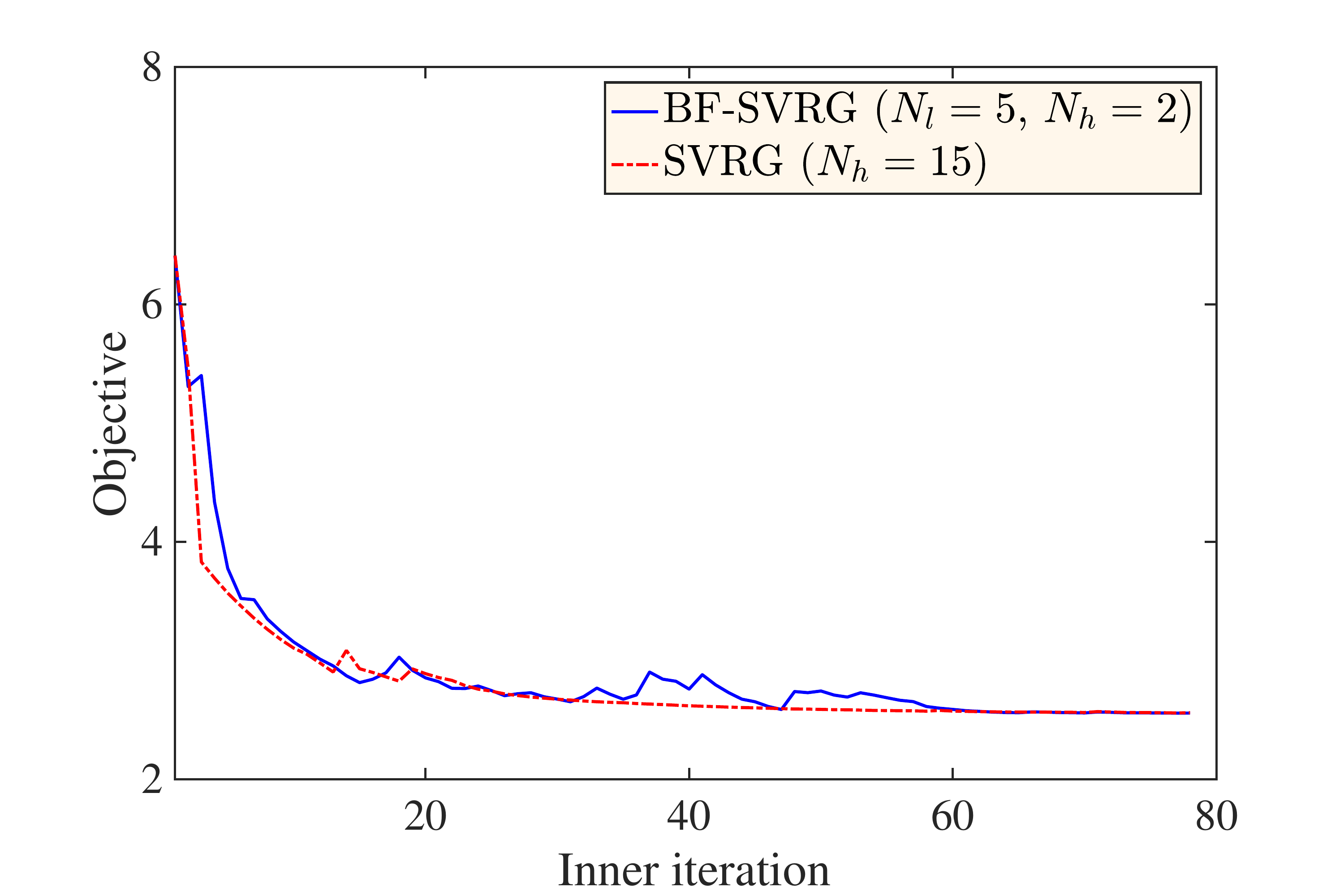}
    \caption{Objective vs. iteration.}\label{fig:ex2iter_eval_bfsvrg}
\end{subfigure}\\
\begin{subfigure}[t]{\textwidth}
    \centering
    \includegraphics[scale=0.3]{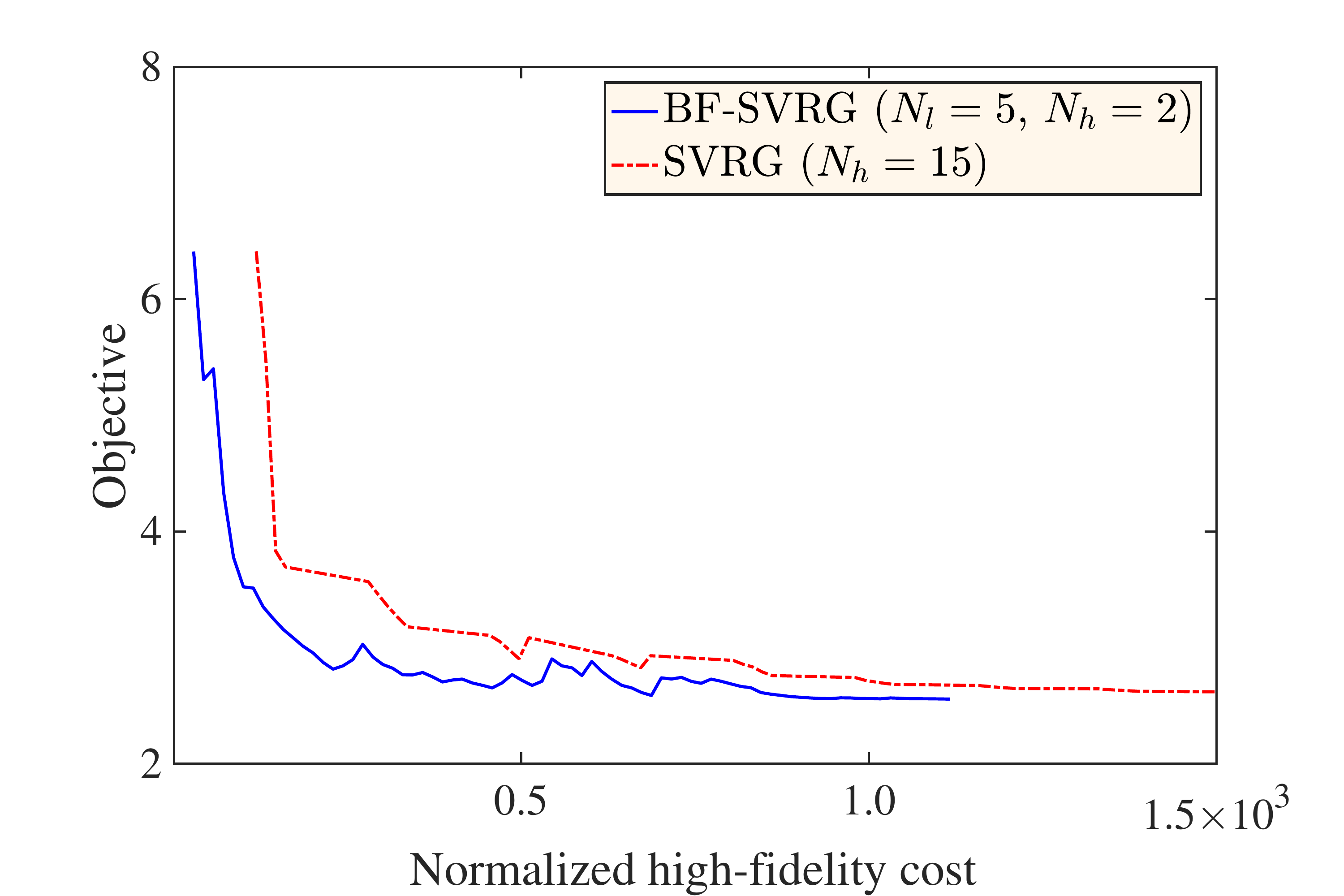}
    \caption{Objective vs. normalized high-fidelity cost.}\label{fig:ex2norm_cost_eval_bfsvrg}
\end{subfigure}
\caption{Comparison of SVRG and the proposed BF-SVRG algorithms for Example II. Figure (a) shows that the BF-SVRG and SVRG algorithms have similar convergence. In Figure (b) normalized costs are calculated with respect to the cost of one high-fidelity model evaluation and shows that the computational cost of the BF-SVRG algorithm is considerably smaller than the SVRG algorithm.}\label{fig:ex2_bfsvrg}
\end{figure}

\begin{figure}[htb!]
\centering
    \includegraphics[scale=0.3]{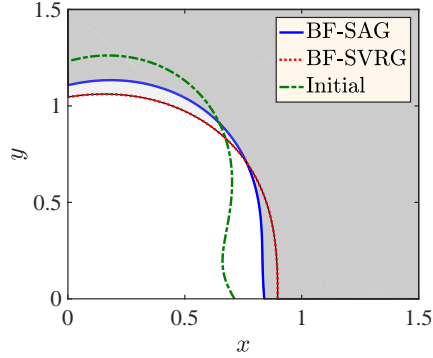}
    \caption{Optimized shape of the hole using the proposed bi-fidelity algorithms. Initial shape is shown using a dashed line. Note that, two-fold symmetry is assumed here.}
    \label{fig:optimized_hole}
\end{figure}

\subsection{Example III (a): Topology Optimization of a Beam under Uncertain Load}
For the third numerical example, we consider topology optimization of a beam problem. 
The design domain is simply supported on the bottom left and right ends with a span of $L$ and is subjected to an uncertain point load $2P$ at the mid span. 
The schematic for this problem is shown in Figure \ref{fig:prims_schematic}. Using the symmetry of the problem, we only consider one-half of the span as our optimization domain. The uncertainty in the load $P$ is given by
\begin{equation}\label{eq:load_unc}
    P(\xi) = P_0(1+0.5\xi),
\end{equation}
where $\xi$ is a uniform random variable in $[0,1]$ and $P_0=1$ is a constant. 
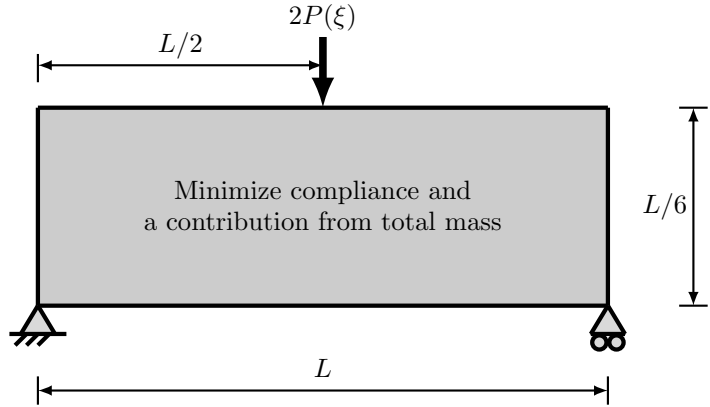
\begin{figure}[htb!]
	\centering
	\begin{tikzpicture}[scale=0.75,every node/.style={minimum size=1cm},on grid]
	\draw [ultra thick,fill=gray!40,draw=none] (0,-3.5) rectangle (10,0);
	\draw [ultra thick] (0,0) -- (10,0);
	\draw [ultra thick] (10,0) -- (10,-3.5);
	\draw [ultra thick] (10,-3.5) -- (0,-3.5);
	\draw [ultra thick] (0,-3.5) -- (0,0);
	
	\draw[ultra thick,fill=gray!30] (9.85,-4.15) circle (0.125);
	\draw[ultra thick,fill=gray!30] (10.15,-4.15) circle (0.125);
	\draw[fill=gray!30]    (10,-3.5) -- ++(0.3,-0.5) -- ++(-0.6,0) -- ++(0.3,0.5);
	\draw [ultra thick] (10,-3.5) -- (10.3,-4);
	\draw [ultra thick] (10.3,-4) -- (9.7,-4);
	\draw [ultra thick] (9.7,-4) -- (10,-3.5);
	
	\draw[fill=gray!30]    (0,-3.5) -- ++(0.3,-0.5) -- ++(-0.6,0) -- ++(0.3,0.5);
	\draw [ultra thick] (0,-3.5) -- (0.3,-4);
	\draw [ultra thick] (0.5,-4) -- (-0.5,-4);
	\draw [ultra thick] (-0.3,-4) -- (0,-3.5);
	\draw [ultra thick] (0,-4) -- (-0.2,-4.2);
	\draw [ultra thick] (-0.2,-4) -- (-0.4,-4.2);
	\draw [ultra thick] (0.2,-4) -- (0,-4.2);
	
	\node[draw=none] at (5, -1.5)  (c)     {Minimize compliance and};
	\node[draw=none] at (5, -2)  (c)     {a contribution from total mass};
	
	\draw[-latex, line width=1mm] (5,1.25) -- (5,0);
	\node[draw = none] at (5,1.6) () {$2P(\xi)$};
	
	
	
	\draw[thick,latex-latex] (0,-5) -- (10,-5);
	\node[draw = none] at (5,-4.6) () {$L$};
	\draw[thick] (0,-4.75) -- (0,-5.25);
	\draw[thick] (10,-4.75) -- (10,-5.25);
	
	\draw[thick,latex-latex] (0,0.75) -- (5,0.75);
	\node[draw = none] at (2.5,1.1) () {$L/2$};
	\draw[thick] (0,0.95) -- (0,0.45);
	\draw[thick] (5,0.95) -- (5,0.45);
	
	\draw[thick,latex-latex] (11.5,0) -- (11.5,-3.5);
	\node[draw = none] at (11,-1.75) () {$L/6$};
	\draw[thick] (11.25,0) -- (11.75,0);
	\draw[thick] (11.25,-3.5) -- (11.75,-3.5);
	\end{tikzpicture}
	\caption{Schematic of the beam problem of Example III (a) (design domain is shown as the shaded region).}
	\label{fig:prims_schematic}
\end{figure}

We optimize the material distribution inside the optimization domain by minimizing a combination of the compliance (\ie\ strain energy) and the mass subject to satisfying the equilibrium equations \cite{sigmund200199,sigmund2013topology}. 
We divide the design domain $\Omega$ into a large number of non-overlapping elements $\{\Omega_i\}_{i=1}^{N_e}$ using a finite element approach{, where $\{v_i\}_{i=1}^{N_e}$ are  the corresponding volumes and $N_e$ is the total number of elements used}. 
We use the Solid Isotropic Material with Penalization (SIMP) approach \citep{bendsoe1989optimal,ZhouSIMP91,sigmund2013topology} to formulate this topology optimization problem, and use some parts of the widely-used 99 line topology optimization code in \cite{sigmund200199}.
In SIMP, the material properties are interpolated by a power-law model in terms of the density $\rho$ of a fictitious porous material, \eg
\begin{equation}\label{eq:simp}
E(\rho_i) = \rho_i^{\beta_\mathrm{P}} E_0;\qquad 0<\rho_i\leq 1;\qquad i=1,2,\dots,N_e,
\end{equation}
where $\beta_\mathrm{P}$ is a penalization parameter and $E_0$ is the bulk material's elastic modulus. For the formulation of the optimization problem considered here and $\beta_\mathrm{P}>1$, intermediate densities are penalized as compared to densities closer to zero or one. We use $\beta_\mathrm{P}=3$ in the present work. 
To avoid a checker-board design, we use filtered values of the design variables $\ppm$ to define the material density $\rhoo$ \cite{bruns2001topology,bourdin2001filters,sigmund2007morphology,andreassen2011efficient}. This density filter is applied to the $e$th element as follows
\begin{equation}
    {\rho_e} = \frac{1}{ \sum_{i=1}^{N_e}{H}_i}\sum_{i=1}^{N_e} {H}_i\theta_i, 
\end{equation}
where the weight ${H}_i = \max\left(0, r_f-d_{ie} \right)$ is the difference between a filter size $r_f$ and the distance $d_{ie}$ between the centers of $i$th and $e$th elements. Herein, we use 1.5 times the element width as $r_f$.
Further use of projections may be needed to achieve a discrete $0-1$ design \cite{sigmund2013topology}. However, we do not use any such projection in this paper. 
We write the optimization problem as 
\begin{equation}\label{eq:top_obj}
\begin{split}
&\mathop{\min}\limits_{\ppm}~~J(\ppm) = \Exp\left[\sum_{i=1}^{N_e} \int_{\Omega_i} W\Big(u_i(\rho_i(\ppm);\Ym),\rho_i(\ppm);\xi\Big)\mathrm{d}V_i\right] + \lambda\sum_{i=1}^{N_e}v_i\rho_i(\ppm)\\
&\text{{subject to}}~~\Km(\rhoo(\ppm);\Ym) \uu(\rhoo(\ppm);\Ym) = \ff(\Ym);\\ 
&\text{\phantom{subject to}}~~ {0}\leq \rho_i(\ppm)\leq {1}, \text{  for }i=1,\dots,N_e,\\
\end{split}
\end{equation} 
where the objective $J(\cdot)$ is the expected value of the integral of the strain energy density $W(\cdot,\cdot;\cdot)$ plus a contribution from the total mass of the structure; $\Km(\cdot;\cdot)$ is the stiffness matrix; $\ff(\cdot)$ is the external force vector; $\lambda$ is the weighting factor for the contribution of the total mass to the objective.
The strain energy density depends on the displacement $\uu(\cdot;\cdot)$ and the material density $\rhoo$. The displacements $\uu(\cdot;\cdot)$ in turn depends on the uncertain variable $\xi$. 

We construct the high-fidelity model by dividing the domain of optimization (\ie~an area of $L/2\times L/6$) into $120\times40$ quadrilateral elements, which results in 4800 optimization variables. For the low-fidelity model, we use $60\times20$ of the same type of elements with 1200 optimization variables. 


To estimate the gradients using the low-fidelity model, we follow the steps shown in Figure \ref{fig:low_fidelity_grads}. 
First, we map the density variable $\rhoo$ from the high-fidelity mesh to the low-fidelity one by averaging, \ie\ a restriction like operation. Note that, we have one element in the low-fidelity mesh in place of four in the high-fidelity mesh. Next, we perform the calculation of the gradients using the low-fidelity mesh. Finally, we map the gradients to the elements in the high-fidelity mesh using a cubic spline interpolation, \ie\ a prolongation like operation. 
Note that, these restriction and prolongation operations are different than used in a multigrid scheme \cite{briggs2000van}. While for the configurations considered here the mapping is simple and can be computed analytically, it can be generalized to any mesh configurations using a proper projection operator, such as an $L^2$ minimization.
Averaging the computation time over 10 runs measured using \texttt{cputime} shows that for a \textsc{Matlab} implementation of the finite element solver and the proposed mapping scheme the cost of such a gradient estimate is 10.47 times cheaper than that of direct calculation of the high-fidelity gradients. This leads to the cost ratio of low- and high-fidelity gradient evaluations $\gamma=0.096$. Note that, the restriction and prolongation costs make up 6.76\% of the total low-fidelity gradient calculation cost. We use the same desktop computer as in the previous example to compute the costs.

\begin{figure}
    \centering
    \includegraphics[scale=0.25]{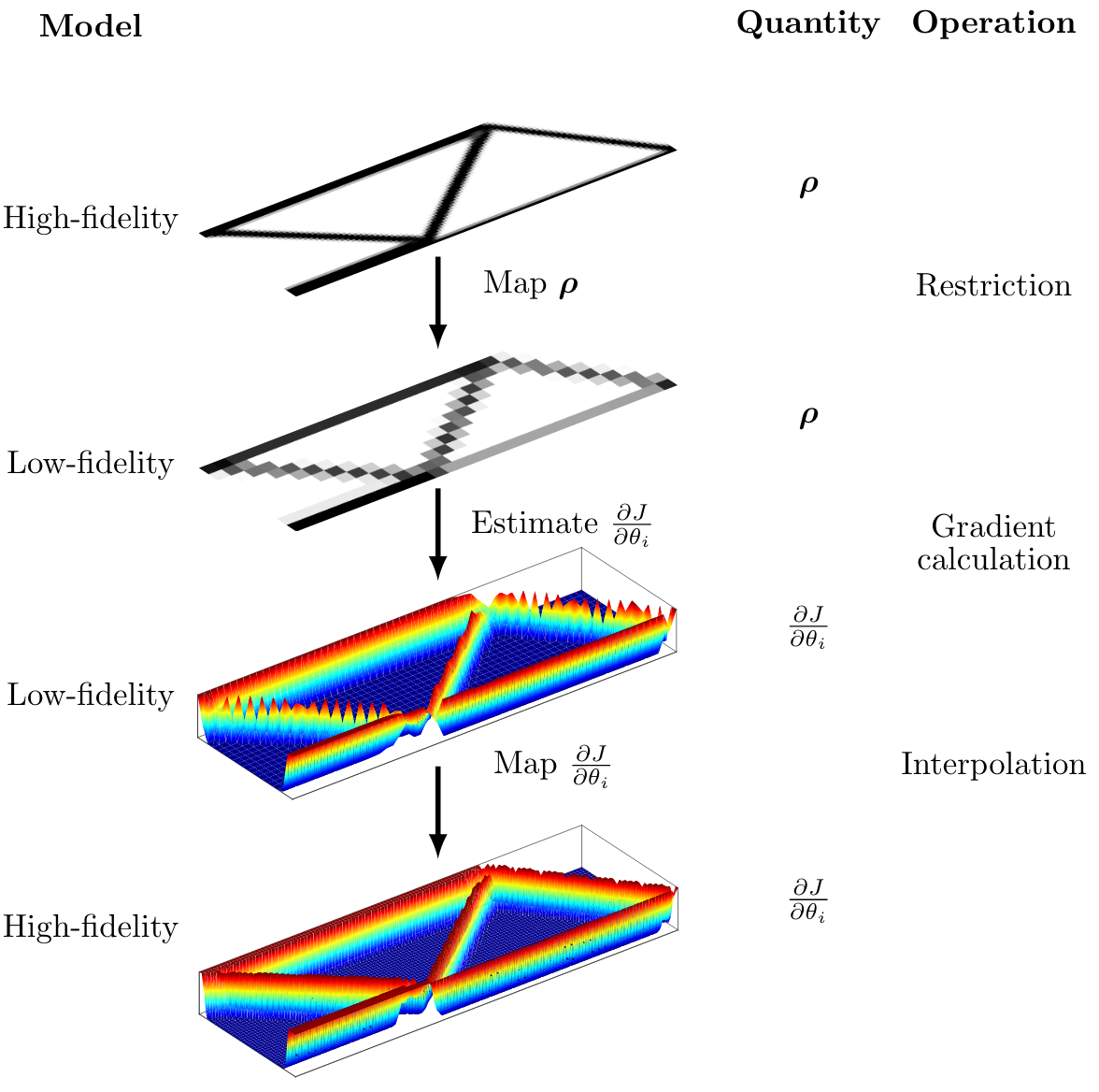}
    \caption{The steps to estimate the gradients using the low-fidelity model are shown here. The density variables $\rho$ for the finer mesh (high-fidelity model) are first mapped onto the coarser mesh (low-fidelity model). Then the coarser mesh is used to calculate the gradients. Finally, the coarse grid gradients are mapped again onto the finer mesh using a cubic spline interpolation.}
    \label{fig:low_fidelity_grads}
\end{figure}

\subsubsection{Results}

We study the proposed algorithms for this example with  a learning rate $\eta=0.05$ and $\lambda = 1$ in \eqref{eq:top_obj}. The final designs obtained from the SAG and BF-SAG algorithms with different number of low-fidelity gradient solves are shown in Figure \ref{fig:ex3_des1}. The first two designs are similar and exhibit a truss-like topology as in the deterministic optimization problem in \cite{sigmund200199}. The final design (Figure \ref{fig:bfsag_des2}) that uses more low-fidelity gradient solves has a smaller mass compared to the other two designs.

\begin{figure}[htb!]
\centering
\begin{subfigure}[t]{\textwidth}
    \centering
    \includegraphics[scale=0.25]{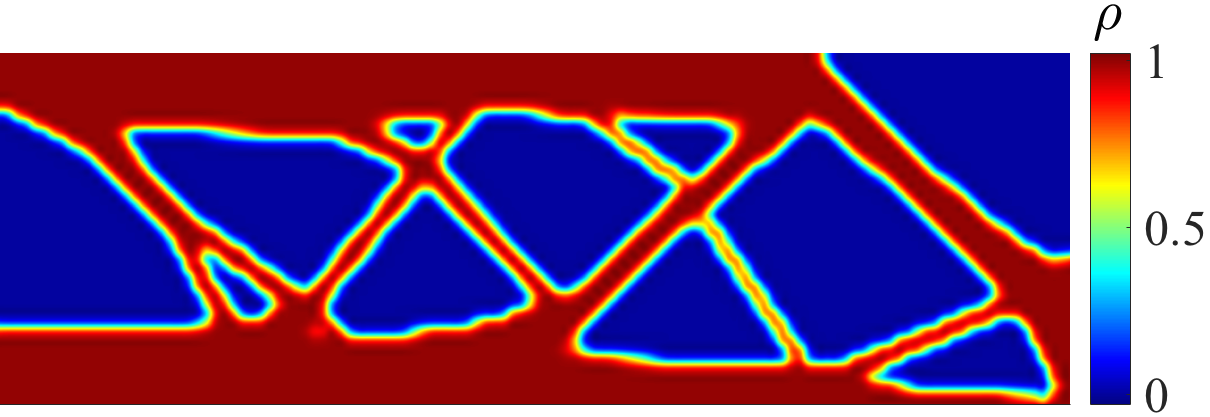}
    \caption{Final design using the SAG algorithm with $N=25$.}\label{fig:sag_des}
\end{subfigure}\\
\begin{subfigure}[t]{\textwidth}
    \centering
    \includegraphics[scale=0.25]{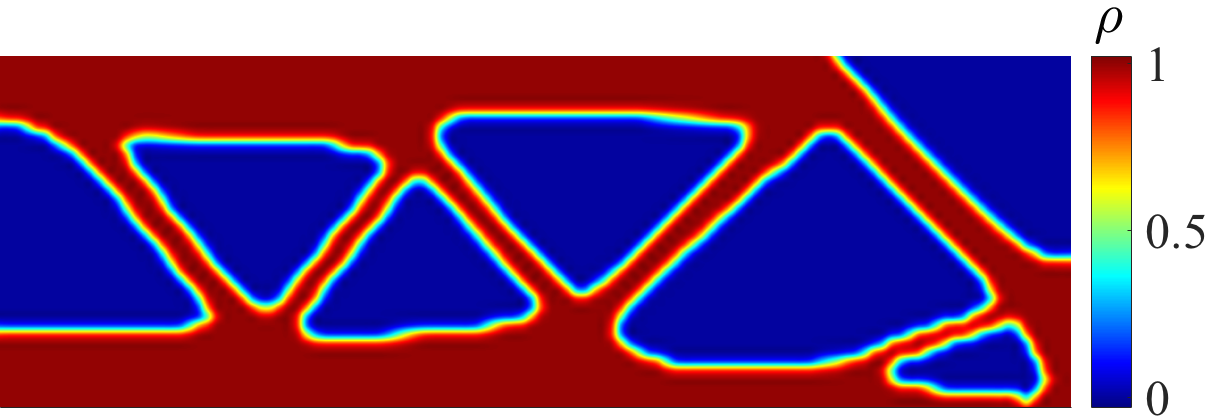}
    \caption{Final design using the BF-SAG algorithm with $N_l=20$ and $N_h=5$.}\label{fig:bfsag_des1}
    \end{subfigure}
    \\
\begin{subfigure}[t]{\textwidth}
    \centering
    \includegraphics[scale=0.25]{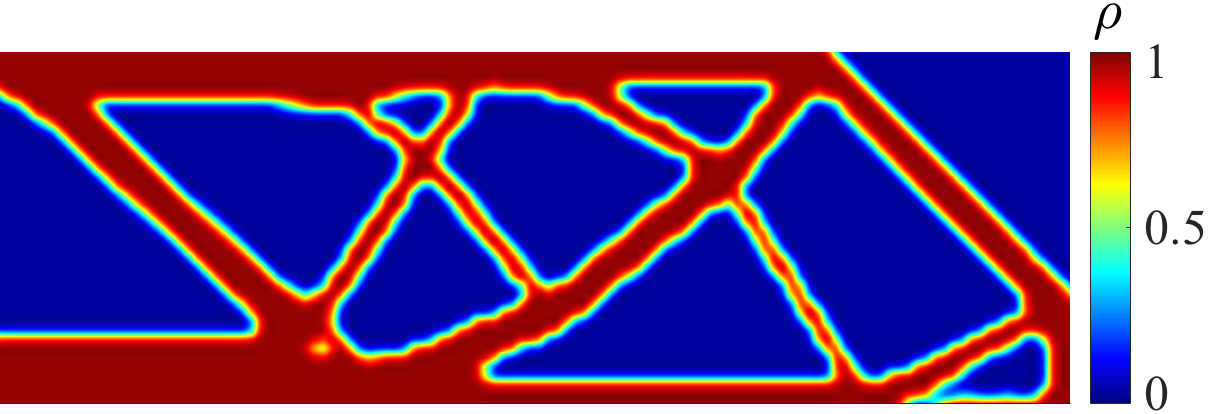}
    \caption{Final design using the BF-SAG algorithm with $N_l=95$ and $N_h=5$.}\label{fig:bfsag_des2}
\end{subfigure}
\caption{Final designs for Example III (a) as obtained from the SAG and the proposed BF-SAG algorithms. Note that, both algorithms use the same number of finite element solves per iteration in subfigures (a) and (b). However, the BF-SAG algorithm only uses 5 high-fidelity solutions compared to 25 in the SAG algorithm. In subfigure (c), we use more low-fidelity solutions per iteration.}\label{fig:ex3_des1}
\end{figure}

The objective is plotted in Figure \ref{fig:Obj1_ex3a} for these two algorithms. The result shows that the performance of the BF-SAG algorithm (solid blue curve) is comparable to the SAG algorithm (dotted green curve), when we use the same number of finite element solves per iteration. Further, if we increase the number of low-fidelity solves per iteration, we can improve the convergence as shown by the (red) dash-dotted curve. 
In terms of the mass of the structure, from Figure \ref{fig:Mass1_ex3a}, we see that the BF-SAG algorithm with more low-fidelity gradient solves produces a structure that has a significantly smaller mass but similar total objective; note that the total objective includes a contribution from the mass. 
\begin{figure}[htb!]
\centering
\begin{subfigure}[t]{\textwidth}
    \centering
    \includegraphics[scale=0.3]{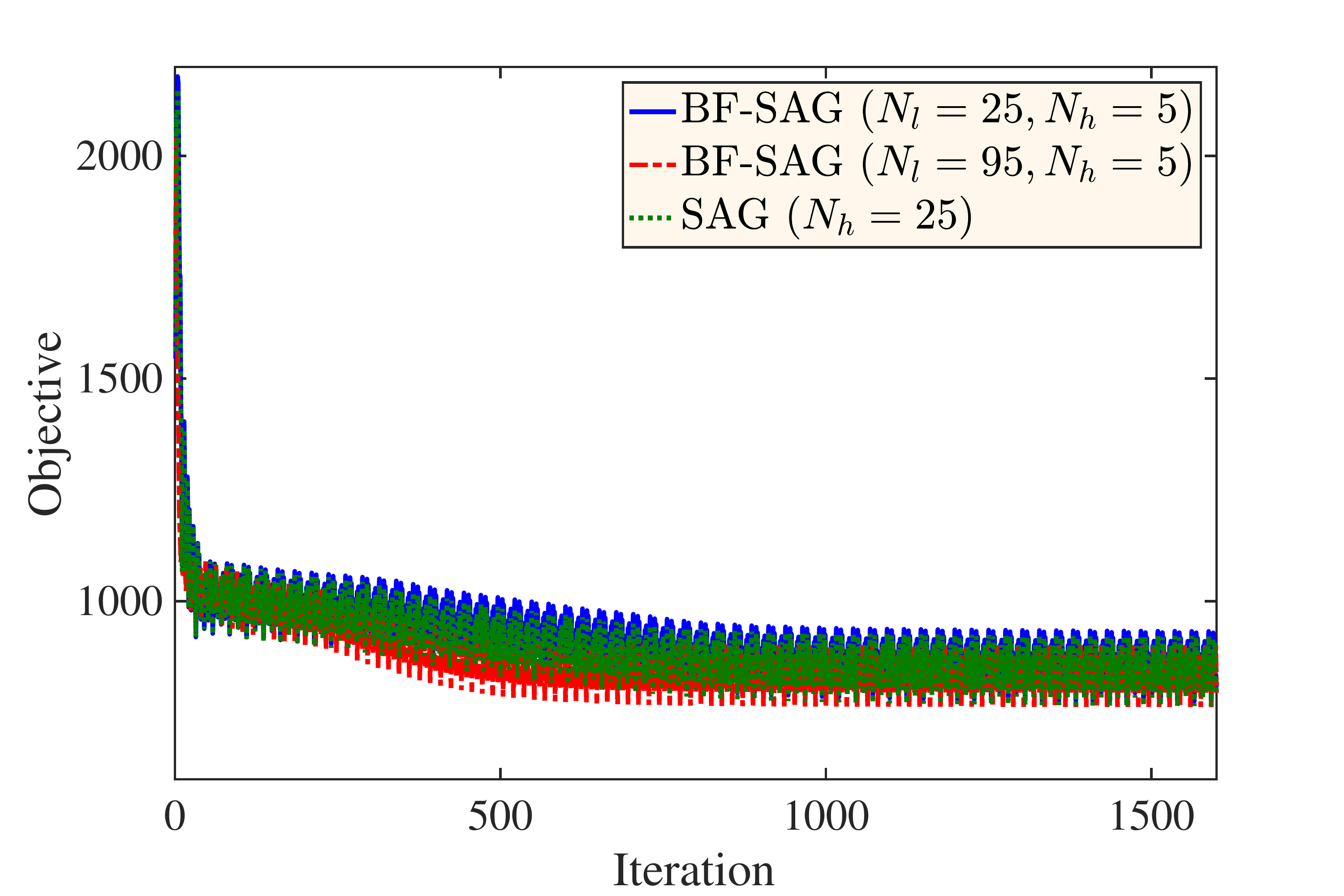}
    \caption{Objective.}\label{fig:Obj1_ex3a}
\end{subfigure}\\
\begin{subfigure}[t]{\textwidth}
    \centering
    \includegraphics[scale=0.3]{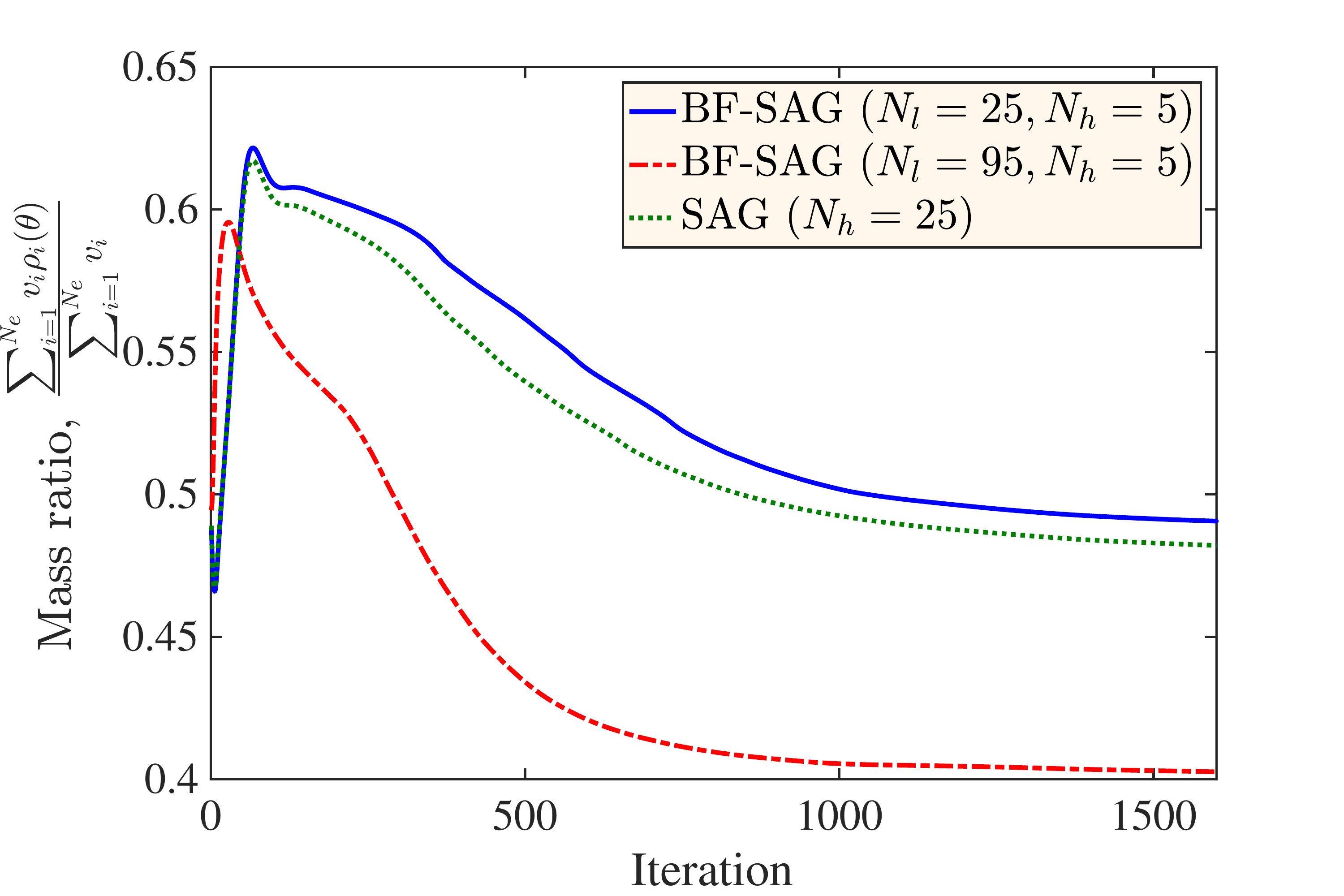}
    \caption{Mass ratio.}\label{fig:Mass1_ex3a}
\end{subfigure}
\caption{
Reduction of objective and mass ratio during the optimization process for two configurations of the BF-SAG algorithm and one configuration of the SAG algorithm for Example III (a).
Note that the faster convergence to the optimum can be achieved by using more low-fidelity models and only a handful of high-fidelity models per iteration.}\label{fig:ex3_Obj1}
\end{figure}

Next, we compare the SVRG and BF-SVRG algorithms for this example. Again, we obtain similar truss-like designs (see Figure \ref{fig:ex3_des2}). 
\begin{figure}[htb!]
\centering
\begin{subfigure}[t]{\textwidth}
    \centering
    \includegraphics[scale=0.25]{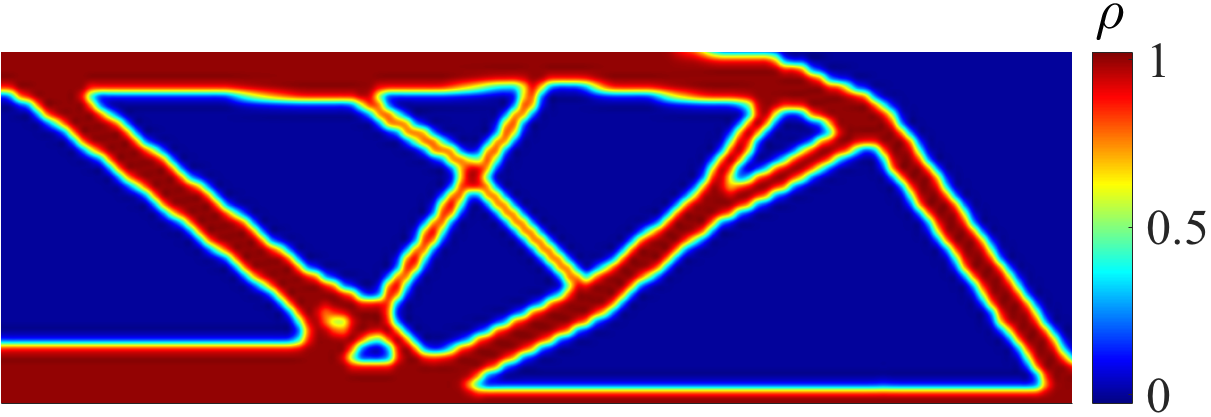}
    \caption{Final design using the SVRG algorithm with $N_h=80$.}\label{fig:svrg_des}
\end{subfigure}\\
\begin{subfigure}[t]{\textwidth}
    \centering
    \includegraphics[scale=0.25]{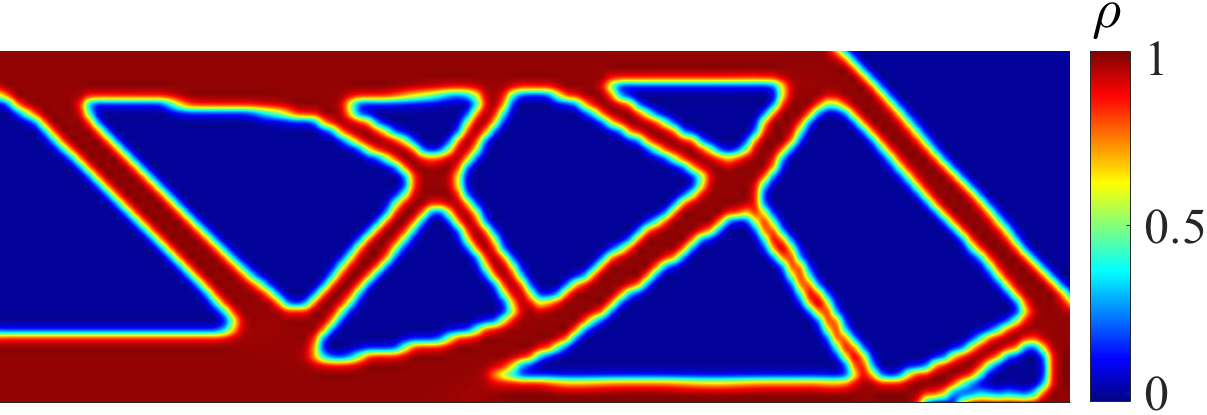}
    \caption{Final design using the BF-SVRG algorithm with $N_l=20$ and $N_h=4$.}\label{fig:bfsvrg_des1}
\end{subfigure}
\\
\begin{subfigure}[t]{\textwidth}
    \centering
    \includegraphics[scale=0.25]{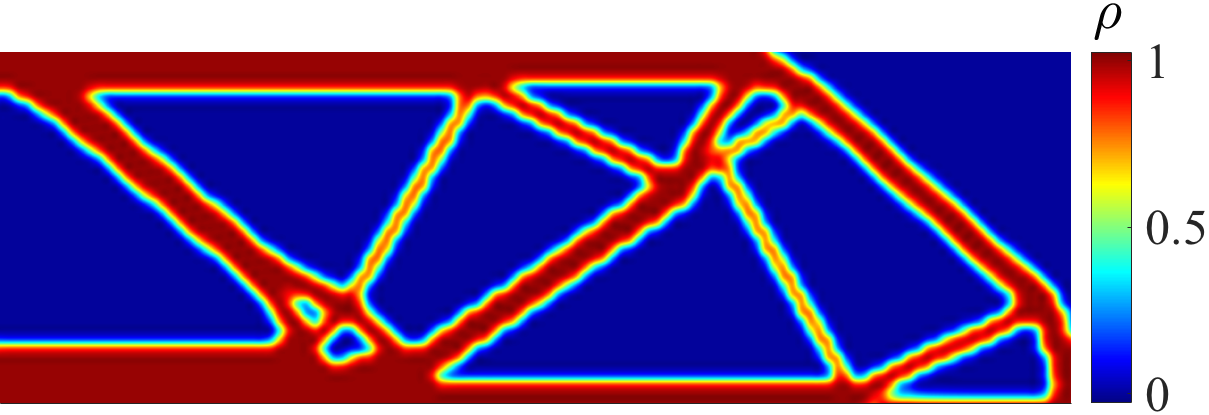}
    \caption{Final design using the BF-SVRG algorithm with $N_l=20$ and $N_h=8$.}\label{fig:bfsvrg_des2}
\end{subfigure}
\caption{Final designs for Example III (a) as obtained from the SVRG and the proposed BF-SVRG algorithms. Note that, in subfigures (a) and (b) both algorithms use the same number of gradient evaluations per outer iteration (see Algorithm \ref{alg:bfsgdcv}). However, in subfigure (c) we use a better estimate of $\boldsymbol{\alpha}$ using more low- and high-fidelity gradient solves.}\label{fig:ex3_des2}
\end{figure}
The plot of the objective in Figure \ref{fig:Obj2_ex3a} shows that SVRG diverges initially as the design changes substantially at the beginning of the optimization. One possible reason might be the poor approximation of the control variate ($\thetaa_\mathrm{prev}$ in Algorithm \ref{alg:svrg}) that leads to a poor design and large compliance values. The use of more low-fidelity gradient samples along with the calculation of an optimal $\boldsymbol{\alpha}$ in \eqref{eq:bfsvrg} improves the convergence of the solution and avoids large objective values. The BF-SVRG algorithm also leads to smaller variations in the objective and hence an improved variance reduction. Similar observations can be made from Figure \ref{fig:Mass2_ex3a} for the mass ratio. These results suggest that the use of the BF-SAG and BF-SVRG algorithms can improve the convergence of the OuU problem when compared to their single fidelity counterparts.

\begin{figure}[htb!]
\centering
\begin{subfigure}[t]{\textwidth}
    \centering
    \includegraphics[scale=0.3]{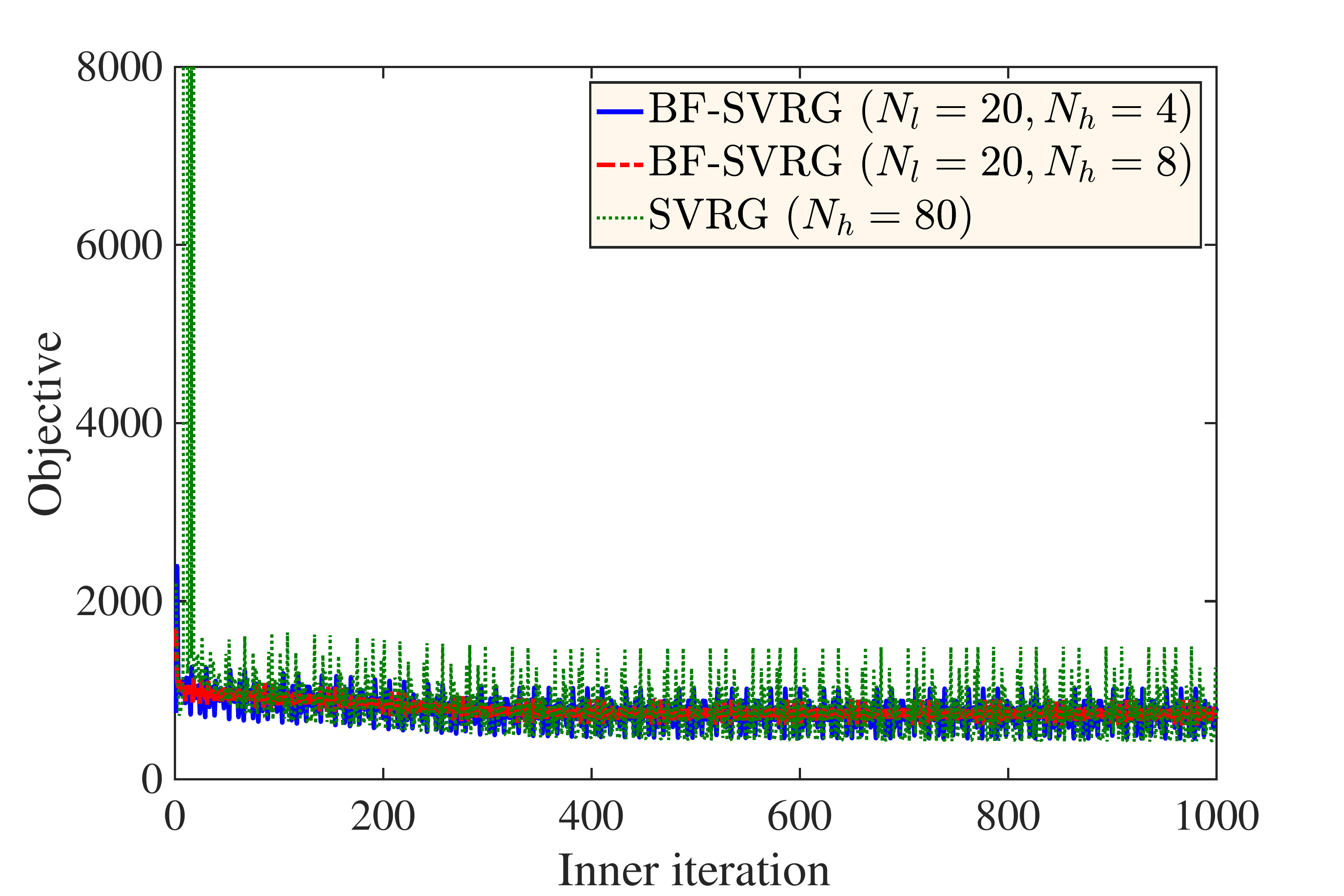}
    \caption{Objective vs. iteration.}\label{fig:Obj2_ex3a}
\end{subfigure}\\
\begin{subfigure}[t]{\textwidth}
    \centering
    \includegraphics[scale=0.3]{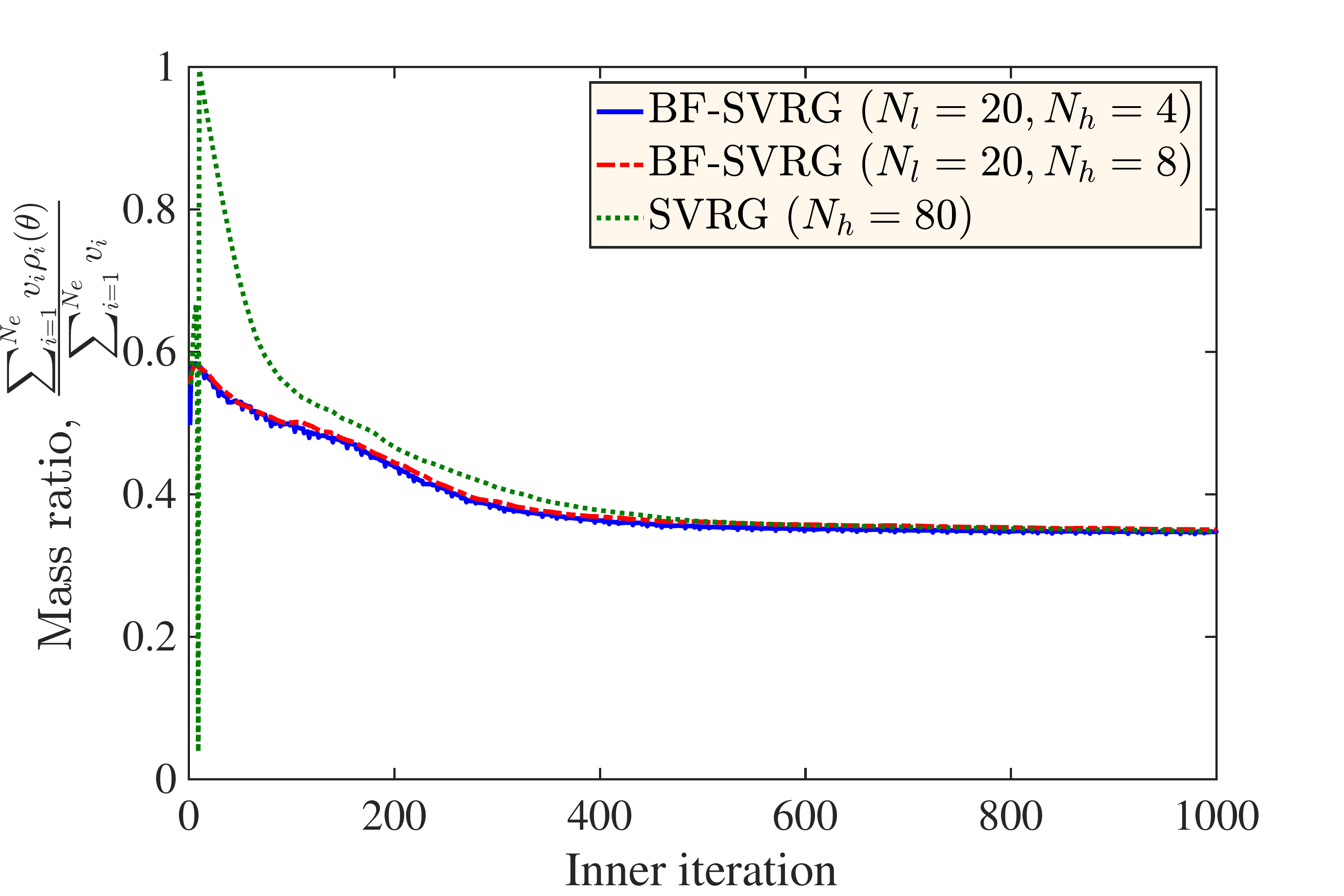}
    \caption{Mass ratio vs. iteration.}\label{fig:Mass2_ex3a}
\end{subfigure}
\caption{
Reduction of objective and mass ratio during the optimization process for two configurations of the BF-SVRG algorithm and one configuration of the SVRG algorithm for Example III (a).
Note that by using a better estimate of $\boldsymbol{\alpha}$ we can reduce the oscillations in the objective.}\label{fig:ex3_Obj2}
\end{figure}

\FloatBarrier

\subsection{Example III (b): Topology Optimization of a Beam under Uncertain Load Magnitude and Direction}

We next consider the same beam problem as in Example III (a) but add another load at a distance $L/8$ from the mid-span when only one-half of the beam is considered due to symmetry; see Figure \ref{fig:prims_schematic_b}. While the magnitude of this force, $P_0=1$ is deterministic, its direction $\phi$ relative to the beam's longitudinal axis is assumed random and given by 
\begin{equation}\label{eq:xi_alpha}
    \phi(\xi_\phi) = \pi/4 + \xi_\phi,
\end{equation}
where $\xi_\phi$ is a uniform random variable in $[-\pi/8,\pi/8]$. Hence, in this example, we have two uncertain parameters --- $\xi$ in \eqref{eq:load_unc} and $\xi_\phi$ in \eqref{eq:xi_alpha}. 

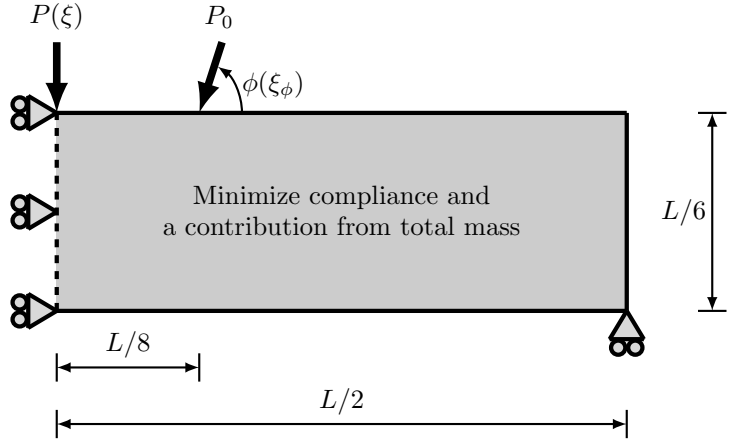
\begin{figure}[htb!]
	\centering
	\begin{tikzpicture}[scale=0.75,every node/.style={minimum size=1cm},on grid]
	\draw [ultra thick,fill=gray!40,draw=none] (0,-3.5) rectangle (10,0);
	\draw [ultra thick] (0,0) -- (10,0);
	\draw [ultra thick] (10,0) -- (10,-3.5);
	\draw [ultra thick] (10,-3.5) -- (0,-3.5);
	\draw [ultra thick,dashed] (0,-3.5) -- (0,0);
	
	\draw[ultra thick,fill=gray!30] (9.85,-4.15) circle (0.125);
	\draw[ultra thick,fill=gray!30] (10.15,-4.15) circle (0.125);
	\draw[fill=gray!30]    (10,-3.5) -- ++(0.3,-0.5) -- ++(-0.6,0) -- ++(0.3,0.5);
	\draw [ultra thick] (10,-3.5) -- (10.3,-4);
	\draw [ultra thick] (10.3,-4) -- (9.7,-4);
	\draw [ultra thick] (9.7,-4) -- (10,-3.5);
	
	\node[draw=none] at (5, -1.5)  (c)     {Minimize compliance and};
	\node[draw=none] at (5, -2)  (c)     {a contribution from total mass};
	
	\draw[-latex, line width=1mm] (0,1.25) -- (0,0);
	\node[draw = none] at (0,1.7) () {$P(\xi)$};
	
	\draw[-latex, line width=1mm] (2.9,1.25) -- (2.5,0);
	\node[draw = none] at (2.85,1.7) () {$P_0$};
	\node[draw = none] at (3.8,0.5) () {$\phi(\xi_\phi)$};
	\draw[thick,-latex] (3.25,0) arc (0:55:1);
	
	\draw[ultra thick,fill=gray!30] (-0.65,-0.15) circle (0.125);
	\draw[ultra thick,fill=gray!30] (-0.65,0.15) circle (0.125);
	\draw[fill=gray!30]    (0,0) -- ++(-0.5,0.3) -- ++(0,-0.6) -- ++(0.5,0.3);
	\draw [ultra thick] (0,0) -- (-0.5,-0.3);
	\draw [ultra thick] (-0.5,-0.3) -- (-0.5,0.3);
	\draw [ultra thick] (-0.5,0.3) -- (0,0);
	\draw[ultra thick,fill=gray!30] (-0.65,-3.35) circle (0.125);
	\draw[ultra thick,fill=gray!30] (-0.65,-3.65) circle (0.125);
	\draw[fill=gray!30]    (0,-3.5) -- ++(-0.5,0.3) -- ++(0,-0.6) -- ++(0.5,0.3);
	\draw [ultra thick] (0,-3.5) -- (-0.5,-3.8);
	\draw [ultra thick] (-0.5,-3.8) -- (-0.5,-3.2);
	\draw [ultra thick] (-0.5,-3.2) -- (0,-3.5);
	
	\draw[ultra thick,fill=gray!30] (-0.65,-1.90) circle (0.125);
	\draw[ultra thick,fill=gray!30] (-0.65,-1.60) circle (0.125);
	\draw[fill=gray!30]    (0,-1.75) -- ++(-0.5,0.3) -- ++(0,-0.6) -- ++(0.5,0.3);
	\draw [ultra thick] (0,-1.75) -- (-0.5,-2.05);
	\draw [ultra thick] (-0.5,-2.05) -- (-0.5,-1.45);
	\draw [ultra thick] (-0.5,-1.45) -- (0,-1.75);
	
	\draw[thick,latex-latex] (0,-5.5) -- (10,-5.5);
	\node[draw = none] at (5,-5.1) () {$L/2$};
	\draw[thick] (0,-5.25) -- (0,-5.75);
	\draw[thick] (10,-5.25) -- (10,-5.75);
	
	\draw[thick,latex-latex] (0,-4.5) -- (2.5,-4.5);
	\node[draw = none] at (1.25,-4.1) () {$L/8$};
	\draw[thick] (0,-4.25) -- (0,-4.75);
	\draw[thick] (2.5,-4.25) -- (2.5,-4.75);
	
	\draw[thick,latex-latex] (11.5,0) -- (11.5,-3.5);
	\node[draw = none] at (11,-1.75) () {$L/6$};
	\draw[thick] (11.25,0) -- (11.75,0);
	\draw[thick] (11.25,-3.5) -- (11.75,-3.5);
	\end{tikzpicture}
	\caption{Schematic for the beam problem in Example III (b) (design domain is shown as the shaded region). Only one-half of the beam is shown because of symmetry.}
	\label{fig:prims_schematic_b}
\end{figure}

\subsubsection{Results}
In this example, we use $\lambda = 0.25$ in the objective function (see (\ref{eq:top_obj})) and a learning rate $\eta=0.05$. The final designs obtained using the SAG and the proposed BF-SAG algorithms are shown in Figure \ref{fig:ex3b_sag_des}. 
The designs differ significantly in the mass they use. Note that, the mass contributes to the objective and larger mass increases the objective value but gives a smaller compliance. As a result, we reach different locally optimum designs. 
Figures \ref{fig:Obj1_ex3b} and \ref{fig:Mass1_ex3b} shows that the design obtained using the BF-SAG algorithm with $N_l=95$ and $N_h=5$ has the smallest objective value but uses more mass. The BF-SAG algorithm with a similar number of gradient evaluations per iteration performs slightly worse than the SAG algorithm that uses only high-fidelity gradients. Interestingly, the SAG algorithm uses more high-fidelity gradients but does not converge to a design that matches the performance of BF-SAG.
Further, in this example, the SVRG algorithm fails to converge. 
The BF-SVRG algorithm, on the other hand, produces meaningful designs as shown in Figures \ref{fig:bfsvrg1_desex3b} and \ref{fig:bfsvrg2_desex3b}. Initially, the designs undergo drastic changes -- in terms of the compliance and objective -- as can be seen in Figure \ref{fig:ex3b_Obj2}. Since the SVRG algorithm uses a control variate of the gradient based on past design parameters this control variate is poorly correlated with the gradient at the current iteration, specially during the initial iterations. We conjecture this to be the main factor in the failure of the SVRG algorithm in this example.

\begin{figure}[htb!]
\centering
\begin{subfigure}[t]{\textwidth}
    \centering
    \includegraphics[scale=0.25]{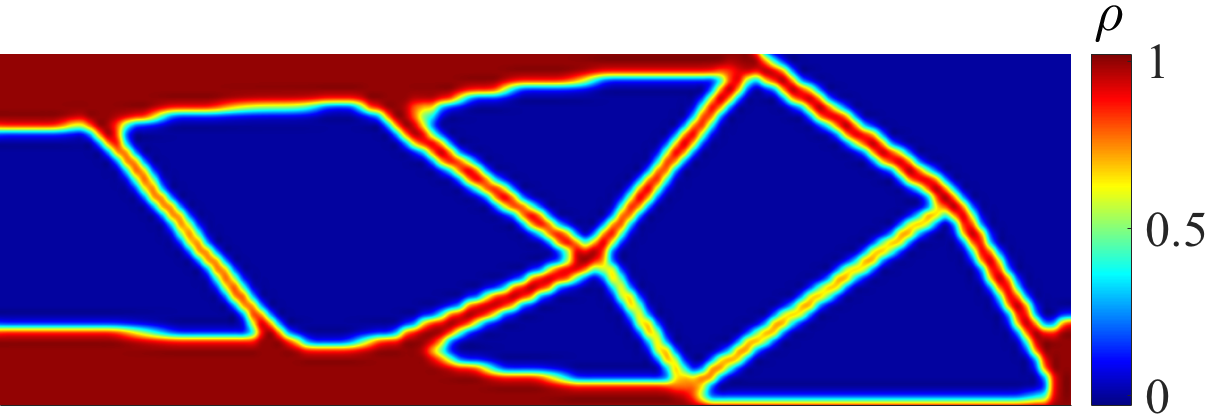}
    \caption{Final design using the SAG algorithm with $N_h=25$.}\label{fig:sag_desex3b}
\end{subfigure}\\
\begin{subfigure}[t]{\textwidth}
    \centering
    \includegraphics[scale=0.25]{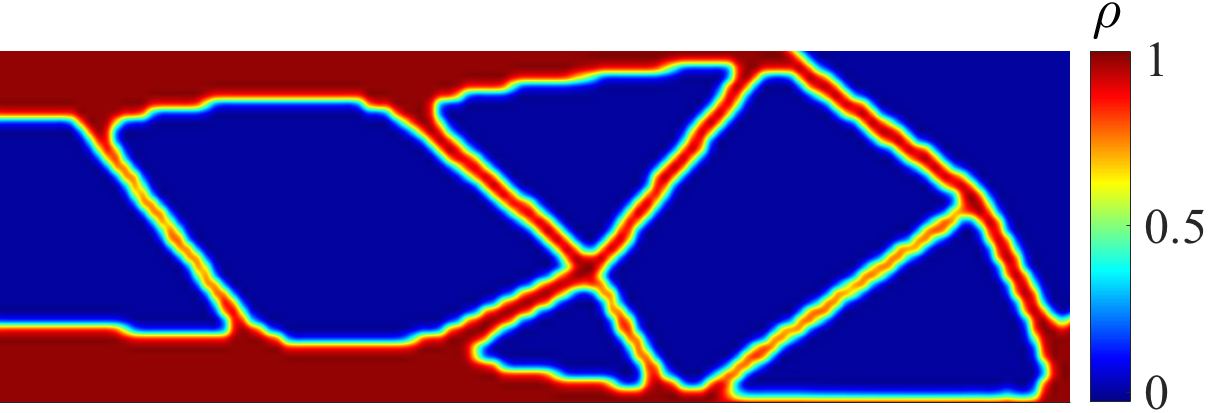}
    \caption{Final design using the BF-SAG algorithm with $N_l=20$ and $N_h=5$.}\label{fig:bfsag_desex3b}
    \end{subfigure}\\
\begin{subfigure}[t]{\textwidth}
    \centering
    \includegraphics[scale=0.25]{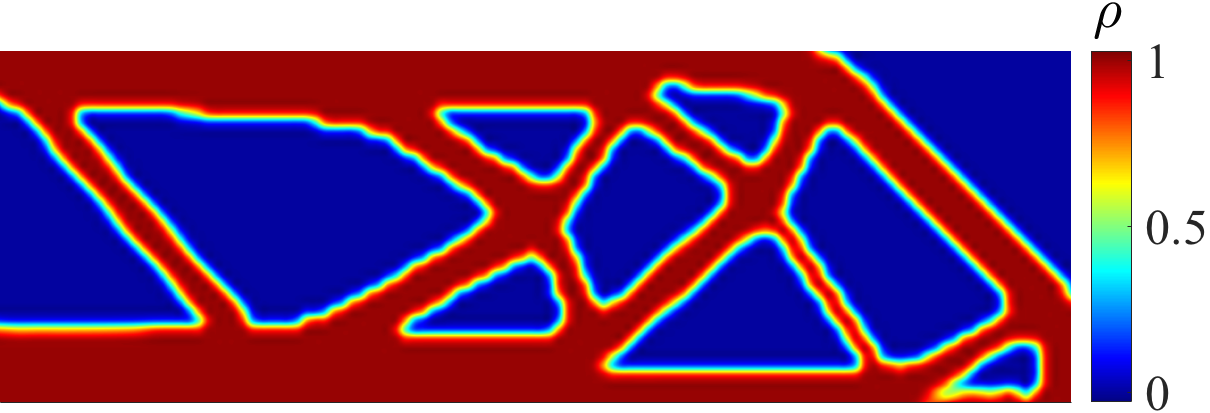}
    \caption{Final design using the BF-SAG algorithm with $N_l=95$ and $N_h=5$.}\label{fig:bfsag2_desex3b}
\end{subfigure}
\caption{Final designs for Example III (b) as obtained from the SAG and the proposed BF-SAG algorithms. Note that, both algorithms use the same number of finite element solves per iteration in subfigures (a) and (b). However, the BF-SAG algorithm only uses 5 high-fidelity solutions compared to 25 in the SAG algorithm. In subfigure (c), we use more low-fidelity solutions per iteration.}\label{fig:ex3b_sag_des}
\end{figure}
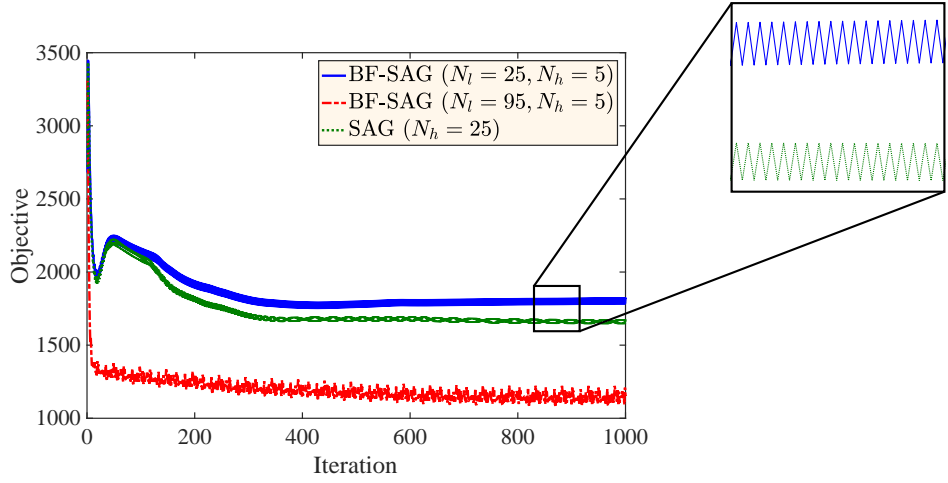
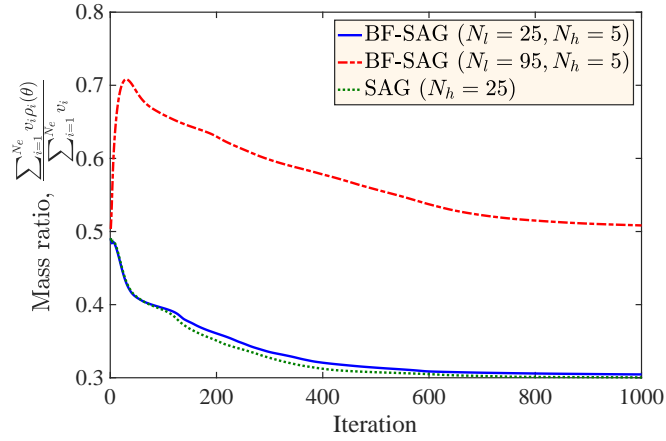
\begin{figure}[htb!]
\centering
\begin{subfigure}[t]{\textwidth}
    \centering
    \begin{tikzpicture}
\node[inner sep=0pt] (russell) at (0,0)
    {\includegraphics[scale=0.3]{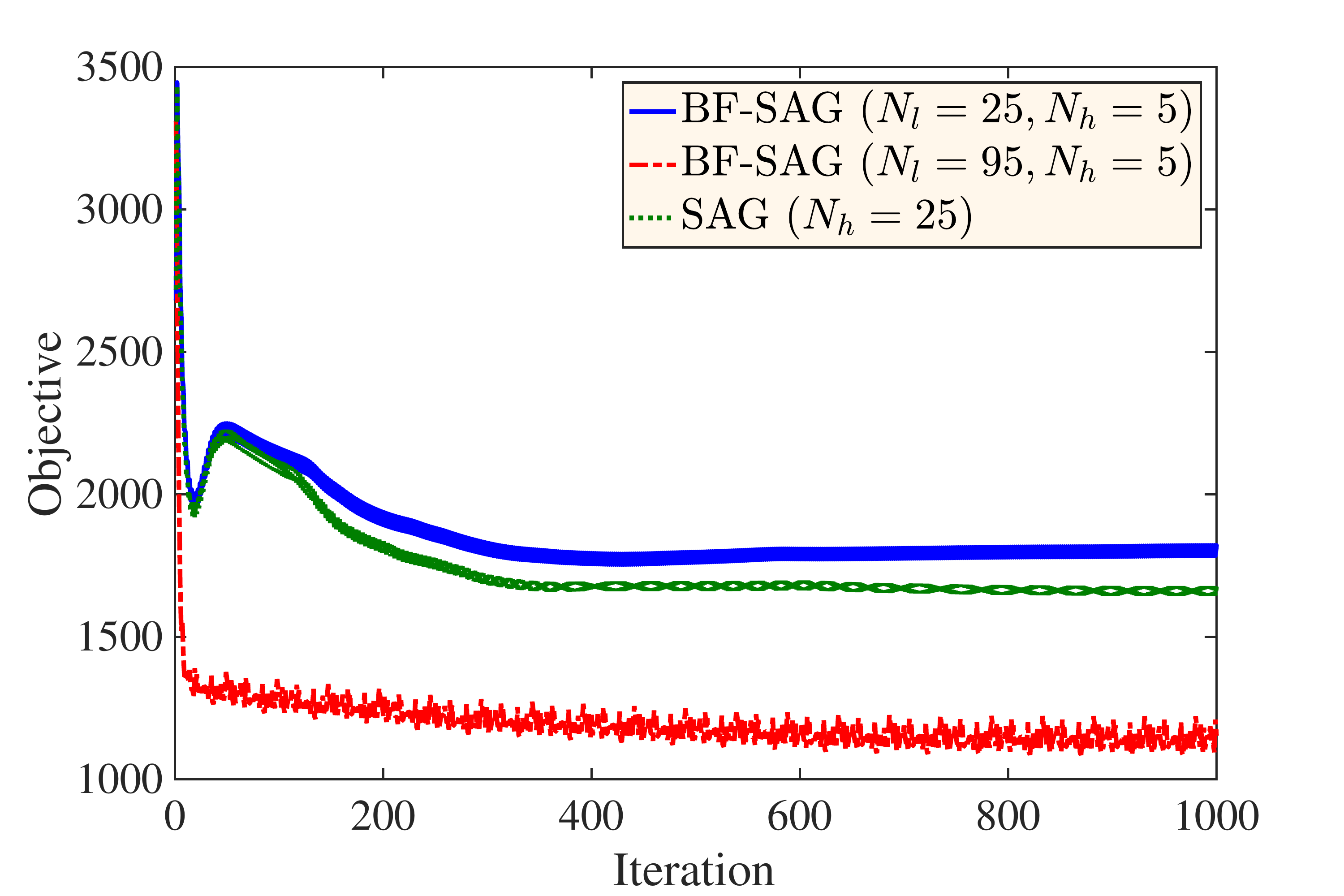}};
\node[inner sep=0pt] (whitehead) at (6.5,2)
    {\includegraphics[scale=0.2]{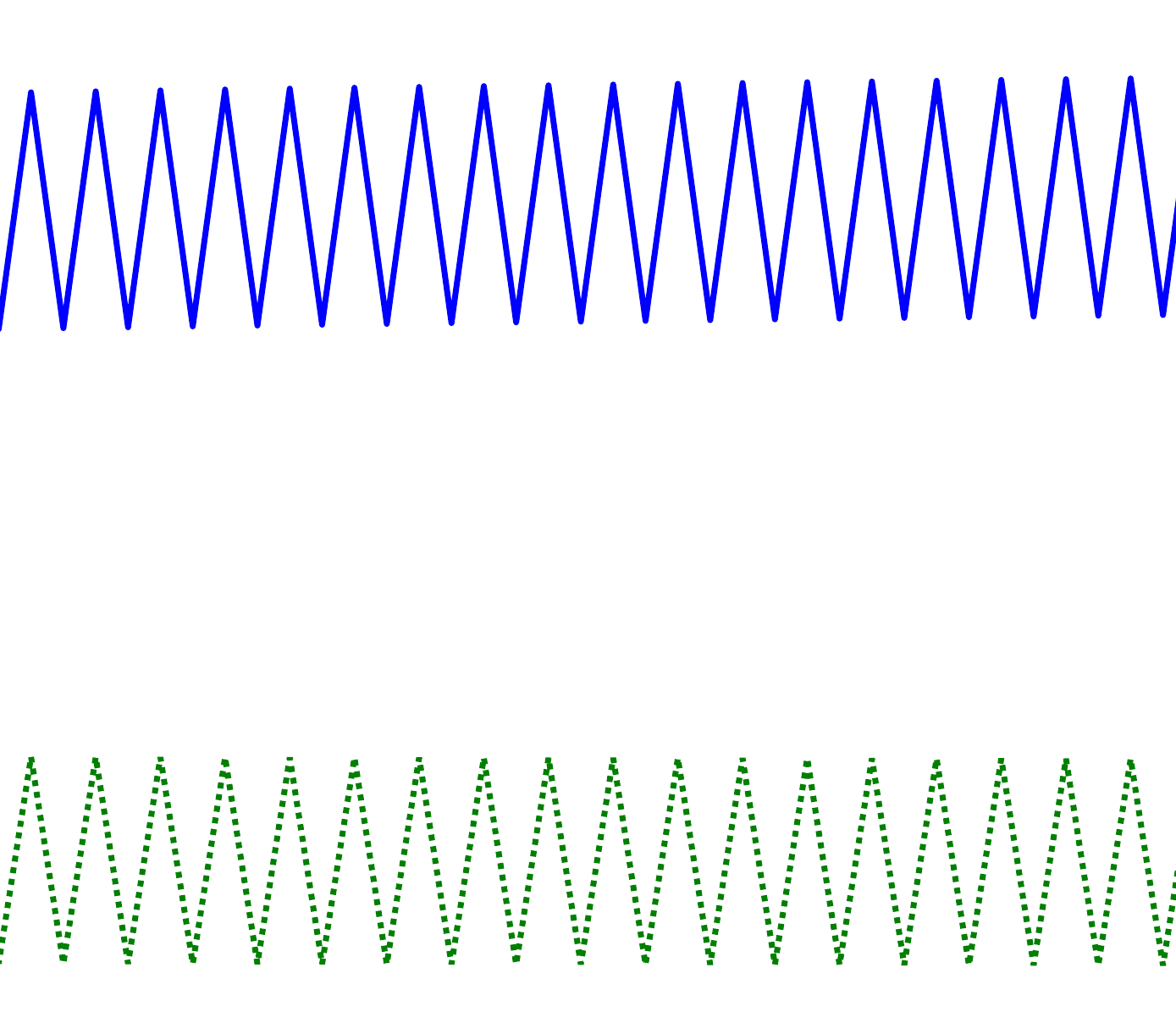}};
\draw[-,thick] (3.1,-1.1) -- (7.9,0.75);
\draw[-,thick] (2.5,-0.5) -- (5.1,3.25);
\draw[draw=black,thick] (5.1,0.75) rectangle ++(2.8,2.5); 
\draw[draw=black,thick] (2.5,-1.1) rectangle ++(0.6,0.6); 
\end{tikzpicture}
    \caption{Objective vs. iteration (Oscillations in the objective are shown in the inset figure).}\label{fig:Obj1_ex3b}
\end{subfigure}\\
\begin{subfigure}[t]{\textwidth}
    \centering
    \includegraphics[scale=0.3]{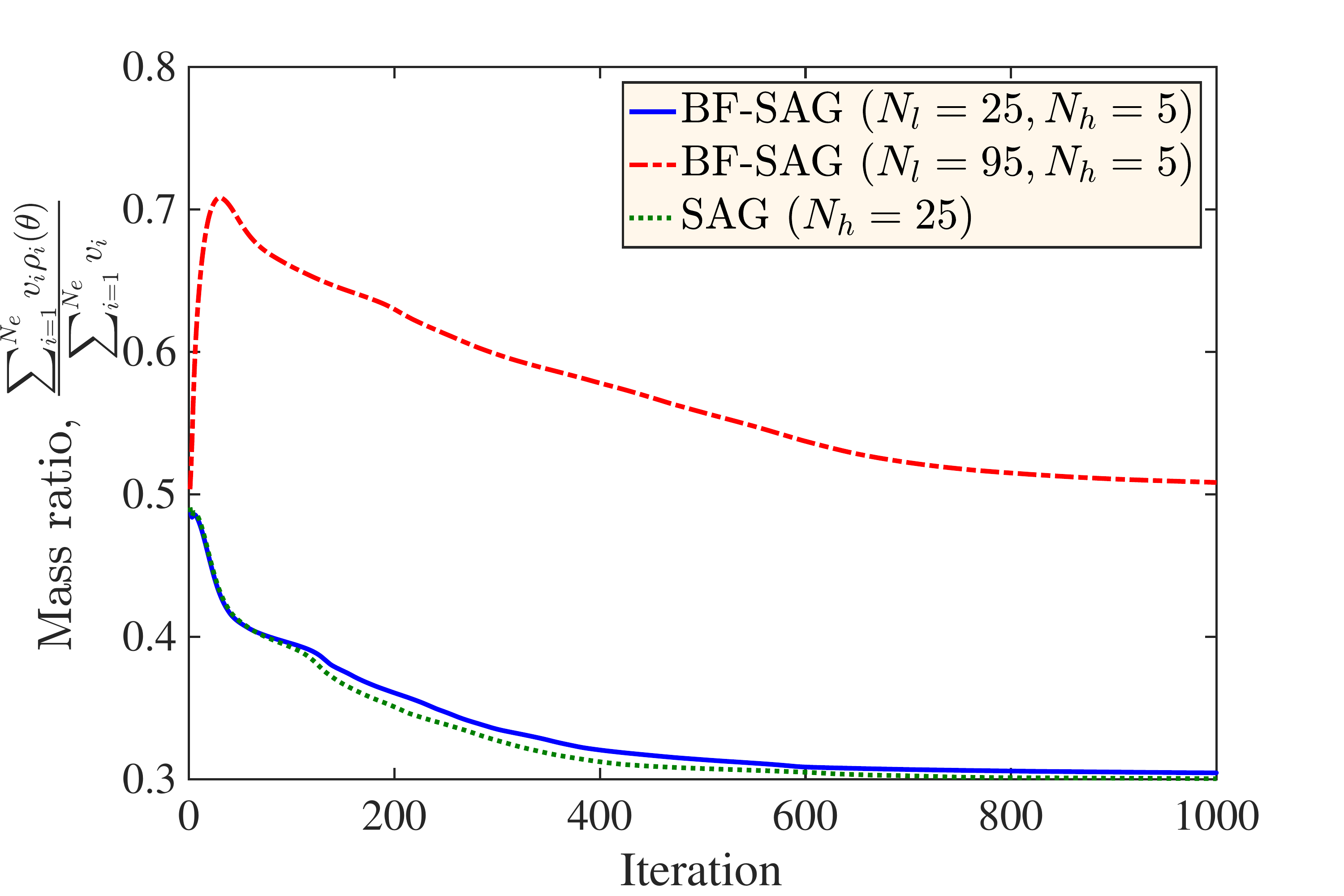}
    \caption{Mass ratio vs. iteration.}\label{fig:Mass1_ex3b}
\end{subfigure}
\caption{Reduction of objective and mass ratio during the optimization process for two configurations of the BF-SAG algorithm and one configuration of the SAG algorithm for Example III (b). Note that by using more low-fidelity models and only a handful of high-fidelity models per iteration we can reach an optimal design that has more mass but significantly smaller objective value. }\label{fig:ex3b_Obj1}
\end{figure}

\begin{figure}[htb!]
\centering
\begin{subfigure}[t]{\textwidth}
    \centering
    \includegraphics[scale=0.25]{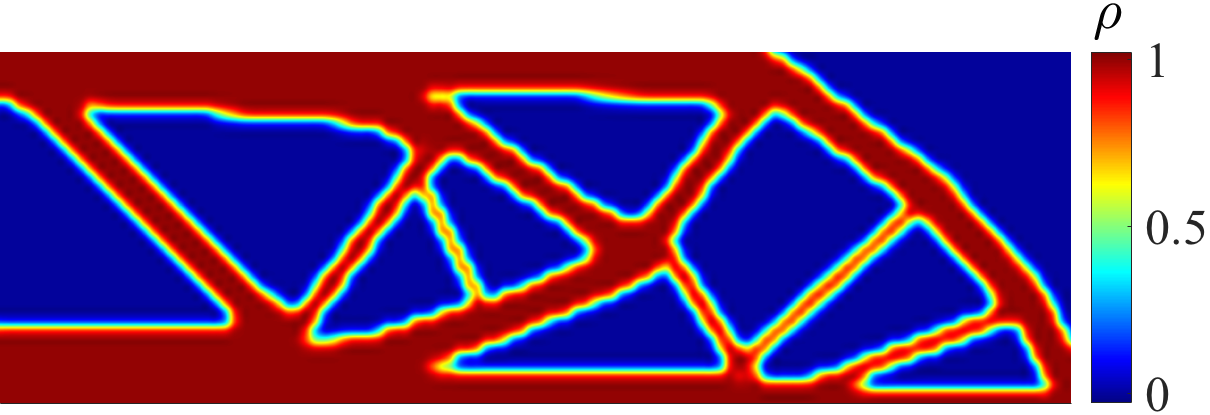}
    \caption{Failed design using the BF-SVRG algorithm with $N_l=20$ and $N_h=4$.}\label{fig:bfsvrg1_desex3b}
\end{subfigure}\\
\begin{subfigure}[t]{\textwidth}
    \centering
    \includegraphics[scale=0.25]{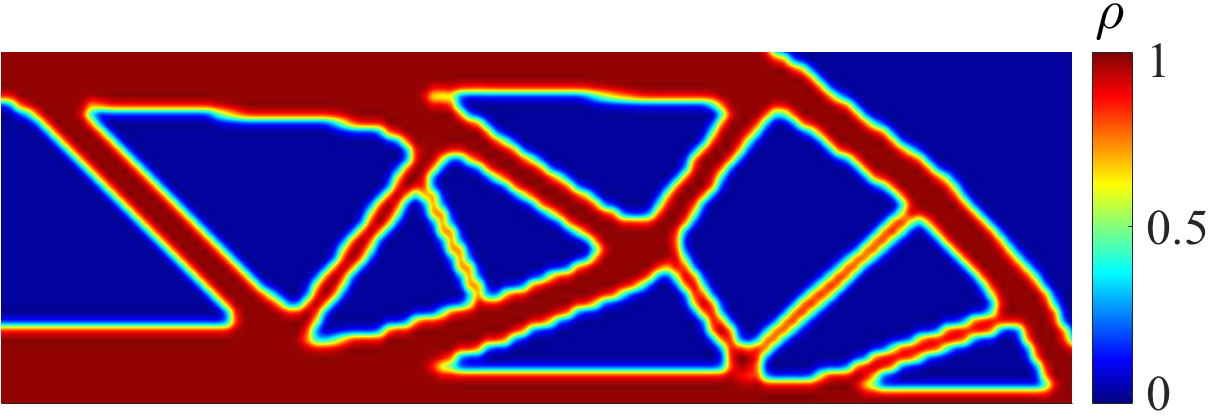}
    \caption{Final design using the BF-SVRG algorithm with $N_l=20$ and $N_h=8$.}\label{fig:bfsvrg2_desex3b}
\end{subfigure}
\caption{Final designs for Example III (b) as obtained from the proposed BF-SVRG algorithm. Note that, the SVRG algorithm fails to produce a design in this example. }\label{fig:ex3b_svrg_des}
\end{figure}

\begin{figure}[htb!]
\centering
\begin{subfigure}[t]{\textwidth}
    \centering
    \includegraphics[scale=0.3]{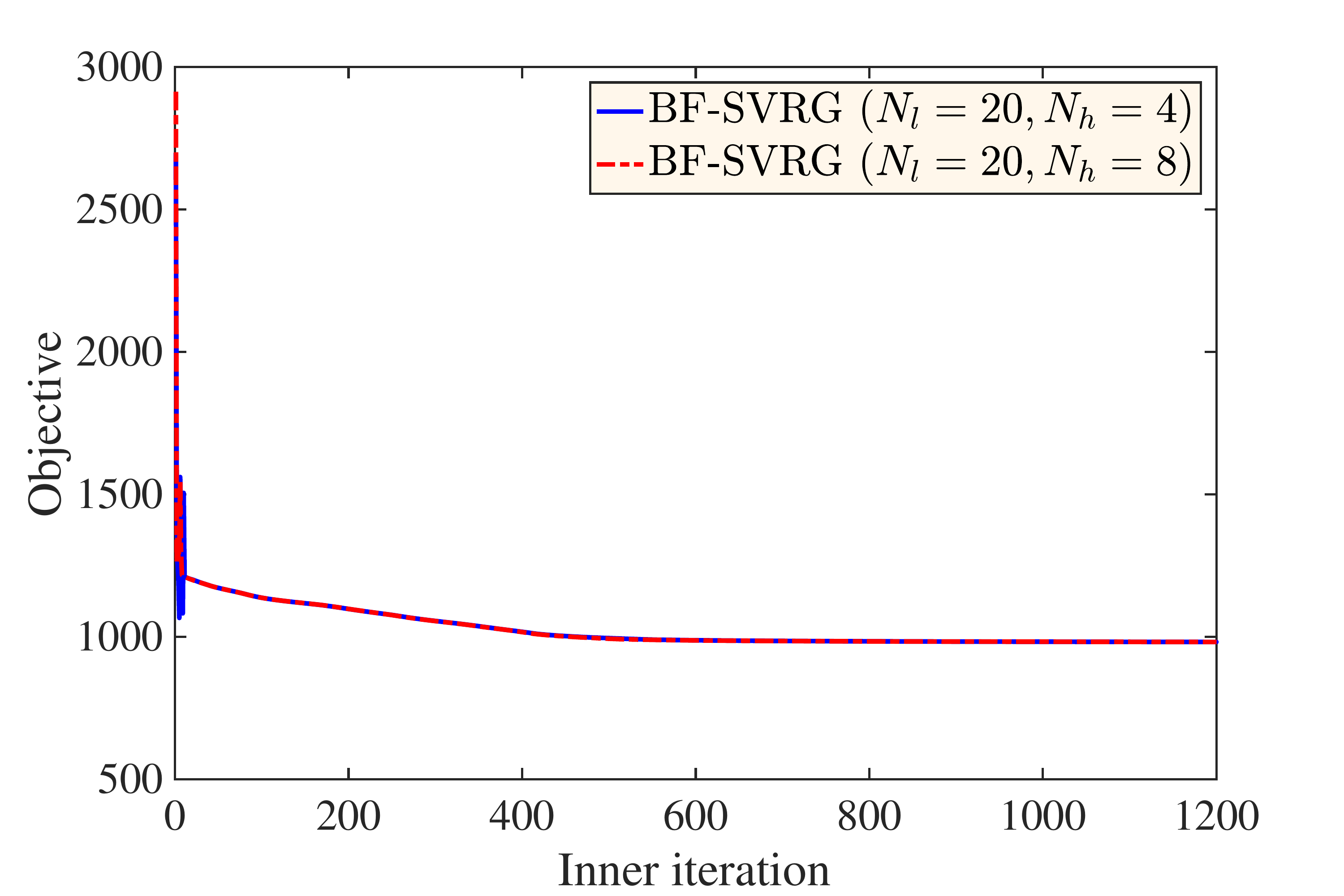}
    \caption{Objective vs. iteration.}\label{fig:compl2_ex3b}
\end{subfigure}\\
\begin{subfigure}[t]{\textwidth}
    \centering
    \includegraphics[scale=0.3]{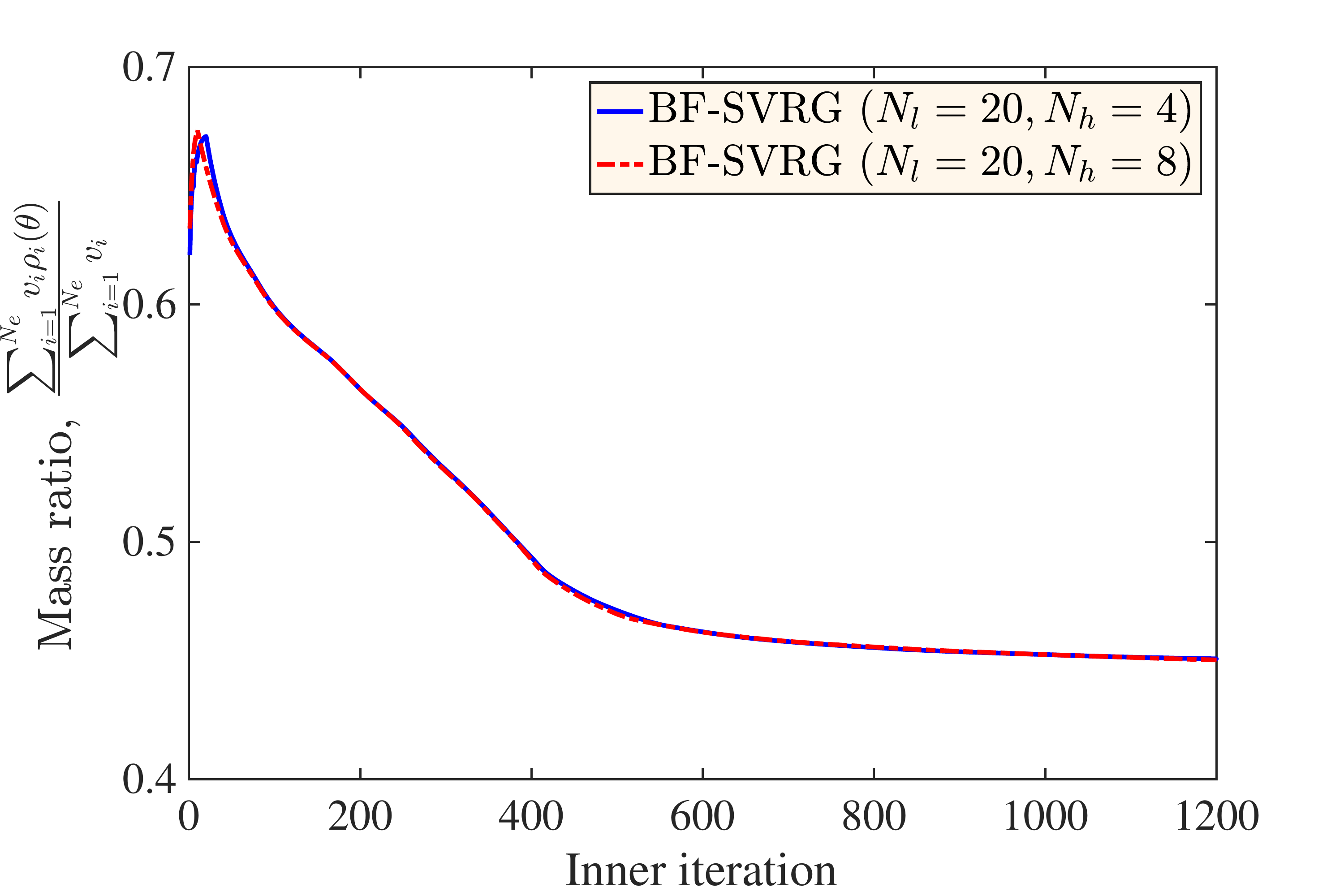}
    \caption{Mass ratio vs. iteration.}\label{fig:Obj2_ex3b}
\end{subfigure}
\caption{Reduction of compliance and objective during the optimization process for two configurations of the BF-SVRG algorithm for Example III (b). 
}\label{fig:ex3b_Obj2}
\end{figure}

\FloatBarrier
\subsection{Example III (c): Design of Beam under Uncertain Load and Material Property}
In this example, we consider the same beam problem as in Example III (a) with uncertainty in the load magnitude (see \eqref{eq:load_unc}). We further assume that the elastic modulus of the material $E_0$ in (\ref{eq:simp}) is uncertain and modeled by a lognormal random field,
\begin{equation}
    E_0(x_1,x_2)=\exp[z(x_1,x_2)],
\end{equation}
 where $z(x_1,x_2)$ is a zero-mean Gaussian field with a covariance function
\begin{equation}\label{eq:correl}
\Exp[z(x_1,x_2)z(y_1,y_2)]=\sigma^2\exp\left(-\frac{|x_1-y_1|}{l_1}-\frac{|x_2-y_2|}{l_2}\right).
\end{equation}
Here, $l_1=l_2=L/40$ and $\sigma=2$ are used. The random field $z(x_1,x_2)$ is expressed using a Karhunen-Lo\'eve expansion truncated at $N_{\max}$th term as follows,
\begin{equation}\label{eq:kl}
z(x_1,x_2)=\sum_{i=1}^{N_{\max}}  \sqrt{\lambda_i} \xi_i \psi_i(x_1,x_2)
\end{equation}
where $\lambda_i$ are eigenvalues and $\psi_i(x_1,x_2)$ are eigenfunctions of the covariance function (\ref{eq:correl}); $\xi_i$ are independent standard normal random variables; $N_{\max}=100$ is selected to capture 99.92\% of total variance of $z$. Figure \ref{fig:E0_realizations} shows three realizations of $E_0$ over the design domain. 
The uncertainty in the load magnitude is assumed same as in \eqref{eq:load_unc}.
Hence, the dimension of the uncertain parameter vector $\Ym$ is 101 in this example, whereas dimension of the optimization variable vector $\ppm$ is 4800 as before.  
\begin{figure}
    \centering
    \includegraphics[scale=0.15]{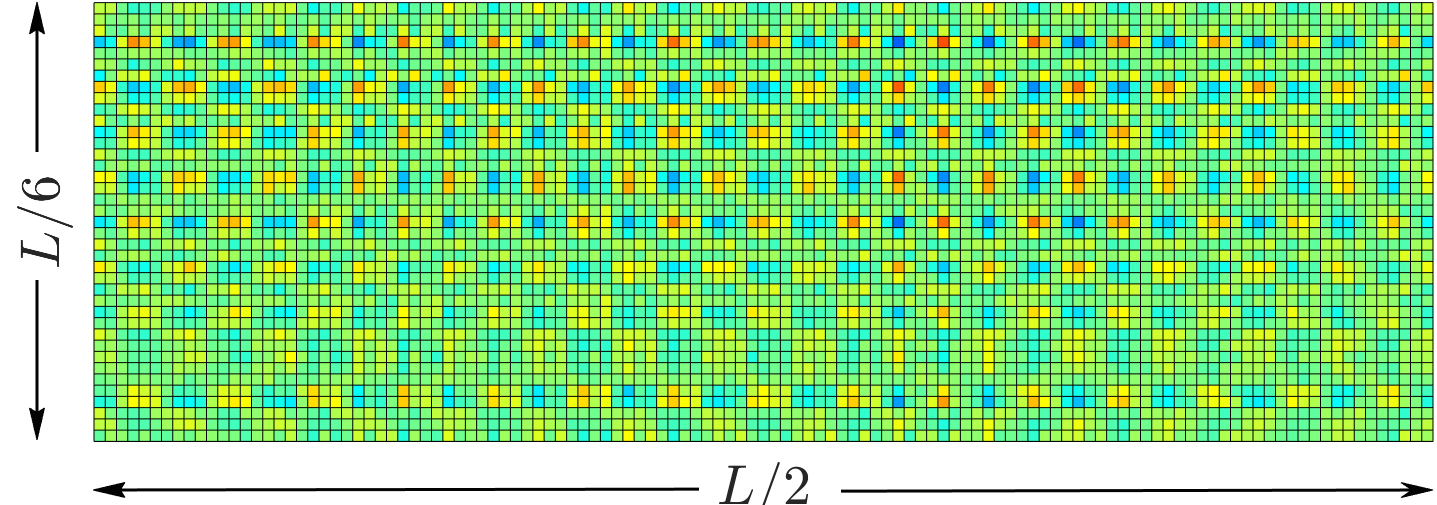}~\includegraphics[scale=0.15]{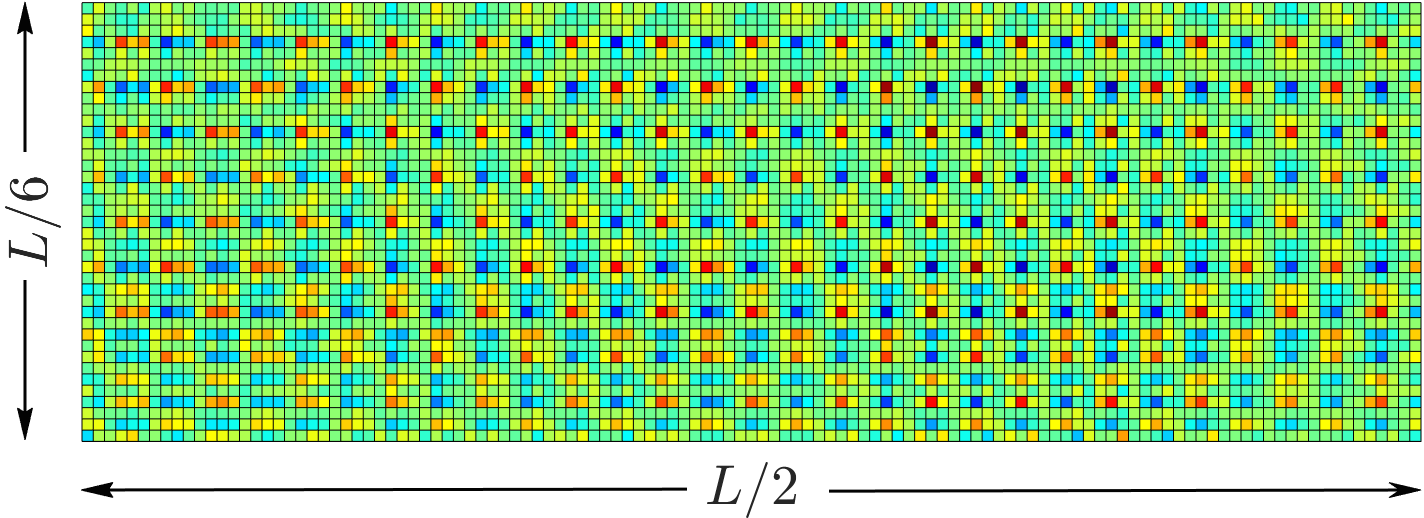}\\~\\
    \includegraphics[scale=0.15]{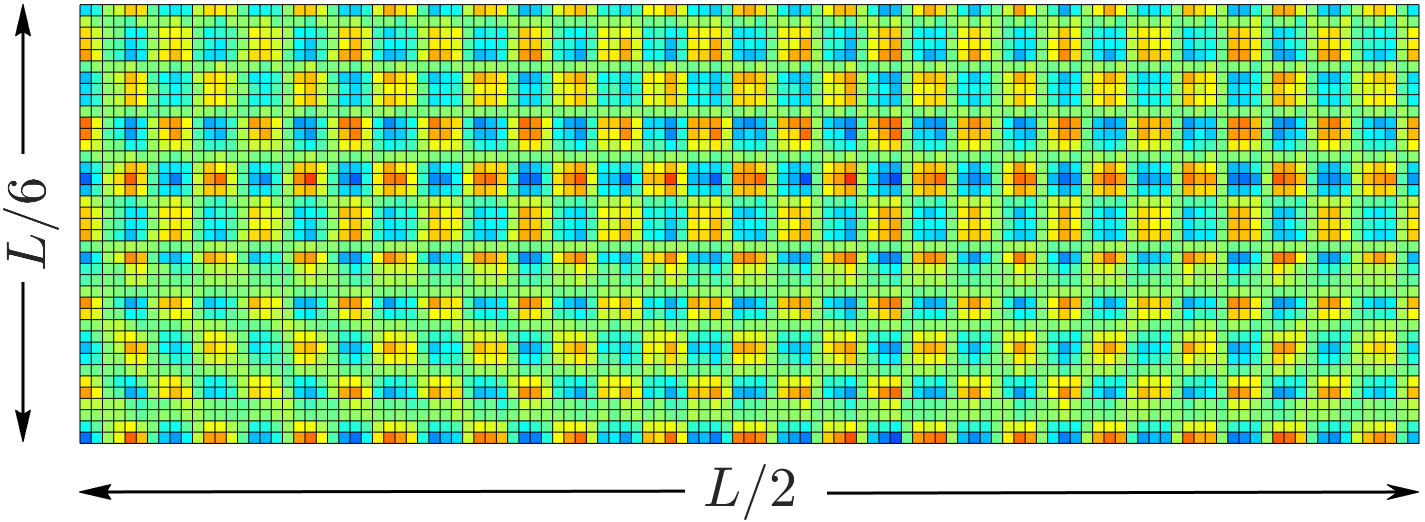}\\~\\
    \includegraphics[scale=0.15]{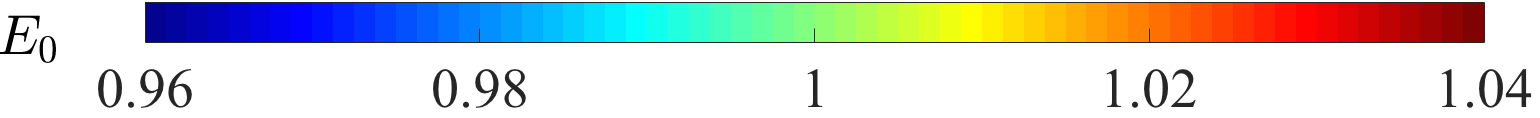}
    \caption{Three realizations of $E_0$ over the design domain in Example III (c).}
    \label{fig:E0_realizations}
\end{figure}

\subsubsection{Results}
We study the proposed algorithms using a learning rate $\eta=0.05$ and $\lambda = 0.25$ in (\ref{eq:top_obj}). The final designs from the SAG and BF-SAG algorithms are shown in Figure \ref{fig:ex3c_des1}. 
It can be noted from the figure that if we use more low-fidelity gradient evaluations per iteration, we obtain a design (Figure \ref{fig:bfsag_des3cb}) that has a smaller number of members but reaches the performance of the SAG algorithm.
This is further evident from the plots of objective and mass ratio for these three cases as shown in Figure \ref{fig:ex3c_Obj1}. Although the mass of the structure obtained from the SAG algorithm is smaller than the other two designs the BF-SAG algorithm with $N_l=95$ and $N_h=5$ produces smaller variation in the objective and mass ratio.  
In this example, the SVRG and BF-SVRG algorithms both produce meaningful designs as shown in Figure \ref{fig:ex3c_des2}. Interestingly, the designs are different visually. However, as displayed in Figure \ref{fig:ex3c_Obj2}, the BF-SVRG algorithm with more low-fidelity gradient evaluations produces the design with smallest variations in the objective.
\begin{figure}[htb!]
\centering
\begin{subfigure}[t]{\textwidth}
    \centering
    \includegraphics[scale=0.25]{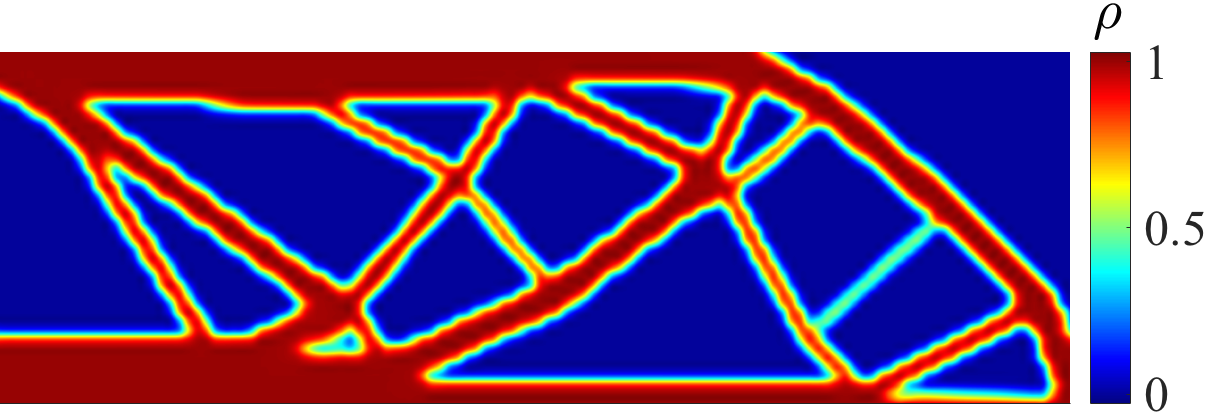}
    \caption{Final design using the SAG algorithm with $N_h=25$.}\label{fig:sag_des3c}
\end{subfigure}\\
\begin{subfigure}[t]{\textwidth}
    \centering
    \includegraphics[scale=0.25]{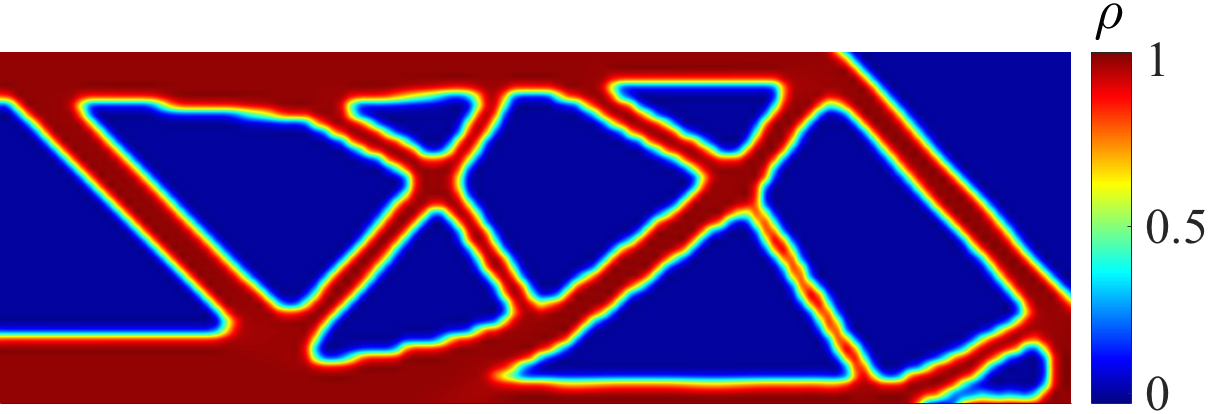}
    \caption{Final design using the BF-SAG algorithm with $N_l=20$ and $N_h=5$.}\label{fig:bfsag_des3ca}
\end{subfigure}
\\
\begin{subfigure}[t]{\textwidth}
    \centering
    \includegraphics[scale=0.25]{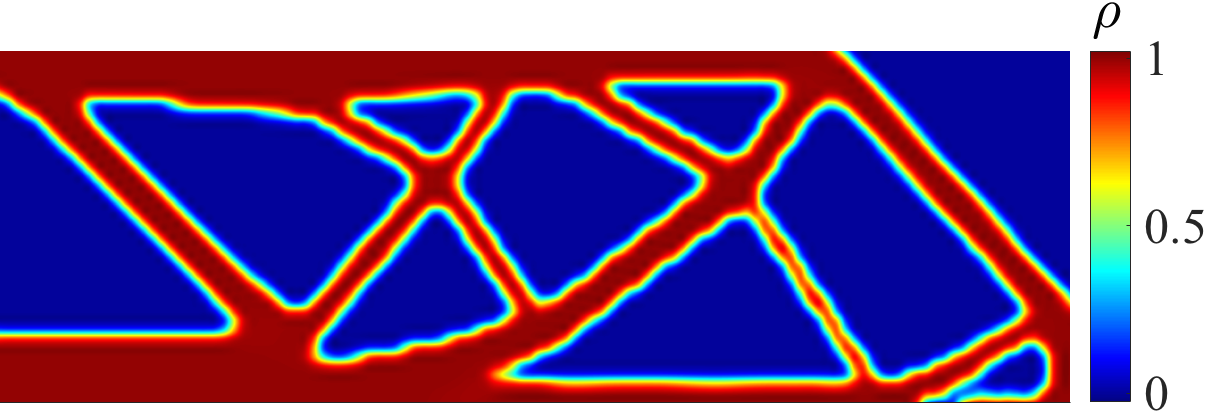}
    \caption{Final design using the BF-SAG algorithm with $N_l=95$ and $N_h=5$.}\label{fig:bfsag_des3cb}
\end{subfigure}
\caption{Final designs for Example III (c) as obtained from the SAG and the proposed BF-SAG algorithms. Note that, in (a) and (b), both algorithms use the same number of gradient evaluations per iteration. However, BF-SAG only uses 5 high-fidelity gradients compared to 25 high-fidelity gradients in the SAG algorithms.}\label{fig:ex3c_des1}
\end{figure}

\begin{figure}[htb!]
\centering
\begin{subfigure}[t]{\textwidth}
    \centering
    \includegraphics[scale=0.3]{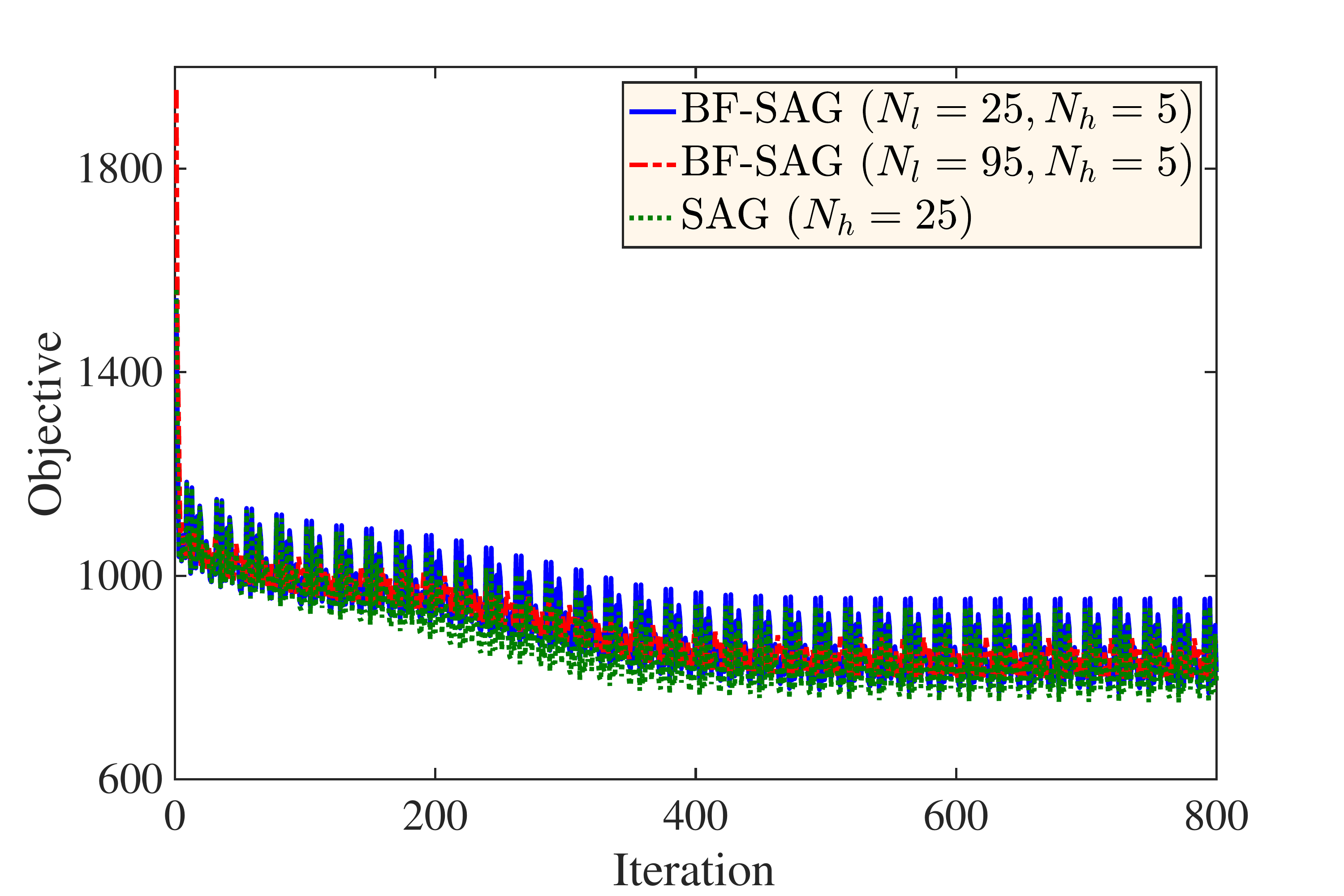}
    \caption{Objective vs. iteration.}\label{fig:Obj}
\end{subfigure}\\
\begin{subfigure}[t]{\textwidth}
    \centering
    \includegraphics[scale=0.3]{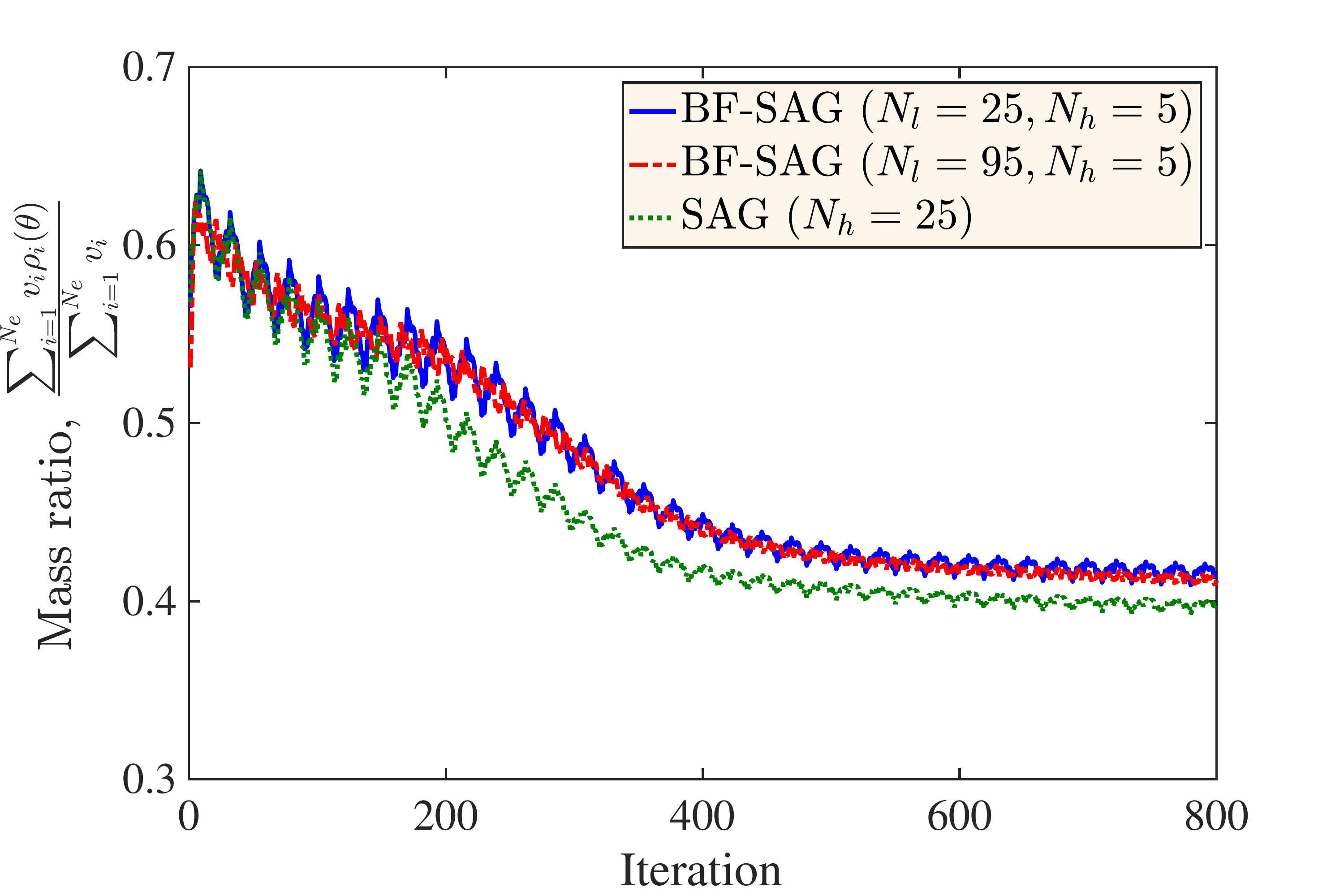}
    \caption{Mass ratio vs. iteration.}\label{fig:Mass}
\end{subfigure}
\caption{Reduction of compliance and objective during the optimization process for two configurations of the BF-SAG algorithm and one configuration of the SAG algorithm in Example III (c). Note that the faster convergence to the optimum can be achieved by using more low-fidelity models and only a handful of high-fidelity models per iteration.}\label{fig:ex3c_Obj1}
\end{figure}

\begin{figure}[htb!]
\centering
\begin{subfigure}[t]{\textwidth}
    \centering
    \includegraphics[scale=0.25]{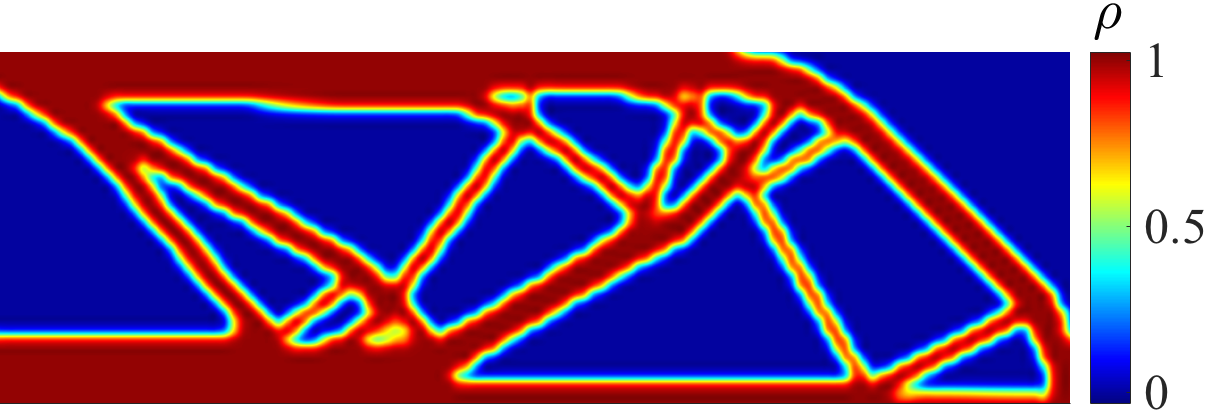}
    \caption{Final design using the SVRG algorithm with $N_h=80$.}\label{fig:svrg_des_rf}
\end{subfigure}\\
\begin{subfigure}[t]{\textwidth}
    \centering
    \includegraphics[scale=0.25]{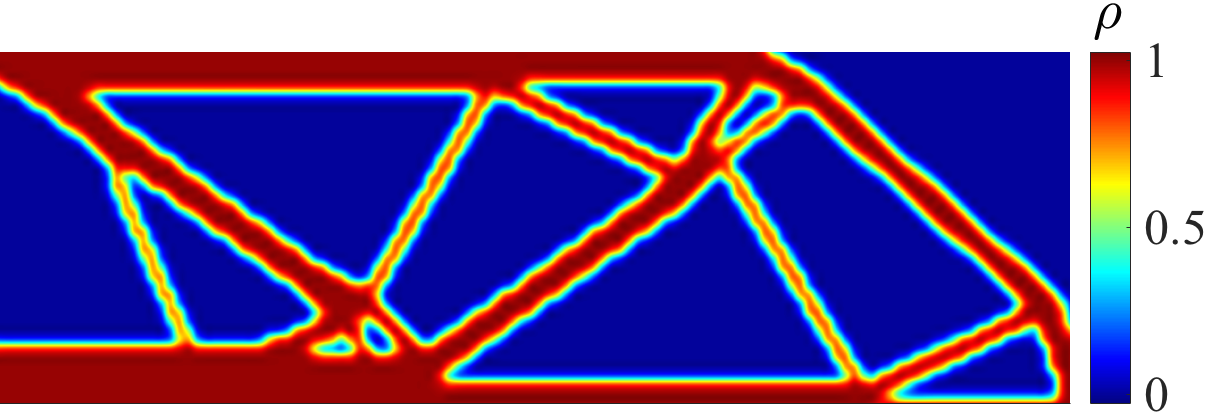}
    \caption{Final design using the BF-SVRG algorithm with $N_l=20$ and $N_h=4$.}\label{fig:bfsvrg_des_rf1}
\end{subfigure}
\\
\begin{subfigure}[t]{\textwidth}
    \centering
    \includegraphics[scale=0.25]{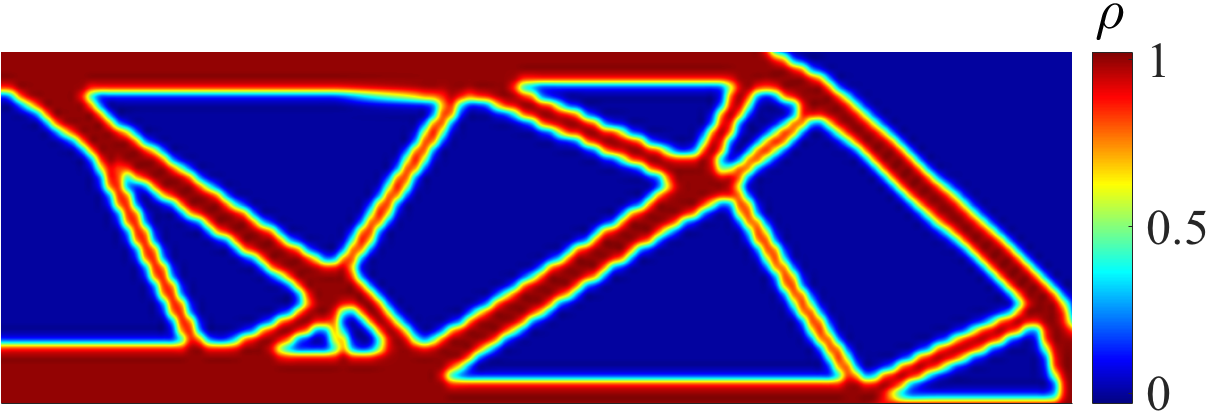}
    \caption{Final design using the BF-SVRG algorithm with $N_l=20$ and $N_h=8$.}\label{fig:bfsvrg_des_rf2}
\end{subfigure}
\caption{Final designs for Example III (c) as obtained from the SVRG and the proposed BF-SVRG algorithms. Note that, in subfigures (a) and (b) both algorithms use the same number of gradient evaluations per outer iteration (see Algorithm \ref{alg:bfsgdcv}). }\label{fig:ex3c_des2}
\end{figure}

\begin{figure}[htb!]
\centering
\begin{subfigure}[t]{\textwidth}
    \centering
    \includegraphics[scale=0.3]{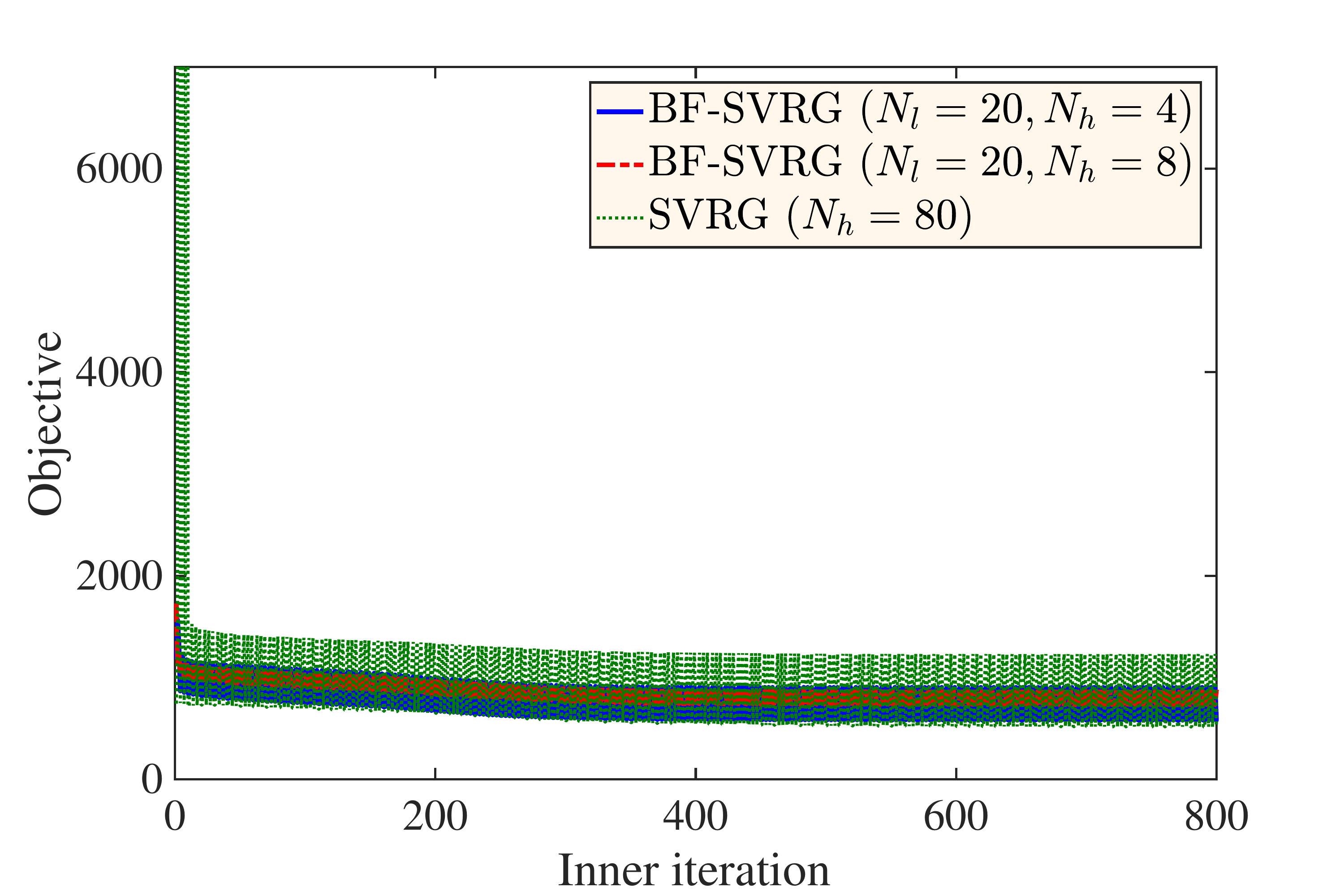}
    \caption{Objective vs. iteration.}\label{fig:Obj_rf}
\end{subfigure}\\
\begin{subfigure}[t]{\textwidth}
    \centering
    \includegraphics[scale=0.3]{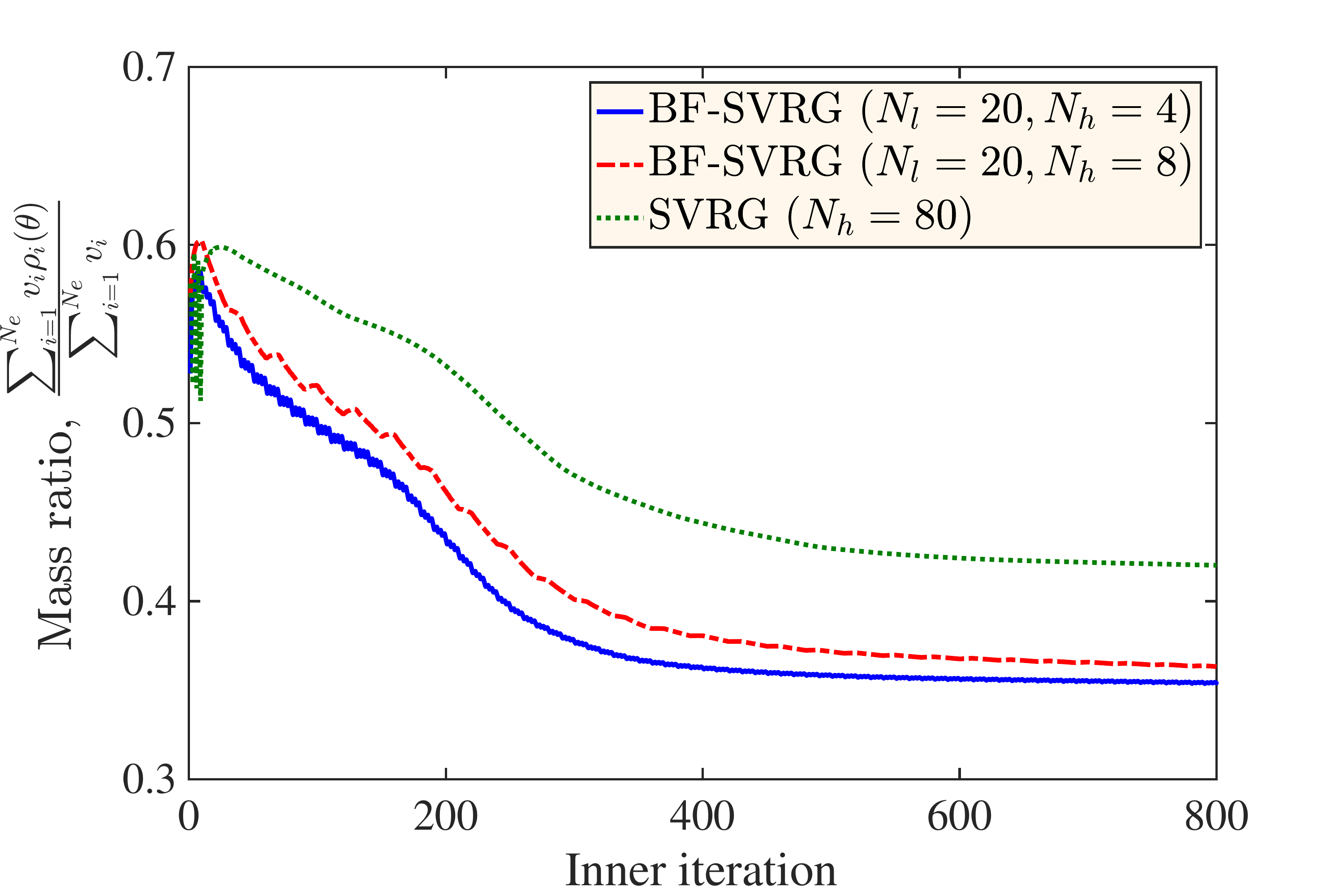}
    \caption{Mass ratio vs. iteration.}\label{fig:Mass_rf}
\end{subfigure}
\caption{Reduction of compliance and objective during the optimization process for two configurations of the BF-SVRG algorithm and one configuration of the SVRG algorithm in Example III (c). Note that by using a better estimate of $\boldsymbol{\alpha}$ we can reduce the oscillations in the objective.
}\label{fig:ex3c_Obj2}
\end{figure}


\section{Conclusions}
In the presence of uncertainty, the cost of design optimization of structures increases many folds. Methods like polynomial chaos or stochastic collocation help in reducing the cost but in the presence of high-dimensional uncertainty their costs increase rapidly as well. 
In this paper, to alleviate this computational burden of  OuU, we propose a bi-fidelity approach with stochastic gradient descent type methods, where most of the gradients are estimated using a low-fidelity model. The gradients are then incorporated into two distinct stochastic gradient descent algorithms. In the first algorithm, we use an average of the gradients, where most of them are updated using the low-fidelity model. In the second algorithm, we use a control variate based on gradients calculated using the low-fidelity model to reduce the variance in the stochastic gradients. 
Linear convergence of these proposed algorithms in ideal conditions are proved.
The efficacy of these algorithms are shown using three numerical examples. After studying the proposed algorithms with a conceptual problem we optimize the shape of a hole in a square plate to minimize the maximum stress in the plate.
In the third example, we apply the proposed algorithms to a topology optimization problem involving uncertainties in load and material properties with the number of uncertain parameters reaching 101. These examples show that using the proposed algorithms we successfully leverage a low-fidelity model to reduce the computational cost of the optimization. In future studies, the proposed bi-fidelity algorithms will be applied to multi-physics optimization problem, where the computational savings will be even more pronounced.

\section{\texorpdfstring{Acknowledgments}{Acknowledgments}}
\label{sec:ack}

The authors acknowledge the support of the Defense Advanced Research Projects Agency's (DARPA) TRADES project under agreement HR0011-17-2-0022. Any opinions, findings, and conclusions or recommendations expressed in this material are those of the authors and do not necessarily reflect the views of the DARPA. The authors would also like to thank Prof. John Evans (CU Boulder) for fruitful discussions regarding the content of this manuscript.

\appendix
\section{Proof of Theorem \ref{theorem1}}\label{proof_thm1}
Assume $\ppm_{k}$ is the vector of optimization parameters after $k$ iterations of algorithm \ref{alg:bfsag}. $\hb_\mathrm{low}(\ppm_k)$ and $\hb_\mathrm{high}(\ppm_k)$ are gradients of the objective with respect to $\ppm_k$ using the low- and high-fidelity models, respectively.
Under the assumption of strong convexity\footnote{A function $J(\ppm)$ is strongly convex with a constant $\mu$ if $J(\ppm)-\frac{\mu}{2}\lVert\ppm\rVert^2$ is convex.}
 of low-fidelity and high-fidelity objectives, 
\begin{equation}
\begin{split}
        (\ppm_k - \ppm^*)^T\hb_\mathrm{low}(\ppm_k) &\ge \mu_\mathrm{low} \lVert \ppm_k - \ppm^* \rVert^2,\\
        (\ppm_k - \ppm^*)^T\hb_\mathrm{high}(\ppm_k) &\ge \mu_\mathrm{high} \lVert \ppm_k - \ppm^* \rVert^2,\\
\end{split}
\end{equation}
where $\mu_\mathrm{low}$ and $\mu_\mathrm{high}$ are constants.
Similarly, if the low- and high-fidelity gradients are Lipschtiz continuous, 
\begin{equation}
\begin{split}
        \lVert \hb_\mathrm{low}(\ppm_k) \rVert ^2 & \leq L^2_\mathrm{low} \lVert \ppm_k - \ppm^* \rVert^2,\\
        \lVert \hb_\mathrm{high}(\ppm_k) \rVert ^2 & \leq L^2_\mathrm{high} \lVert \ppm_k - \ppm^* \rVert^2,\\
\end{split}
\end{equation}
where $L_\mathrm{low},L_\mathrm{high}$ are the Lipschitz constants for high- and low-fidelity gradients, respectively.
The parameters are updated in Algorithm 4 using
\begin{equation}
    \ppm_{k+1} = \ppm_{k} - \eta \widehat{\hb}_k.
\end{equation}
The expected value of the gradient $\widehat{\hb}_k$ at iteration $k$ is
\begin{equation}
\Exp[\widehat{\hb}_k|\ppm_{k}] =  p_l \hb_\mathrm{low}(\ppm_{k}) + p_h \hb_\mathrm{high}(\ppm_{k}) + (1- p_l -p_h) \dm_{k-1},
\end{equation}
where $p_l = N_l/N$ and $p_h = N_h/N$.

Next, we evaluate the following expectation
\begin{equation}
    \begin{split}
        \Exp[\lVert \ppm_{k+1} - \ppm^*\rVert^2\lvert \ppm_{k}] & = \Exp[\lVert \ppm_{k} - \ppm^* -\eta \widehat{\hb}_k\rVert^2\lvert \ppm_{k}]\\
        & = \lVert \ppm_{k} - \ppm^* \rVert ^2 - 2 \eta (\ppm_{k} - \ppm^*)^T\Exp[\widehat{\hb}_k| \ppm_{k}] + \eta^2 \Exp[\lVert \widehat{\hb}_k\rVert^2| \ppm_{k} ]\\
        & \leq \lVert \ppm_{k} - \ppm^* \rVert ^2 - 2 \eta (\ppm_{k} - \ppm^*)^T\Exp[\widehat{\hb}_k| \ppm_{k}] + \eta^2 L^2 \lVert \ppm_{k} - \ppm^* \rVert ^2,
    \end{split}
\end{equation}
where $L^2 = \max \left\{p_l(1-p_l-p_h)^{(k-j)}\frac{L_\mathrm{low}^2\lVert \ppm_j - \ppm^* \rVert^2}{\lVert \ppm_k - \ppm^* \rVert^2},p_h(1-p_l-p_h)^{(k-j)}\frac{L_\mathrm{high}^2\lVert \ppm_j - \ppm^* \rVert^2}{\lVert \ppm_k - \ppm^* \rVert^2}\right\}$ for $j=1,\dots,k$. 
Using the strong convexity property of $J(\ppm)$, 
\begin{equation}
    \begin{split}\label{eq:proof_bfsag} 
        \Exp[\lVert \ppm_{k+1} - \ppm^*\rVert^2\lvert \ppm_{k}] & \leq (1 - 2 \eta \mu + \eta^2 L^2) \lVert \ppm_{k} - \ppm^* \rVert ^2,\\
      \Exp[\lVert \ppm_{k+1} - \ppm^*\rVert^2\lvert \ppm_{0}]  &\leq (1 -  \mu^2/L^2)^k \lVert \ppm_{0} - \ppm^* \rVert ^2,
    \end{split}
\end{equation}
where $\mu = \min \left\{p_l(1-p_l-p_h)^{(k-j)}\frac{\mu_\mathrm{low}\lVert \ppm_j - \ppm^* \rVert^2}{\lVert \ppm_k - \ppm^* \rVert^2},p_h(1-p_l-p_h)^{(k-j)}\frac{\mu_\mathrm{high}\lVert \ppm_j - \ppm^* \rVert^2}{\lVert \ppm_k - \ppm^* \rVert^2}\right\}$ for $j=1,\dots,k$; and learning rate is chosen as $\eta = \mu/L^2$ subject to $\mu^2/L^2\leq 1$. This completes the proof of Theorem \ref{theorem1}.

The constants $\mu$ and $L^2$ in \eqref{eq:proof_bfsag} are affected by the parameter update history as mentioned in Section \ref{sec:bfsag}. To see this, let us define
\begin{equation}
\begin{split}
    c^k_{\min} &= \underset{j}{\min} \left\{ (1-p_l-p_h)^{k-j} \lVert \thetaa_j-\thetaa^* \rVert^2 \right\};\\
    c^k_{\max} &= \underset{j}{\max} \left\{ (1-p_l-p_h)^{k-j} \lVert \thetaa_j-\thetaa^* \rVert^2 \right\}; \quad \text{for $j=1,\dots,k-1$}.\\
\end{split}
\end{equation}
Hence, the constants $\mu$ and $L^2$ can be written as
\begin{equation}
    \begin{split}
    \mu &= \frac{c^k_{\min}}{\lVert \thetaa_k-\thetaa^*\rVert^2} \min \left\{ p_l\mu_\mathrm{low},p_h\mu_\mathrm{high} \right\}; \\
        L^2 &= \frac{c^k_{\max}}{\lVert \thetaa_k-\thetaa^*\rVert^2} \max \left\{ p_lL^2_\mathrm{low},p_hL^2_\mathrm{high} \right\}. \\
    \end{split}
\end{equation}
Note that, if $p_l$ and $p_h$ are fixed $\mu$ depends on $c^k_{\min}$, \ie ~on $\min \left\{ (1-p_l-p_h)^{k-j} \lVert \thetaa_j-\thetaa^* \rVert^2 \right\}$ for $j=1,\dots,k$. Further, $(1-p_l-p_h)^{k-j}$ increases with $j$ since $1-p_l-p_h<1$ but $\lVert \thetaa_j-\thetaa^* \rVert^2$ depends on the parameter updates $\left\{\ppm_j\right\}_{j=1}^k$. Similarly, $L^2$ depends on $c^k_{\max}$ and in turn on $\left\{\ppm_j\right\}_{j=1}^k$. Hence, the parameter update history affects $\mu$ and $L^2$.

\section{Proof of Theorem \ref{theorem2}}\label{proof_thm2} 
Using the assumption of strong convexity of objective using high-fidelity models, 
\begin{equation}
\begin{split}
        (\ppm_k - \ppm^*)^T\hb_\mathrm{high}(\ppm_k) &\ge \mu_\mathrm{high} \lVert \ppm_k - \ppm^* \rVert^2,\\
\end{split}
\end{equation}
where $\mu_\mathrm{high}$ is a constant.
Similarly, if the high-fidelity gradients are Lipschtiz continuous 
\begin{equation}
\begin{split}
        \lVert \hb_\mathrm{high}(\ppm_k) \rVert ^2 & \leq L^2_\mathrm{high} \lVert \ppm_k - \ppm^* \rVert^2,\\
\end{split}
\end{equation}
where $L_\mathrm{high}$ is the Lipschitz constant.
For the inner iterations, we can evaluate the following expectation
\begin{equation}
    \begin{split}
        \Exp[\lVert & \ppm_{k+1} - \ppm^*\rVert^2\lvert \ppm_{k},\ppm_\mathrm{prev}]  \\
        & = \Exp\left[\left\lVert \ppm_{k} - \ppm^* -\eta \left[\widehat\hb_\mathrm{high}-\frac{\boldsymbol\alpha}{N_h}\sum_{i=1}^{N_h}\left(\hb_\mathrm{low}(\thetaa_\mathrm{prev};\Ym_{i})-\widehat\hb_\mathrm{low}\right)\right]\right\rVert^2\Bigg\lvert \ppm_{k},\ppm_\mathrm{prev}\right]\\
        & = \lVert \ppm_{k} - \ppm^* \rVert ^2 - 2 \eta (\ppm_k-\ppm^*)^T\hb_\mathrm{high}(\ppm_k)+\eta^2\lVert \hb_\mathrm{high}(\ppm_k)\rVert^2 \\
        & \qquad \qquad + \eta^2 \sum_{q=1}^{n_{\thetaa}}\Var\left( \widehat h_{\mathrm{high},q}-\frac{{\alpha}_{qq}}{N_h}\sum_{i=1}^{N_h}\left(h_{\mathrm{low},q}(\thetaa_\mathrm{prev};\Ym_{i})-\widehat h_{\mathrm{low},q}\right)\Bigg\lvert \ppm_k,\ppm_\mathrm{prev} \right),\\
    \end{split}
\end{equation}
where ${h}_q$ is the gradient with respect to $\theta_q$ and $\mathrm{Var}(\cdot)$ denotes variance of its argument.
Note that, if $\widehat{\hb}_\mathrm{low}=\Exp[\hb_\mathrm{low}(\ppm;\Ym)]$ exactly,
\begin{equation}
    \begin{split}
        \Var\Bigg( \widehat h_{\mathrm{high},q}-\frac{\alpha_{qq}}{N_h}\sum_{i=1}^{N_h}\big(h_{\mathrm{low},q}&(\thetaa_\mathrm{prev};\Ym_{i})-\widehat h_{\mathrm{low},q}\big) \Bigg\lvert \ppm_k, \ppm_\mathrm{prev}\Bigg)\\
        & = \frac{1}{N_h} (1-\rho_{hl,q}^2)\Var(h_{\mathrm{high},q}( \ppm_k;\Ym)), \\
    \end{split}
\end{equation}
where $\alpha_{qq}=\mathrm{Cov}(h_{\mathrm{low},q}(\ppm_\mathrm{prev};\Ym),h_{\mathrm{high},q}(\ppm_k;\Ym))/\Var(h_{\mathrm{low},q}(\ppm_\mathrm{prev};\Ym))$ and the correlation coefficient $\rho_{hl,q}=\mathrm{Cov}(h_{\mathrm{low},q}(\ppm_\mathrm{prev};\Ym),h_{\mathrm{high},q}(\ppm_k;\Ym))/\sqrt{\Var(h_{\mathrm{low},q}(\ppm_\mathrm{prev};\Ym))\Var(h_{\mathrm{high},q}(\ppm_{k};\Ym))}$. 
On the other hand, if we use $N_l$ samples to estimate $\widehat{\hb}_\mathrm{low}$, \ie~$\widehat{\hb}_\mathrm{low}=\frac{1}{N_l}\sum_{i=1}^{N_l}\hb_\mathrm{low}(\thetaa_\mathrm{prev};\Ym_i)$ then we can write
\begin{equation}
    \begin{split}
        \Var\Bigg( \widehat h_{\mathrm{high},q}-\frac{\alpha_{qq}}{N_h}\sum_{i=1}^{N_h}\big(h_{\mathrm{low},q}&(\thetaa_\mathrm{prev};\Ym_{i})-\widehat h_{\mathrm{low},q}\big) \Bigg\lvert \ppm_k, \ppm_\mathrm{prev}\Bigg)\\
        & = \frac{1}{N_h} \left(1-\frac{\rho_{hl,q}^2}{1+N_h/N_l}\right)\Var(h_{\mathrm{high},q}( \ppm_k;\Ym)), \\
    \end{split}
\end{equation}
where the coefficient $\alpha_{qq}$ is obtained by minimizing the mean-square error in $\widehat{h}_{\mathrm{high},q}$ \cite{pasupathy2012control,fairbanks2017low}, \textit{i.e.},
\begin{equation}
\alpha_{qq}=\frac{\mathrm{Cov}(h_{\mathrm{low},q}(\ppm_\mathrm{prev};\Ym),h_{\mathrm{high},q}(\ppm_k;\Ym))}{\Var(h_{\mathrm{low},q}(\ppm_\mathrm{prev};\Ym))}\left( \frac{1}{1+N_h/N_l} \right),
\end{equation}
 and the correlation coefficient $\rho_{hl,q}$ is same as before.
Next, let us assume 
\begin{equation}
\frac{1}{N_h} \left(1-\frac{\rho_{hl,q}^2}{1+N_h/N_l}\right)\Var(h_{\mathrm{high},q}(\ppm_k;\Ym)) \leq L^2_\mathrm{high} \delta_{k,q} \lVert \ppm_{k} - \ppm^* \rVert ^2 \end{equation}
for some constants $\delta_{k,q}$. 
Further, assume $\delta_k = \max \{1,\delta_{k,q} \}$ for $q=1,\dots,n_{\thetaa}$.
Hence,
\begin{equation}
    \begin{split}
        \Exp[\lVert \ppm_{k+1} - \ppm^*\rVert^2\lvert \ppm_{k},\ppm_\mathrm{prev}] & \leq (1-2\eta \mu_\mathrm{high}+ \eta^2 L^2_\mathrm{high}) \lVert \ppm_{k} - \ppm^* \rVert ^2 \\
        & \qquad \qquad \qquad \qquad + \eta^2 L^2_\mathrm{high}\delta_k \lVert \ppm_{k} - \ppm^* \rVert ^2\\
        & \leq (1-2\eta \mu_\mathrm{high}+ 2 \eta^2 L^2_\mathrm{high}\delta_k) \lVert \ppm_{k} - \ppm^* \rVert ^2
    \end{split}
\end{equation}
At $k$th inner iteration let us use the learning rate $\eta = \frac{\mu_\mathrm{high}}{2L^2_\mathrm{high}\delta_k}$. This leads to
\begin{equation}
    \begin{split}
        \Exp[\lVert \ppm_{k+1} - \ppm^*\rVert^2\lvert \ppm_{k},\ppm_\mathrm{prev}] 
        & \leq \left(1-\frac{\mu^2_\mathrm{high}}{2L^2_\mathrm{high}\delta_k}\right) \lVert \ppm_{k} - \ppm^* \rVert ^2,\\
      \Exp[\lVert \ppm_{k+1} - \ppm^*\rVert^2\lvert \ppm_\mathrm{prev}]  & \leq \left(1-\frac{\mu^2_\mathrm{high}}{2L^2_\mathrm{high}\underline{\delta}}\right)^k \lVert \ppm_\mathrm{prev} - \ppm^* \rVert ^2,
    \end{split}
\end{equation}
where $\underline{\delta}=\min \{\delta_i\}_{i=1}^k$; and $\ppm_1=\ppm_\mathrm{prev}$. 
Similarly, for $j$th outer iteration,
\begin{equation}\label{eq:proof_end}
    \begin{split}
        \Exp[\lVert \widetilde\ppm_{j} - \ppm^*\rVert^2\lvert \ppm_0]  & \leq\left(1-\frac{\mu^2_\mathrm{high}}{2L^2_\mathrm{high}\delta}\right)^{jm} \lVert \ppm_{0} - \ppm^* \rVert ^2,
    \end{split}
\end{equation}
where $\delta = \min \{\underline \delta_i\}_{i=1}^j$ subjected to $\frac{\mu^2_\mathrm{high}}{2L^2_\mathrm{high}\delta}\leq 1$ and this proves Theorem \ref{theorem2}. 

Note that, if $\rho_{hl}$ is close to 1, \ie\  the low- and the high-fidelity models are highly correlated, then $\delta_k$  can be assumed small. This implies that $\delta$ will be close to 1 and, thus, we can use a larger learning rate $\eta$. This, in turn, leads to a smaller right hand in \eqref{eq:proof_end} and a tighter bound on $\Exp[\lVert \widetilde\ppm_{j} - \ppm^*\rVert^2\lvert \ppm_0]$.
\section*{References}
\biboptions{sort&compress}
\bibliography{mybibfile}

\end{document}